\documentclass[12pt]{article}%

\usepackage[margin=1.2in]{geometry}
\usepackage{graphicx,todonotes}
\usepackage{graphics}

\usepackage{color}
\usepackage{amsmath}
\usepackage{amsfonts}
\usepackage{amssymb}
\usepackage{enumerate}
\usepackage[normalem]{ulem}
\usepackage{hyperref}

\definecolor{PineGreen}{cmyk}{0.9,0,1,0.40}

\newtheorem{theorem}{Theorem}[section]

\newtheorem{definition}[theorem]{Definition}

\newtheorem{example}[theorem]{Example}

\newtheorem{lemma}[theorem]{Lemma}
\newtheorem{assumption}[theorem]{Assumption}
\newtheorem{notation}[theorem]{Notation}

\newtheorem{proposition}[theorem]{Proposition}

\newtheorem{remark}[theorem]{Remark}

\newenvironment{proof}[1][Proof]{\textbf{#1.} }{\ \rule{0.5em}{0.5em}}

\newcommand{\bfone}{{\rm \bf 1}}
\newcommand{\E}{{\rm \bf E}}

\newcommand{\Prob}{{\rm \bf P}}
\newcommand{\prob}{{\rm \bf P}}
\newcommand{\LCP}{{\rm LCP}}

\newcommand{\conv}{{\rm conv\ }}

\newcommand{\range}{{\rm range}}
\newcommand{\dN}{{\mathbb N}}

\newcommand{\dR}{{\mathbb R}}
\newcommand{\dZ}{{\mathbb Z}}

\newcommand{\calA}{{\cal A}}
\newcommand{\calB}{{\cal B}}

\newcommand{\calD}{{\cal D}}

\newcommand{\calF}{{\cal F}}
\newcommand{\calG}{{\cal G}}
\newcommand{\calH}{{\cal H}}

\newcommand{\calM}{{\cal M}}

\newcommand{\calQ}{{\cal Q}}

\newcommand{\calT}{{\cal T}}

\newcommand{\ep}{\varepsilon}

\newcommand{\Leb}{{\rm Leb}}

\newcounter{figurecounter}
\title{Multiplayer Stopping Games with a Conditionally Atomless Filtration%
\thanks{Solan acknowledges the support of the Israel Science Foundation, Grant \#211/22.}}
\author{Robert Samuel Simon%
\thanks{Department of Mathematics, London School of Economics
Houghton Street, London WC2A 2AE, UK.
e-mail: R.S.Simon@lse.ac.uk.}, 
Eilon Solan%
\thanks{School of Mathematical Sciences, 
Tel Aviv University, 6997800, Israel.
e-mail: eilonsolan@gmail.com.}, 
and Grzegorz Tomkowicz\thanks{Centrum Edukacji G2, ul. Moniuszki 9, 41-902 Bytom, Poland.
e-mail: gtomko@vp.pl.}}

\date{\today}

\begin{document}

\maketitle

\begin{abstract}
Stopping games are the multiplayer generalization of optimal stopping problems,
in which each player independently chooses when to stop,
and the payoff depends on who stops first and when.
While the existence of $\varepsilon$-equilibria in two-player stopping games
with integrable payoffs is known,
the multiplayer case has remained open.

We prove that every multiplayer stopping game in discrete time
with integrable payoff process admits an $\varepsilon$-equilibrium,
for every $\varepsilon > 0$,
provided the filtration is conditionally atomless.
This condition is satisfied, in particular, when the underlying probability space
is $\{0,1\}^{\dN}$ endowed with the uniform measure and its natural filtration.

Our method also shows that every three-player stopping game,
without the restriction to conditionally atomless filtration, 
admits an $\varepsilon$-equilibrium,
for every $\varepsilon > 0$.

The proof combines techniques developed for two-player stopping games and for
multiplayer quitting games,
with a key new insight:
the richness encoded in a conditionally atomless filtration can substitute for
an external public randomization device,
allowing players to coordinate without an explicit coordination signal.
\end{abstract}

\bigskip
\noindent\textbf{Keywords:}
Stopping games, equilibrium.

\bigskip
\noindent\textbf{MSC2020 Classification:}
Primary: 91A55, 91A25. 
Secondary: 60G40, 91A06.

\newpage
\section{Introduction}

Stopping games are the multiplayer analogue of optimal stopping problems.
An $n$-player stopping game is played over a filtered probability space.
Each player chooses a stopping time,
and the payoff depends on the first stage at which at least one player stops and on the set of players who stop at that stage.

Two-player zero-sum stopping games in discrete time 
in which at each stage only one player is allowed
to stop were introduced by \cite{Dynkin69}, who proved the existence of the value and identified (pure) $\ep$-optimal strategies.
\cite{Neveu75} extended the model to allow simultaneous stopping and provided sufficient conditions for the existence of the value.
\cite{RSV01} proved that every two-player zero-sum stopping game with integrable payoff process admits a value in randomized strategies, and \cite{SS04} showed that every two-player non-zero-sum stopping game with integrable payoff process admits an approximate equilibrium in randomized strategies.
Beyond two players, the existence problem
has remained open without additional assumptions
on the game.
\cite{HHM24} established the existence of 
$\ep$-equilibria in pure stopping times under
the condition that at every stage, each player weakly prefers that any other player stops alone, rather than that the two stop together.
\cite{Heller12}
establishes the existence of sequential normal-form
correlated approximate equilibria.

In this paper, we prove that every discrete-time multiplayer stopping game with integrable payoff process admits an $\ep$-equilibrium in randomized strategies whenever the filtration is conditionally atomless.
Following \cite{Delbaen2020}, 
recall that a sigma-algebra $\calF$ is \emph{conditionally atomless relative to} a sub-sigma-algebra $\calG$ if every $A \in \calF$ has a subset $B \in \calF$ such that
$0 < \E[\mathbf{1}_B \mid \calG] < \E[\mathbf{1}_A \mid \calG]$,
almost surely on the set $\{\E[\mathbf{1}_A \mid \calG] > 0\}$.
A filtration $(\calF_t)_{t \in \dN}$ is \emph{conditionally atomless} if $\vee_{t \in \dN}\calF_t$ is conditionally atomless relative to $\calF_k$ for every $k \in \dN$.
The condition concerns the information revealed
over the entire future, rather than the information
revealed at a single stage.
In particular, every multiplayer stopping game with integrable payoff process defined on the probability space $\{0,1\}^\dN$, endowed with the uniform distribution and the natural filtration, admits an $\ep$-equilibrium in randomized strategies for every $\ep > 0$.

The tools that we develop will allow us to deduce that every three-player stopping game admits an $\ep$-equilibrium for every $\ep > 0$,
whether or not the filtration is conditionally atomless.

Thus, unlike earlier results, 
we do not impose ordering conditions on the payoffs,
nor do we introduce external correlation.

Our proof combines ideas from \cite{SS04} and \cite{SS20}, together with a new insight.
\cite{SS20} showed that every multiplayer quitting game---that is, a stopping game with trivial filtration and time-homogeneous payoff process---admits a sunspot $\ep$-equilibrium for every $\ep > 0$.
A sunspot $\ep$-equilibrium is an $\ep$-equilibrium in an extended game in which, at the beginning of every stage, the players observe a public signal uniformly distributed on $[0,1]$, independently of
past signals and of the players' private randomizations.

The key new ingredient in our proof is that conditional atomlessness of the filtration can replace the public randomization device used in sunspot equilibria.
More precisely, the public lotteries
can be approximated by events observable after
suitable finite blocks of stages. We use this
randomness from the original filtration to implement
the sunspot construction within the stopping-game
method of \cite{SS04}, without adding an external
public randomization device.

This result should be contrasted with the continuous-time case.
While two-player non-zero-sum stopping games in continuous time with integrable and right-continuous payoff process and right-continuous filtration admit an approximate equilibrium in randomized strategies (\cite{LS13}),
three-player stopping games in continuous time with integrable payoff process may fail to admit an approximate equilibrium, even when the payoff process is time-homogeneous.
\cite{LSV05} provided such an example with a trivial
filtration. Their nonexistence argument also applies
when the players observe an arbitrary common
filtration.

The paper is organized as follows.
Section~\ref{sec:model:result} presents the model and the main results.
Section~\ref{sec:preliminaries} introduces the techniques used in the proof.
The proof of the main result appears in Section~\ref{sec:proof},
and in Section~\ref{sec:proof:absorbing-eq} we establish an intermediate result used in that proof.
In Section~\ref{section:three} we prove that every three-player stopping game admits an approximate equilibrium. 
Section~\ref{sec:discussion} concludes.

%%%%%%%%%%%%%%%%%
\section{The Model and the Main Results}
\label{sec:model:result}

For every vector $x = (x^i)_{i \in I}$
and every $i \in I$
we denote $x^{-i} := (x^j)_{j \neq i}$.
For every finite set $A$,
denote by $\Delta(A)$ the set of probability distributions on $A$.

%%%%%%%%%%%%%%%%%
\subsection{Stopping Games}
\label{sec:model}

We consider a discrete-time multiplayer stopping game with set of stages $\dN = \{0,1,2,\dots\}$.

\begin{definition}[Stopping game]
\label{def:stopping:game}
A \emph{stopping game} is a 6-tuple $\Gamma =
(\Omega, I, \calF,\prob,{\mathbb F},X)$ where:
\begin{itemize}
\item $(\Omega,\calF,\prob)$ is a probability space.
\item $I = \{1,\dots,n\}$ is a finite set of $n$ players with $n \geq 2$.
\item ${\mathbb F} = (\calF_t)_{t \in \dN}$ is a filtration over $(\Omega,\calF,\prob)$.
\item $X = (X_t)_{t \in \dN}$ is an ${\mathbb F}$-adapted $\mathbf{R}^{n(2^{n}-1)}$-valued payoff process,
where $X_t = (X_{J,t}^i)_{\varnothing \neq J \subseteq I,i \in I}$.
We will require that the payoff process is integrable:
$\E\left[\sup_{t \in \dN, \varnothing \neq J \subseteq I, i \in I}|X_{J,t}^i|\right] < \infty$.
\end{itemize}
\end{definition}

The interpretation of the data of the game is that $\calF_t$ is the information the players possess at stage $t$,
and $X^i_{J,t}$ is player~$i$'s terminal payoff 
when no player stopped before stage $t$,
and the players in $J$ are those who stop at stage $t$.

\subsection{Strategies, payoffs, and equilibria}

We allow players to randomize and
use the concept of behavioral strategies.

\begin{definition}[Strategy]
\label{def:strategy}
Let $\Gamma =
(\Omega, I, \calF,\prob,{\mathbb F},X)$ be a stopping game.
A \emph{(behavioral) strategy} for player $i \in I$ is an $\mathbb{F}$-adapted
process $x^i=(x^i_t)_{t \in \dN}$ with values in
$[0,1]$.
\end{definition}

At stage $t$, if the game has not stopped earlier, player $i$ stops with probability $x^i_t$.
The fact that $x^i_t$ is $\calF_t$-measurable captures the property that this probability is determined by the information that player~$i$ possesses at stage $t$.

To realize behavioral strategies on the underlying probability space, 
we require that $(\Omega,\calF,\prob)$ is rich enough to support a sequence $(\zeta^i_{t})_{t \in \dN, i \in I}$ of i.i.d.~random variables, uniformly
distributed over $[0,1]$,
independent of $\vee_{t \in \dN} \calF_t$. 
If $(\Omega,\calF,\prob)$ does not support such a sequence,
extend the probability space to 
$(\Omega \times [0,1]^{I \times \dN}, \calF \times \calB^{I \times \dN},
\prob \otimes \Leb^{I \times \dN})$,
where $\calB$ is the Borel sigma algebra on $[0,1]$,
$\Leb$ is the uniform measure on $[0,1]$,
and $\prob \otimes \Leb^{I \times \dN}$ is the independent product of $\prob$ with countably many copies of $\Leb$.
These variables are only a technical device for implementing randomized stopping decisions.
%\footnote{In the introduction, 
%we mentioned stopping games that are defined on the probability space $\{0,1\}^\dN$, endowed with the uniform distribution and the natural filtration.
%Since we incorporate the randomization of the players within the probability space, the probability space of this example should have been $\{0,1\}^\dN \times [0,1]^{I \times \dN}$.
%Since we did not discuss the implementation of behavioral strategies in the introduction, we elected to keep the discussion there at an intuitive level.}
Intuitively, under the strategy $x^i = (x^i_t)_{t \in \dN}$, 
player~$i$ stops at the first stage $t$ at which $\zeta^i_t \leq x^i_t$.

\begin{definition}[$t^i$, $\theta$, $J_\theta$]
\label{def:theta}
Let $\Gamma$ be a stopping game,
and let $x = (x^i)_{i \in I}$ be a strategy profile.
We denote by
\[ t^i = t^i(x^i) := \min\{t \in \dN \colon \zeta^i_{t} \leq x^i_{t}\} \]
the \emph{stopping stage} of player~$i$ under $x^i$,
with the convention that the minimum of an empty set is $\infty$.
The first stage at which at least one player stops is
\[ \theta = \theta(x) := \min\{ t^i \colon i \in I\}. \]
On $\{\theta < \infty\}$, the set of players
stopping at stage $\theta$ is
\[ J_\theta = J_\theta(x) := \{i \in I \colon \theta = t^i\} \subseteq I. \] 
\end{definition}

Since $(\zeta^i_t)_{t\in \dN,i \in I}$ are independent and independent of $\vee_{t \in \dN} \calF_t$,
for every stage $t$, conditional on $\calF_t$ and on the event that no player has stopped before stage $t$, the stopping decisions of the players at stage $t$ are independent.

The stopping time $t^i$ of each player $i$ depends on her strategy $x^i$.
When calculating the probability and the expectation of events that depend on these stopping times,
we will write 
$\Prob_x$ and $\E_x$ for $\prob$ and $\E$,
to emphasize the strategy profile that the players follow.

\begin{definition}[Expected payoff, $\gamma^i(x)$]
\label{def:payoff}
The \emph{expected payoff} to player $i$ under strategy profile $x = (x^i)_{i \in I}$ is
\[
   \gamma^i(x)
   := \E_x\bigl[X_{J_\theta,\theta}^i \cdot \bfone_{\{\theta<\infty\}}\bigr].
\]
\end{definition}

Note that $\gamma^i(x)$ is well-defined
because the payoff process is integrable.
Throughout the paper we will define various games $\widetilde \Gamma$ that differ from $\Gamma$ in their parameters.
We will denote by $\widetilde\gamma^i(x)$ player~$i$'s expected payoff in these games.

\begin{remark}[Final payoff]
The convention that payoffs equal $0$ when no player ever stops is
without loss of generality.
Any bounded final payoff $\chi = (\chi_i)_{i \in I}$
measurable with respect to $\vee_{t \in \dN} \calF_t$ (payoff when $\theta=\infty$) can be
reduced to the zero case by subtracting the conditional expectation
$\E[\chi^i \mid \calF_t]$ from each $X^{i}_{J,t}$
for every nonempty subset $J \subseteq I$,
every $i \in I$, and every $t \in \dN$; see \cite[Section~3.4]{LS05} for details.
\end{remark}

\begin{definition}[$\ep$-equilibrium]
Given $\ep\geq 0$, a strategy profile $x^*$ is an \emph{$\ep$-equilibrium}
if for every $i \in I$ and every strategy $x^i$ of player~$i$,
\[
   \gamma^i(x^*) \geq \gamma^i(x^i,x^{*,-i}) - \ep.
\]
\end{definition}

\subsection{Conditionally Atomless Filtrations}

The concept of a conditionally atomless filtration,
which we define now,
will play an important role in the paper.
A sigma-algebra $\calF$ is \emph{conditionally atomless} relative to a sub-sigma-algebra $\calG$ if for every $A \in \calF$ there is $B \subset A$ such that $B \in \calF$ 
and 
$0 < \E[\mathbf{1}_B \mid \calG] < \E[\mathbf{1}_A \mid \calG]$
on the set $\{\E[\mathbf{1}_A \mid \calG] > 0\}$.
This concept was defined and studied by
\cite{Delbaen2020}.

We now extend this concept to filtrations.

\begin{definition}[Conditionally atomless filtration]
Let $(\Omega,\calF,\prob)$ be a probability space,
and let $\mathbb{F} = (\calF_t)_{t \in \dN}$ be a filtration.
We say that $\mathbb{F}$ is \emph{conditionally atomless} if
$\vee_{t \in \dN} \calF_t$ is conditionally atomless relative to $\calF_k$, 
for every $k \in \dN$.
\end{definition}

\begin{example}[Conditionally atomless filtrations]
\label{example:cond:atom}
Let $(\Omega,\calF,\prob)$ be the probability space that corresponds to countably many tosses of a fair coin:
$\Omega = \{0,1\}^{\mathbf{N}}$,
$\calF$ is the Borel sigma-algebra generated by the product topology,
and $\prob$ is the product of countably many copies of a uniform distribution on $\{0,1\}$.
For every $t \in \dN$,
let $\calF_t$ be the sigma-algebra generated by the cylinder sets defined by the first $(t+1)$ coordinates.
Then $(\calF_t)_{t \in \dN}$ is conditionally atomless.

Let $\prob'$ be 
the product measure corresponding to countably many tosses of an unfair coin with success probability $1/3$.
The filtration $(\calF_t)_{t \in \dN}$ is conditionally atomless
for $(\Omega,\calF,\prob')$ as well.
\end{example}

\begin{example}[Non-conditionally atomless filtrations]
Let $(\Omega,\calF,(\calF_t)_{t \in \dN})$ be as in Example~\ref{example:cond:atom}.
Let $\prob''$ be
the probability space that corresponds to the independent tossing of unfair coins: 
coin $t$ has parameter $\frac{1}{2^{t+1}}$.
Since $\sum_{t \in \dN} \frac{1}{2^{t+1}} < \infty$,
$\prob''$ contains an atom at $(0,0,\dots)$, and hence
$(\calF_t)_{t \in \dN}$ is not conditionally atomless for $(\Omega,\calF,\prob'')$.

Note that if $\Omega$ is countable,
then no filtration that is defined on $\Omega$ 
is conditionally atomless.
\end{example}

\cite[Lemma~3]{Delbaen2020} implies that if $\calF$ is conditionally atomless relative to $\calG$,
then there is a uniform random variable on $[0,1]$ that is independent of $\calG$.
We therefore have the following.
\begin{lemma}[Conditional atomlessness and uniform random variable]
\label{lemma:dalbaen}
Let $(\Omega,\calF,\prob)$ be a probability space,
and let $(\calF_t)_{t \in \dN}$ be a conditionally atomless filtration.
Then, for every $k \in \dN$
there is a random variable 
$U_k \sim \mathrm{Uniform}[0,1]$ that is 
measurable with respect to $\vee_{t \in \dN} \calF_t$ and
independent of $\calF_k$.
\end{lemma}

The conclusion of Lemma~\ref{lemma:dalbaen} captures the idea that the players' future observations
are ``rich enough'' to simulate any auxiliary randomization device at
every stage:
Under conditional atomlessness, the uniform signals required by the sunspot equilibrium
construction of~\cite{SS20} can be synthesized from the filtration
itself.

%For the two probability spaces presented in Example~\ref{example:cond:atom},
%the standard filtration $\mathbb{F} = (\calF_t)_{t \in \dN}$,
%where $\calF_t$ is the algebra generated by the cylinder sets determined by the first $t+1$ coordinates,
%is conditionally atomless.

%Note that in the context of stopping games,
%since 
%$\vee_{t \in \dN} \calF_t$ 
%is independent of $(\zeta^i_t)_{i \in I, t \in \dN}$,
%the random variable $U_k$ in Lemma~\ref{lemma:dalbaen}
%is also independent of $(\zeta^i_t)_{i \in I, t \in \dN}$.
%and therefore also of $(t^i)_{i\in I}$, $\theta$, and $J_\theta$.

\subsection{Main Results}

Our main result is the following.

\begin{theorem}[Existence of $\ep$-equilibrium]
\label{thm:main}
Every multiplayer stopping game with 
a conditionally atomless filtration 
(and an integrable payoff process) 
admits an $\varepsilon$-equilibrium,
for every $\varepsilon > 0$.
\end{theorem}

The condition that the payoff process is integrable ensures that the expected payoff is well-defined.
Shmaya and Solan \cite{SS04} proved Theorem~\ref{thm:main} when $|I|=2$,
without the condition that the filtration is conditionally atomless.
The tools we will develop to solve Theorem~\ref{thm:main} will allow us to solve the case of stopping games with three players and with any filtration (not necessarily conditionally atomless).

\begin{theorem}
\label{theorem:three}
Every three-player stopping game 
(with an integrable payoff process) 
admits an $\varepsilon$-equilibrium,
for every $\varepsilon > 0$.
\end{theorem}

We do not know whether Theorem~\ref{theorem:three} holds when $|I| > 3$.
The difficulty in proving Theorems~\ref{thm:main} and~\ref{theorem:three} stems from the fact that the horizon is infinite.
When the horizon is finite,
an $\ep$-equilibrium exists by a standard backward induction argument.

\subsection{Sketch of the Proof of Theorem~\ref{thm:main}}

We now provide the main ideas of the proofs.

\paragraph{Reductions.}
We first make a few simplifying assumptions on the game, which are without loss of generality (see Section~\ref{sec:reduction}).
Their proof follows standard approximation arguments.
\begin{itemize}
\item 
The payoff process is bounded (Section~\ref{sec:reduction:bounded}) has a finite range (Section~\ref{section:finite:range}),
and all payoffs are an integer multiple of a certain real number (Section~\ref{section:reduction:unique}).
\item 
Payoffs are generic, in the sense that if for some $g = (g^i_J)_{i \in I, \varnothing \subset J \subseteq I}$ and $g' = (g'^i_J)_{i \in I, \varnothing \subset J \subseteq I}$ 
in the range of $X$ we have 
$g^i_J = g'^i_{J'}$ for some $i \in I$ and $ \varnothing \subset J,J' \subseteq I$,
then $g=g'$ and $J = J'$ (Section~\ref{section:reduction:unique}).
\item
Any payoff vector in the range of the payoff process,  
almost surely appears infinitely often along the play (Section~\ref{section:ii:payoff}).
\end{itemize}

We will denote by 
$y^i := \max_{g} g^i_{\{i\}}$
the maximal payoff that player~$i$ can achieve by stopping alone, where $g$ ranges over all payoff vectors in the range of $X$.
We will also denote by $Z_i \in \dR^n$ the payoff vector when player $i$ stops and obtains $y^i$.
By the reductions,
this maximum $y^i$ is attained, the vector $Z_i$ is unique, and both are independent of $\omega$.
We will also denote by $R$ the $n \times n$ matrix  
whose $i$'th column is $Z_i$.

In our proof, the construction of an $\ep$-equilibrium will depend on $R$,
and on whether this matrix is close to a Q-matrix.
Q-matrices are related to linear complementarity problems (see Section~\ref{subsec:lcp}),
and their relation to stopping games has been studied in \cite{SS20}.

\paragraph{A certain principal submatrix of $R$ is close to a Q-matrix.}

When a certain principal submatrix of $R$ is close to a Q-matrix,
we follow the construction of \cite{SS20}.
In this construction, 
at every stage of the $\ep$-equilibrium, at most one player stops with a small probability,
while the remaining players stop with probability 0.
In the model of \cite{SS20},
it is assumed that players observe at the beginning of every stage a public signal that is uniformly distributed in $[0,1]$,
and both the identity of the player who stops with positive probability as well as the probability by which she stops
are determined by the sequence of past signals.
In our model,
the information available in the filtration will replace the public signals.

\paragraph{A certain principal submatrix of $R$ is not close to any Q-matrix.}

When the relevant principal submatrix of
$R$ is not close to any Q-matrix, we adapt the construction of \cite{SS04}.
Roughly, as in \cite{SS04}, 
we will construct an increasing sequence of bounded $\mathbb F$-adapted stopping times $\phi_0 < \phi_1 < \dots$
such that in the \emph{periodic} stopping game that starts at stage $\phi_k$ and restarts when reaching stage $\phi_{k+1}$ there is a $O(\sqrt\ep)$-equilibrium that terminates with positive probability.
Moreover, the payoff in these $O(\sqrt\ep)$-equilibria,
which are $\calF_{\phi_k}$-measurable functions, does not vary much with $k$.
A concatenation of these $O(\sqrt\ep)$-equilibria,
possibly with an addition of punishment strategies, 
will yield an $O(\sqrt\ep)$-equilibrium in the game that starts at stage $\phi_0$.
Since $\phi_0$ is bounded,
a standard backward induction argument will produce an $O(\sqrt\ep)$-equilibrium in the original game.

\paragraph{The proof of Theorem~\ref{theorem:three}.}
The only place in the proof of Theorem~\ref{thm:main} where the conditional atomlessness of the filtration takes part is when the principal submatrix of $R$ that we consider is close to a Q-matrix.
In the presence of three players, 
one can circumvent the use of the filtration 
by constructing a periodic $\ep$-equilibrium as in \cite{FTV97}, see also \cite{Solan03}.

%%%%%%%%%%%%%%%%%%%%%%%%%%%%%%%%%%%%%%%%%%%%%%%
\section{Preliminaries}
%%%%%%%%%%%%%%%%%%%%%%%%%%%%%%%%%%%%%%%%%%%%%%%
\label{sec:preliminaries}

This section presents the tools used in the proof of Theorems~\ref{thm:main} and~\ref{theorem:three}.
Section~\ref{sec:atomless2} shows how conditional atomlessness can be used as a source of randomness. Sections~\ref{subsec:lcp}--\ref{subsec:sunspot} introduce 
the necessary background in the linear complementarity, quitting games, and sunspot equilibria. Sections~\ref{sec:trees}--\ref{sec:nt:approx} develop stopping games on a finite tree and some approximation tools.

The notion of Q-matrices, which is presented in Section~\ref{subsec:lcp},
will be used to partition the probability space into two sets,
according to whether 
a certain principal submatrix of
$R$ is close to some Q-matrix.
Sections~\ref{sec:atomless2} and~\ref{subsec:sunspot} will be used to construct an approximate equilibrium when the relevant principal submatrix of 
$R$ is close to some Q-matrix,
and Sections~\ref{sec:trees}--\ref{sec:nt:approx} will be used to construct an approximate equilibrium when the relevant principal submatrix of 
$R$ is not close to any Q-matrix.
Section \ref{subsec:sunspot} is based on \cite{SS20}.
Sections~\ref{sec:trees}--\ref{sec:nt:approx} extend several definitions and results of \cite{SS04} from two-player stopping games to multiplayer stopping games.

We alert the reader that we will define several auxiliary games.
Besides the original stopping game $\Gamma$,
we will present \emph{quitting games},
which are stopping games with trivial filtration
and time-homogeneous payoff process.
These games will be denoted by $\calQ$
and used to construct approximate equilibria when a certain principal submatrix of 
$R$ is close to some Q-matrix.
We will also define stopping games on a finite tree,
which are a simple class of stochastic games,
and denote them by $\calT$.
These games will be used to define approximate equilibria when the relevant principal submatrix of $R$ is not close to any Q-matrix. 

%%%%%%%%%%%%%%%%%%%
\subsection{List of Notations}

To help the reader, we start with a list of notations, so that the reader can always easily find the meaning of each notation.

\bigskip
\noindent\textbf{Stopping Game $\Gamma =
(\Omega, I, \calF,\prob,{\mathbb F},X)$:}
\[\begin{tabular}{l|l|l}
Notation & Meaning & Page\\
\hline
$(\Omega,\calF,\prob)$ & Probability space & \pageref{def:stopping:game} \\
$I = \{1,\dots,n\}$ & Set of players & \pageref{def:stopping:game} \\
$\mathbb F = (\calF_t)_{t \in \dN}$ & Filtration & \pageref{def:stopping:game} \\
$X^i_{J,t}$ & Player~$i$'s payoff when set of players $J$ stop at stage $t$ & \pageref{def:stopping:game} \\
$x^i_t$ & Probability of stopping of player~$i$ at stage $t$ & \pageref{def:strategy}\\
$\theta$ & Stage in which the game stops & \pageref{def:theta}\\
$J_\theta$ & Set of players who stop at stage $\theta$ & \pageref{def:theta}\\
$\gamma^i(x)$ & Player~$i$'s expected payoff under strategy profile $x$ & \pageref{def:payoff}\\
$y^i$ & Player~$i$'s maximal payoff if she stops alone\\
& \qquad at some far future stage & \pageref{def:y:stopping}\\
$Z_{i}$ & The payoff vector when player~$i$ stops alone optimally & \\
& \qquad at some far future stage & \pageref{def:y:stopping}\\
$R$ & The matrix of single stoppings & \pageref{def:z:stopping}\\
$(\mathbb{F}, \tau, \tau')$-strategy &
Strategy between stages $\tau$ and $\tau'$ (excluded)
& \pageref{def:f:strategy}\\
$\bigl((\mathcal{G}_k)_{k \in \dN}, (q_{G,k})_{k \in \dN, G \in \calG_k}\bigr)$ & $\delta$-approximation of $\Gamma$ between $\tau$ and
$\tau'$ & \pageref{def:delta-approx}\\
\end{tabular}
\]

\noindent\textbf{Quitting Game $\calQ = (I,X)$:}
\[\begin{tabular}{l|l|l}
Notation & Meaning & Page\\
\hline
$I = \{1,\dots,n\}$ & Set of players & \pageref{def:stopping:game} \\
$\theta$ & Stage in which the game stops & \pageref{def:sunspot:theta}\\
$X^i_{J}$ & Player~$i$'s payoff when set of players $J$ stop at stage $\theta$ & \pageref{def:quitting}\\
$x^i$ & A stationary strategy for player $i$ & \pageref{def:quitting:strategy}\\
$0^i$ & The strategy of player~$i$ in which she never stops & \pageref{equ:45:a}\\
$1^i$ & The strategy of player~$i$ in which she always stops & \pageref{equ:45:a}\\
$\calQ_R$ & The quitting game that corresponds to the payoff matrix of\\
& \qquad single stoppings $R$ & \pageref{def:quitting:qr}\\
$R_{\bullet i}$ & The $i$'th column of the matrix $R$ & \pageref{def:quitting:ri}\\
$I^{\mathrm{sp}}_R$ & The set of supporters for $R$ & \pageref{def:player:nor}\\
$R^{\mathrm{sp}}_\calQ$ & The principal submatrix of $R_\calQ$ indexed by the set of supporters & \pageref{def:matrix:nor}\\
$\calQ_{\mathrm{sun}}$ & The game with a public correlation device & \pageref{def:sunspot}\\
$\xi^i$ & Player~$i$'s strategy in $\calQ_{\mathrm{sun}}$ & \pageref{def:sunspot:strategy}\\
$\gamma^i_\lambda(x)$ & Player~$i$'s expected $\lambda$-discounted payoff under strategy profile $x$ & \pageref{def:payoff:discounted}\\
$\gamma_{\mathrm{sun}}^i(\xi)$ & Player~$i$'s expected payoff in $\calQ_{\mathrm{sun}}$ under strategy profile $\xi$ & \pageref{def:sunspot:i:delta}\\
$(\widehat i_k)_{k \in \dN}, (\widehat \delta_k)_{k \in \dN}$ & The choices of players and stopping probabilities under the\\
& \qquad sunspot $\ep$-equilibrium in $\calQ_{\mathrm{sun}}$ & \pageref{def:sunspot:i:delta}\\ 
\end{tabular}
\]

\noindent\textbf{Stopping Games on finite trees $\calT = \bigl(I, S, S_1, r,
    (C_s)_{s \in S \setminus S_1},
    (p_s)_{s \in S \setminus S_1},
    (X_s)_{s \in S \setminus S_1}\bigr)$:}
\[\begin{tabular}{l|l|l}
Notation & Meaning & Page\\
\hline
$r$ & The root & \pageref{def:tree-game} \\
$S$ & The set of nodes & \pageref{def:tree-game} \\
$S_0$ & The set of non-leaf nodes & \pageref{def:tree-game} \\
$y^i_\calT$ & Player~$i$'s maximal payoff if she stops alone at some node &\pageref{equ:y:finite:tree}\\
$R_\calT$ & Matrix of payoff vectors under single stopping &\pageref{def:Z:tree}\\
$Z_{\calT,i}$ & The payoff vector when player~$i$ stops alone optimally &\pageref{def:Z:tree}\\
$\delta_\calT$ & The difference between $y^i_\calT$ and the second highest & \pageref{def:tree:delta}\\
& \qquad payoff player~$i$ can obtain by stopping alone\\
$\alpha^i$ & Stationary strategy of player~$i$ in $\calT$ & \pageref{subsection:tree:strategy}\\
$\theta_{\mathrm{fr}}$ & Last stage of the first round & \pageref{def:tree:theta}\\
$\pi(\alpha)$ & Probability of termination under $\alpha$ in the first round & \pageref{def:tree:theta}\\
$\rho(\alpha)$ & Expected payoff under $\alpha$ in the first round (unconditional) & \pageref{def:tree:theta}\\
$\gamma(\alpha)$ & Expected payoff under $\alpha$ in the first round (conditional) & \pageref{def:tree:theta}\\
$\prob_{s,\alpha}$ & Probability distribution induced by the initial node $s$ & \pageref{def:tree:prob}\\
& \qquad and strategy profile $\alpha$\\
$\mathbb{T}_\ep$ & Family of games that satisfy some simplifying properties & \pageref{def:family:t}\\
$\calT(\tau, \tau'; (\mathcal{G}_k), (q_{G,k}), G)$ & 
The game on a tree that corresponds to
a $\delta$-approximation & \pageref{def:tree:approx}\\
\end{tabular}
\]

\subsection{Simple Observations}

In our arguments we will repeatedly use 
the following property.
For every probability space $(\Omega,\calF,\prob)$,
every filtration $(\calF_t)_{t \in \dN}$,
every event $A \in \vee_{t \in \dN} \calF_t$,
and every $\ep > 0$,
there is $t_0 \in \dN$ sufficiently large and $B \in \calF_{t_0}$
such that $\prob(A \triangle B) \leq \ep$,
where $A \triangle B = (A \setminus B) \cup (B \setminus A)$.
A recursive argument then implies that for every finite $ \left(\vee_{t \in \dN} \calF_t\right)$-measurable partition $(A_i)_{i=1}^n$ of $\Omega$
there are $t_0 \in \dN$ sufficiently large and an
$\calF_{t_0}$-measurable partition $(B_i)_{i=1}^n$
such that $\prob(A_i \triangle B_i) \leq \ep$ for every $i \in \{1,\dots,n\}$.

We will also use three well-known observations regarding equilibria, which are summarized here for future reference.
The first observation states that if we change the payoff process by at most 1 on a set with measure $\ep$,
and by at most $\ep$ on the remaining set,
then any $\ep$-equilibrium in the game with the new payoff process is a $4\ep$-equilibrium in the original game.
The second observation states that 
for any bounded stopping time $\tau$,
any $\ep$-equilibrium in the game that starts at stage $\tau$ (that is, when the players are restricted to strategies that stop with probability 0 before stage $\tau$) can be extended to an $\ep$-equilibrium in the whole game (that is, without the restriction on the players' strategy sets).
The third observation states that if for some finite $\calF_0$-measurable partition $\calD = (D_1,\dots,D_L)$ of $\Omega$,
the stopping game restricted to each $D_\ell$ admits an $\ep$-equilibrium,
then so does the whole game.

For a stopping time $\tau$,
denote by $\calF_\tau$ 
the sigma-algebra generated by the sets
$\{ F \in \calF_k \colon k \in \dN, \tau = k \hbox{ on } F\}$.

\begin{lemma}
\label{lemma:simple}
let $\ep > 0$,  
and let $\Gamma = (\Omega,I,\calF,\prob,\mathbb{F}, X)$ be a stopping game.
\begin{itemize}
\item[O1)]
Let $\widetilde X$ be a payoff process such that 
(i) $\prob(|X_t - \widetilde X_t| > \ep \hbox{ for some } t \in \dN) \leq \ep$ and
(ii)
$\|X - \widetilde X\|_\infty \leq 1$
or $\|X\|_\infty, \|\widetilde X\|_\infty \leq 1$.
Then any $\ep'$-equilibrium in $\widetilde\Gamma = (\Omega,I,\calF,\prob,\mathbb{F}, \widetilde X)$
is a $(\ep' + 6\ep)$-equilibrium in $\Gamma$.
\item[O2)]
Let $\tau$ be a bounded stopping time, and let
$\widetilde\Gamma = (\Omega,I,\calF,\prob,\widetilde{\mathbb{F}}, \widetilde X)$ 
be the game where 
\[ \widetilde\calF_t := \calF_{\tau + t}, \qquad \forall t \in \dN \]
and
\[ \widetilde X_{J,t} := X_{J,\tau + t}, \qquad \forall t \in \dN. \]
Then any $\ep$-equilibrium in $\widetilde\Gamma$ can be extended to an $\ep$-equilibrium in $\Gamma$.
\item[O3)]
Let $\calD = (D_1,\dots,D_L)$ be a finite $\calF_0$-measurable partition of $\Omega$.
If for every $\ell=1,\dots,L$ the game 
$\Gamma_\ell = (\Omega,I,\calF,\prob(\cdot\mid D_\ell),\mathbb{F}, X)$ admits an $\ep$-equilibrium,
then $\Gamma$ admits an $\ep$-equilibrium.
\end{itemize}
\end{lemma}

Assertion (O1) holds because 
(i) affects each of $\gamma^i(x^i,x^{*,-i})$ and $\widetilde\gamma^i(x^i,x^{*,-i})$ by at most $\ep$,
and (ii) affects the difference between the two expected payoffs by at most $2\ep$.

Assertion (O2) holds because we can supplement any $\ep$-equilibrium $x^*$ in $\widetilde\Gamma$ by an initial play before stage $\tau$,
in which the players play a 0-equilibrium in the stopping game with random but bounded horizon $\tau$,
where the payoff if no player stops before stage $\tau$ is given by the payoff under $x^*$ in $\widetilde\Gamma$.

Assertion (O3) holds because 
if $x^*_\ell$ is an $\ep$-equilibrium in $\Gamma_\ell$,
then the strategy profile that coincides with $x^*_\ell$ on $D_\ell$ for each $\ell=1,\dots,L$ 
is an $\ep$-equilibrium in $\Gamma$.

%%%%%%%%%%%%%%%%%%%%%%%
\subsection{Conditional Atomlessness and Lotteries}
\label{sec:atomless2}

In this section, we prove
the next result,
which shows how the property of conditional atomlessness allows us to replace a lottery by the uncertainty in the probability space.
This result is the only place in our proof where the conditional atomlessness property is used.

\begin{lemma}[Simulating a lottery]
\label{lem:simulate}
Let $(\Omega,\calF,\prob)$ be a probability space,
let $\mathbb F = (\calF_t)_{t \in \dN}$ be a conditionally atomless filtration, let $\delta > 0$,
let $n \in \dN$,
let $t_0 \in \dN$,
and let $\nu = (\nu_i)_{i \in I} \in \Delta(\{1,2,\dots,n\}$).
There is $t_1 > t_0$ and an $\calF_{t_1}$-measurable partition $(C_i)_{i=1}^n$
such that for every $A\in \calF_{t_0}$ we have
$\bigl| \prob(A \cap C_i) - \prob(A) \nu_i\bigr| \leq \delta$.
\end{lemma}

The lemma states that the information provided by the filtration between stages $t_0$ and $t_1$ can be used to mimic a lottery on the set $\{1,\dots,n\}$ in a way that is almost independent from the information available at time $t_0$.

Note that if Lemma~\ref{lem:simulate} holds for $t_1$,
then it holds for every integer larger than $t_1$.

\bigskip

\begin{proof}
By Lemma~\ref{lemma:dalbaen},
there is a random variable $U \sim \mathrm{Uniform}[0,1]$ that is $(\vee_{t \in \dN} \calF_t)$-measurable and
independent of $\calF_{t_0}$.
Denote 
\[ B_i := \left\{ \sum_{j < i} \nu_j \leq U < \sum_{j \leq i} \nu_j \right\}, \ \ \ \forall i=1,\dots,n.  \]
Then
\[ \prob(B_i) = \prob\left(\sum_{j < i} \nu_j \leq U < \sum_{j \leq i} \nu_j\right) = \nu_i, \ \ \ \forall i=1,\dots,n. \]
Since $U$ is $(\vee_{t \in \dN} \calF_t)$-measurable,
we have $B_i \in \vee_{t \in \dN} \calF_t$ for every $i=1,\dots,n$.
Moreover, the events $(B_i)_{i=1}^{n}$ form a partition of $\Omega$.
Hence,
there are $t_1 > t_0$ and a partition $(C_i)_{i=1}^{n}$ of $\Omega$
which is $\calF_{t_1}$-measurable such that
\begin{equation}
\label{equ:eq:0}
\Prob(B_i \triangle C_i) \leq \delta, \ \ \ \forall i = 1,\dots,n.
\end{equation}
Lemma~\ref{lem:simulate} follows since
for every event $A \in \calF_{t_0}$
and every $i \in \{1,\dots,n\}$,
\begin{align}
\bigl| \prob(A \cap C_i) - \prob(A) \nu_i\bigr|
&\leq 
\bigl| \prob(A \cap B_i) - \prob(A) \nu_i\bigr|
\label{equ:eq:1}
+ \bigl| \prob(A \cap C_i) - \prob(A \cap B_i) \bigr|\\
\label{equ:eq:2}
&= \bigl| \prob(A \cap C_i) - \prob(A \cap B_i) \bigr|\\
&\leq \Prob(A \cap (B_i \triangle C_i))\\
\label{equ:eq:3}
&\leq \Prob(B_i \triangle C_i) \leq \delta,
\end{align}
where Eq.~\eqref{equ:eq:1} follows from the triangle inequality,
Eq.~\eqref{equ:eq:2} holds since $U$ is independent of $\calF_{t_0}$ and since $\prob(B_i) = \nu_i$,
and Eq.~\eqref{equ:eq:3} holds by Eq.~\eqref{equ:eq:0}.
\end{proof}

\subsection{Linear complementarity problems and Q-matrices}
\label{subsec:lcp}

In this section, we present a version of linear complementarity problems that is suitable for our study, and the corresponding concept of Q-matrices.

Denote by $\calM_{n \times n}[-1,1]$
the set of all $n \times n$ matrices $R=(R_{ij})_{1 \leq i,j\leq n}$
whose entries lie in $[-1,1]$.
Throughout the paper we restrict attention to matrices in this set.

\begin{definition}[Linear complementarity problems]
\label{def:lcp}
Let $n$ be a positive integer,
let $R = (R_{ij})_{1 \leq i,j \leq n}\in \calM_{n \times n}[-1,1]$,
and let
$q \in [-1,1]^n$.
The \emph{linear complementarity problem} $\LCP(R,q)$ consists of
finding $w \in [-1,1]^n \text{ and }
z=(z_0,z_1,\dots,z_n) \in \Delta(\{0,1,\dots,n\})$ such that
\begin{equation}
\label{eq:LCP}
w = z_0 q + \sum_{i=1}^n z_i R_{\bullet i},\qquad
w_i \geq R_{ii} \text{ for all } i \in \{1,\dots,n\},\qquad
z_i = 0 \text{ or } w_i = R_{ii} \text{ for all } i \in \{1,\dots,n\},
\end{equation}
where $R_{\bullet i}$ is the $i$'th column of $R$.
\end{definition}

\begin{remark}[Comparison with the standard formulation of LCPs]
\label{rem:lcp-shift}
The concept of linear complementarity problems defined here differs from the standard one 
(see \cite{CPS09}) in three respects.
In the standard formulation, 
\begin{itemize}
\item $z_0 = 1$ and $z_1,\dots,z_n$ are required to be non-negative.
\item The entries of $R$, $q$, and $w$ are not restricted to $[-1,1]$ but can be any real number.
\item 
$w_i$ is compared to 0 
(in the second and third conditions in \eqref{eq:LCP})
rather than to $R_{ii}$.
\end{itemize}
While linear complementarity problems that exhibit the first and third bullets above 
were studied in \cite{SS20},
our paper is the first in which the second point is also present.
\end{remark}

\begin{definition}[Q-matrix]
An $n \times n$ matrix $R\in \calM_{n \times n}[-1,1]$ is a \emph{Q-matrix} if for every
$q \in [-1,1]^n$ the problem $\LCP(R,q)$ has at least one solution.
\end{definition}

The next result states that if $R$ is not close to a $Q$-matrix,
then there is $q \in [-1,1]^n$ for which $\LCP(R,q)$ has no approximate solution.
This result is the reason for restricting the entries of $R$, $q$, and $w$ to compact sets in the definition of linear complementarity problems.

\begin{lemma}
\label{lemma:Qmatrix}
Let $n \geq 1$.
For every $\ep > 0$ there is $\delta = \delta_{\LCP}(\ep) > 0$
such that if $R \in \calM_{n \times n}[-1,1]$ is not $\ep$-close to any Q-matrix in the $L_\infty$ norm,
then
there is
$q \in [-1,1]^n$ such that there are no
$w \in [-1,1]^n$
and
$z=(z_0,z_1,\dots,z_n) \in \Delta(\{0,1,\dots,n\})$
for which
\begin{align}
\nonumber
w &= z_0 q + \sum_{i=1}^n z_i R_{\bullet i},\\
\label{equ:lcp:approx}
w_i &\geq R_{ii}-\delta \text{ for all } i \in \{1,\dots,n\},\\
z_i &= 0 \text{ or } |w_i - R_{ii}|\leq\delta \text{ for all } i \in \{1,\dots,n\}.
\nonumber
\end{align}
\end{lemma}

\begin{proof}
The proof uses a standard compactness-continuity argument.
Suppose by contradiction that the lemma is not true.
Then there is $\ep > 0$ such that for every $\delta > 0$
there are an $n \times n$ matrix $R^\delta \in \calM_{n\times n}[-1,1]$ that is not $\ep$-close to any Q-matrix,
and for every $q \in [-1,1]^n$ there are 
$w^{\delta,q} \in [-1,1]^n$
and
$z^{\delta,q} \in \Delta(\{0,1,\dots,n\})$
for which \eqref{equ:lcp:approx} holds for $R^\delta$.

Let $R$ be an accumulation point of $(R^\delta)_{\delta > 0}$,
which exists since $\calM_{n\times n}[-1,1]$ is compact.
If $R$ were a Q-matrix, then 
$R^\delta$ would be $\ep$-close to a Q-matrix for every $\delta$ sufficiently close to 0, which is a contradiction.
Hence $R$ is not a Q-matrix.

Fix $q \in [-1,1]^n$, and let $(w,z)$ be such that 
$(R,w,z)$ is 
an accumulation point of solutions
$(R^\delta,w^{\delta,q},z^{\delta,q})_{\delta > 0}$ of \eqref{equ:lcp:approx} as $\delta$ goes to $0$, which exists since $(w^{\delta,q})_{\delta > 0}$ and $(z^{\delta,q})_{\delta > 0}$ lie in compact sets.
Taking the limit in \eqref{equ:lcp:approx} 
for $R^\delta$, $w^{\delta,q}$, and $z^{\delta,q}$,
as $\delta$ goes to 0,
we deduce that $(w,z)$ is a solution of $\LCP(R,q)$.
Since this holds for every $q \in [-1,1]^n$,
$R$ is a Q-matrix, a contradiction.
\end{proof}

%When $q_i \geq R_{ii}$ for all $i$,
%the linear complementarity problem $\LCP(R,q)$ has a trivial solution, namely, 
%$z_0 = 1$ and $w = q$.

\subsection{Quitting games}
\label{subsec:quitting}

A \emph{quitting game} is the simplest class of stopping games: the filtration is trivial and the payoff process is time-homogeneous.

\begin{definition}[Quitting game]
\label{def:quitting}
A \emph{quitting game} $\calQ$ is a stopping game
$(\Omega,I,\calF,\prob,\mathbb{F}, X)$
that satisfies the following conditions:
\begin{itemize}
\item No uncertainty: $\calF_t = \{ \varnothing, \Omega\}$ for every $t \in \dN$.
\item Payoffs are time-homogeneous: $X_t = X_{t+1}$ for every $t \in \dN$.
We assume without loss of generality that the payoffs of all players lie in $(-1,1)$.
\end{itemize}
\end{definition}

Since the filtration is trivial and the payoffs are time-homogeneous, we identify a quitting game with the pair $\calQ = (I,X)$,
where $X = (X_J)_{\varnothing \neq J \subseteq I}$
and $X_J \in (-1,1)^n$ for every $\varnothing \neq J \subseteq I$.
Here, $X_J$ 
is the payoff vector when $J$ is the set of stopping players.%
\footnote{Recall that $X_{J,t}$ is the payoff vector in the stopping game when the set of stopping players in stage $t$ is $J$.}

Quitting games were introduced by \cite{SV01},
who provided sufficient conditions that guarantee the existence of an $\ep$-equilibrium for every $\ep > 0$.
Additional sufficient conditions have been provided by \cite{Simon07} and \cite{SS20}.
Note that the trivial filtration of a quitting game is not conditionally atomless.

Since in quitting games the filtration is trivial,
a behavioral strategy of a player reduces to a sequence 
$x^i = (x_t^i)_{t \in \dN}$,
where $x^i_t \in [0,1]$ for every $t \in \dN$,
with the interpretation that $x^i_t$ 
is the probability with which player $i$ stops at stage $t$, conditional on no player having stopped earlier.
\label{def:quitting:strategy}

A strategy $x^i$ is \emph{stationary} (or \emph{time-homogeneous})
if $x^i_t$ does not depend on
$t$. 
In this case, $x^i$ will represent the real number in $[0,1]$
such that $x^i_t = x^i$ for all $t \in \dN$.

For every quitting game $\calQ = (I, X)$,
define the matrix $R_\calQ \in \calM_{n \times n}[-1,1]$ of payoff vectors that correspond to single stoppings by letting its
$i$'th column be $X_{\{i\}}$, for each $i \in I$.
\label{def:quitting:rq}

Conversely,
every $n \times n$ matrix $R$ defines a quitting game $\calQ_R$, where the payoff vector if player~$i$ stops alone is $R_{\bullet i}$, the $i$'th column of $R$.\label{def:quitting:ri}
The payoff when more than one player stops is arbitrary;
we set it to $0$ for concreteness.
This choice will not affect the analysis: 
whenever we use $\calQ_R$,
we will consider only equilibria that involve termination by a single player. 

\begin{definition}[The quitting game $\calQ_R$]
\label{def:quitting:qr}
Let $R \in \calM_{n\times n}[-1,1]$.
Denote by $\calQ_R$ the quitting game 
$(I,X)$
where,
\[
X_{\{i\}} = R_{\bullet i} \quad \forall i \in I,
\qquad
X_J = (0,\dots,0) \quad \forall J \subseteq I \hbox{ with } |J| \ge 2.
\]
\end{definition}

\subsection{Sunspot equilibrium in quitting games}
\label{subsec:sunspot}

While the concept of Nash equilibrium is the most prevalent in game theory,
various alternative concepts have been proposed and studied.
One such concept is \emph{sunspot equilibrium},
which was introduced in \cite{CassShell}.
In this section, we define the concept of sunspot equilibrium in quitting games,
and recall a result of \cite{SS20},
which states that every quitting game admits a sunspot $\ep$-equilibrium for every $\ep > 0$.
We will use sunspot equilibria in quitting games
as a tool in our construction of an $\ep$-equilibrium in stopping games.

%%--------------------------------------------------------------------
\subsubsection{Sunspot equilibrium}
%%--------------------------------------------------------------------

A \emph{sunspot $\ep$-equilibrium} 
is an $\ep$-equilibrium in an extended game with a public
correlation device: at the beginning of every stage $t$, all players
observe a public signal $\psi_t$ that is drawn uniformly from
$[0,1]$, independently of all past signals and of the players' private randomizations.

\begin{definition}[Sunspot equilibrium]
\label{def:sunspot}
Let $\calQ$ be a quitting game.
The game $\calQ_{\mathrm{sun}}$ 
is the extension of $\calQ$ in which, at the beginning of every stage $t$, 
all players observe a public signal $\psi_t$ that is drawn uniformly from
$[0,1]$, independently of all past signals and of $(\zeta^i_t)_{t \in \dN, i \in I}$.
\end{definition}

A strategy of player $i$ in $\calQ_{\mathrm{sun}}$ is a sequence of measurable
functions $\xi^i = (\xi_t^i)_{t \in \dN}$, where
$\xi_t^i : [0,1]^{t+1} \to [0,1]$
for every $t \in \dN$.
The interpretation of $\xi^i$ is that, if no
player has stopped before stage $t$, 
then player $i$ stops at stage $t$ with
probability $\xi_t^i(\psi_0,\dots,\psi_t)$.
The domain of $\xi^i_t$ is $[0,1]^{t+1}$ rather than $\Omega \times [0,1]^{t+1}$
because the filtration in quitting games is trivial.
\label{def:sunspot:strategy}

Denoting by $\theta$ the (random) stage in which the game is stopped,\label{def:sunspot:theta}
and by $J_\theta$ the (random) set of players who stop at stage $\theta$,
the expected payoff under a strategy profile $\xi = (\xi^i)_{i \in I}$ is
\label{def:sunspot:i:delta}
$\gamma_{\mathrm{sun}}^i(\xi) := \E_{\xi}[X_{J_\theta}^i \cdot\mathbf{1}_{\theta < \infty}]$. 
A strategy profile $\xi^*$ is a
\emph{sunspot $\ep$-equilibrium} of $\calQ$ if it is an
$\ep$-equilibrium in $\calQ_{\mathrm{sun}}$,
namely, if
$\gamma_{\mathrm{sun}}^i(\xi^*) \geq \gamma_{\mathrm{sun}}^i(\xi^i, \xi^{*,-i}) - \ep$,
for every $i \in I$
and every strategy $\xi^i$ of player $i$ in $\calQ_{\mathrm{sun}}$.

%%--------------------------------------------------------------------
\subsubsection{The main result of Solan and Solan (2020)}
%%--------------------------------------------------------------------

The main result of \cite{SS20} is the following.

\begin{theorem}[{\cite[Theorem~1]{SS20}}]
\label{thm:SS20-main}
Every quitting game admits a sunspot $\ep$-equilibrium for
every $\ep > 0$.
\end{theorem}

In the proof of Theorem~\ref{thm:main} we will not use Theorem~\ref{thm:SS20-main} as a black box.
Rather we will use the explicit structure of the sunspot $\ep$-equilibrium strategy profiles.
We therefore describe the structure of these strategy profiles here.

For a quitting game $\calQ$,
the structure of the sunspot $\ep$-equilibrium constructed by \cite{SS20} depends on whether a certain principal submatrix of the matrix $R_\calQ$ is a Q-matrix.
In Section~\ref{sec:normal} we identify the principal submatrix of $R_\calQ$ that we will need to consider.
The construction when this principal submatrix is a Q-matrix
will be described in Section~\ref{sec:ss:Q}.
This construction will be extended to the case that this principal submatrix is close to a Q-matrix in Section~\ref{sec:ss:Q:close}.
Finally, in Section~\ref{sec:ss:nonQ} we will handle the case when this principal submatrix is not close to any Q-matrix.

%%--------------------------------------------------------------------
\subsubsection{Free-riders and supporters}
\label{sec:normal}
%%--------------------------------------------------------------------

We say that a player $i \in I$ is a free-rider in a quitting game $(I,X)$ if
conditioned on the event that at most one player ever stops,
player $i$ does not wish to stop.
That is, if 
$X_{\{i\}}^i < 0$
(player~$i$ prefers that no-one stops rather than she stops alone)
and
$X_{\{i\}}^i < X_{\{j\}}^i$ for every $j \neq i$
(player $i$'s payoff by stopping along is lower than her payoff when any other player $j \neq i$ stops alone).

Since we will need the concept of free-riders in stopping games, where the payoff of a player when she stops alone is a stochastic process,
we find it more convenient to define free-riders for matrices.

\begin{definition}[Free-riders and supporters]
\label{def:normal}
Let $R$ be an $n \times n$ matrix.
Player~$i \in I$ is a \emph{free-rider} for $R$ if 
$R_{ii} < 0$ and 
$R_{ii} < R_{ij}$ for every $j \neq i$.
Player $i$ is a \emph{supporter} for $R$ otherwise.
\end{definition}

For a quitting game $\calQ = (I,X)$ we have
$R_{\calQ,ij} = X_{\{j\}}^i$ for every $1 \leq i,j \leq n$.
Hence the concept of free-riders as defined in Definition~\ref{def:normal} for the matrix $R_\calQ$ coincides with the intuitive definition described before the definition.

Let $I^{\mathrm{sp}}_R$ be the set of supporters for $R$.
The next result states that if all players are free-riders, 
then an $\ep$-equilibrium exists for every $\ep > 0$.\label{def:player:nor}

\begin{lemma}[Existence of a $0$-equilibrium when all players are free-riders; \cite{SS20}]
\label{lemma:abnormal}
Let $\calQ$ be a quitting game.
If $I^{\mathrm{sp}}_{R_\calQ} = \varnothing$, then $\calQ$ admits a 0-equilibrium.
\end{lemma}

Indeed, 
the strategy profile in which no player ever stops is a 0-equilibrium in $\calQ$.
This holds because:
\begin{itemize}
\item 
Under the strategy profile in which no player ever stops, 
the payoff of all players is 0.
\item 
A unilateral deviation by player $i$ can only lead to unilateral stopping by that player, which yields to player~$i$ the payoff $X_{\{i\}}^i = R_{ii}<0$ since the player is a free-rider. 
\end{itemize}
Hence no player can profit from deviating.

\begin{notation}[The principal submatrix $R^{\mathrm{sp}}_\calQ$]
\label{def:matrix:nor}
Denote by $R^{\mathrm{sp}}_\calQ$ the principal submatrix of $R_\calQ$ that is indexed by $I^{\mathrm{sp}}_{R_\calQ}$;
that is, it is derived from $R_\calQ$
by deleting all rows and columns that 
correspond to free-riders.
\end{notation}

%%--------------------------------------------------------------------
\subsubsection{The construction of a sunspot $\ep$-equilibrium when $R^{\mathrm{sp}}_\calQ$ is a Q-matrix}
\label{sec:ss:Q}
%%--------------------------------------------------------------------

The space $[0,1]^\dN$ is a probability space when equipped with the product sigma-algebra and the countable product of the Lebesgue measure.
Endow it with the natural filtration ${\mathbb F}=(\calF_t)_{t \in \dN}$;
that is, 
$\calF_t$ is the sigma-algebra generated by the cylinder sets defined by the first $(t+1)$ coordinates for every $t \in \dN$.

When $R^{\mathrm{sp}}_\calQ$ is a Q-matrix, 
the sunspot $\ep$-equilibrium $\xi^*$ constructed by \cite{SS20} has the following structure.
There exist two adapted processes on $[0,1]^\dN$:
a process $(\widehat i_t)_{t \in \dN}$ taking values in $I \cup \{\varnothing\}$, 
and a process $(\widehat\delta_t)_{t \in \dN}$ taking values in $(0,\ep]$,
such that the following strategy profile $\xi^*$ is  a sunspot $\ep$-equilibrium
and terminating with probability 1 in all subgames:%
\footnote{Since $\widehat i_t$
and $\widehat \delta_t$ are $\calF_t$-measurable, 
they can be identified with functions of $(\psi_0,\dots,\psi_t)$.}
\[ \xi^{*,i}_{t}(\psi_0,\dots,\psi_t) := 
\left\{\begin{array}{ll}
0, & i \neq \widehat i_t(\psi_0,\dots,\psi_t),\\
\widehat \delta_t(\psi_0,\dots,\psi_t), & i = \widehat i_t(\psi_0,\dots,\psi_t).
\end{array}
\right.
\]
That is,
under $\xi^*$, after every finite sequence of signals $(\psi_0,\dots,\psi_t)$, 
player~$\widehat i_t(\psi_0,\dots,\psi_t)$ stops with probability $\widehat\delta_t(\psi_0,\dots,\psi_t)$
and all other players stop with probability 0.

The construction of the processes $(\widehat i_t)_{t \in \dN}$ and $(\widehat\delta_t)_{t \in \dN}$ is far from trivial, but,
as \cite{SS20} show,
they can be defined in such a way that the following properties hold:%
\footnote{The keen reader will note the discrepancy between the conditions and conclusions of \cite{SS20} and those we listed here.
With regard to the conditions, \cite{SS20} define the set $I_{R_\calQ}^{\mathrm{sp}}$ differently.
In that work, this set contains all players whose minmax value is non-positive. 
Following their proof closely reveals that their result holds under our definition of $R_\calQ^{\mathrm{sp}}$ as well.

With regard to the conclusion, 
for every given $\eta > 0$ \cite{SS20} construct $\xi^*$ so that under $\xi^*$ the play terminates with probability at least $1-\eta$.
Applying a standard limiting argument to that construction yields an $\ep$-equilibrium $\xi^*$ under which the play terminates with probability 1.

Another discrepancy concerns the definition of the linear complementarity problem.
As mentioned in Remark~\ref{rem:lcp-shift},
in the definition of linear complementarity problems of \cite{SS20}, the vectors $q$ and $w$ can be any vectors in $\dR^n$,
and are not restricted to $[-1,1]^n$.
When the payoffs in the quitting game lie in $(-1,1)$,
as we assumed,
the proof in \cite{SS20} does not use vectors $q$ and $w$ that do not lie in $[-1,1]^n$,
and hence our definition of linear complementarity problems suffices to apply \cite{SS20}.}
\begin{enumerate}
\item 
$\widehat i_t$ 
is a measurable function of $(\widehat i_0,\dots,\widehat i_{t-1},\psi_t)$,
and $\widehat \delta_t$ is a function of $(\widehat i_0,\dots,\widehat i_{t-1},\widehat i_t)$, 
for every $t \in \dN$.
\item
For every $t \in \dN$,
no player can profit more than $\ep$ by deviating from $\xi^{*,i}$ in the subgame that starts at $t$.
\end{enumerate}

\subsubsection{The construction of a sunspot equilibrium when $R^{\mathrm{sp}}_\calQ$ is $\ep$-close to a Q-matrix}
\label{sec:ss:Q:close}

We now argue that the construction in Section~\ref{sec:ss:Q} applies, albeit with a larger $\ep$, 
when $R^{\mathrm{sp}}_\calQ$ is close to a Q-matrix.

Suppose that $R^{\mathrm{sp}}_\calQ$ is $\ep$-close to a
Q-matrix $\widetilde R$ of dimensions $|I^{\mathrm{sp}}_\calQ| \times |I^{\mathrm{sp}}_\calQ|$.
Extend $\widetilde R$ to an $(n \times n)$-matrix $R'$ by filling the missing entries with their values in $R_\calQ$.
Any free-rider according to $R_\calQ$ is also a free-rider according to $R'$.
Hence, 
the set of free-riders for $R'$ weakly contains $I^{\mathrm{sp}}_\calQ$.

In Section~\ref{sec:ss:Q} we described the sunspot $\ep$-equilibrium for the quitting game $\calQ_{R'}$.%
\footnote{We do not rule out the case that some players who are supporters for $R_\calQ$ become free-riders for $R'$.
When this happens, $R'^{\mathrm{sp}}$ is a principal matrix of $\widetilde R$,
and $R'^{\mathrm{sp}}$ might be far from all $Q$-matrices.
Nevertheless, the construction of \cite{SS20} applies to the game $\calQ_{R'}$ with the matrix $\widetilde R$ taking the role of the principal submatrix of supporters.}
The payoffs upon single stopping in the games $\calQ_{R'}$ and $\calQ$ differ by at most $\ep$.
Moreover, the sunspot $\ep$-equilibrium $\xi^*$ that was constructed in Section~\ref{sec:ss:Q} for $\calQ_{R'}$ involves only single stoppings, and a player who stops at some stage,
does so with probability at most $\ep$.
Hence, if a player deviates from $\xi^*$,
the difference between her payoffs in the games $\calQ_{R'}$ and $\calQ$ is at most $2\ep$.
As in the proof of Lemma~\ref{lemma:simple}(O1),
$\xi^*$ is a sunspot $4\ep$-equilibrium in $\calQ$.

%%--------------------------------------------------------------------
\subsubsection{The construction of a sunspot equilibrium when $R^{\mathrm{sp}}_\calQ$ is not a Q-matrix}
\label{sec:ss:nonQ}
%%--------------------------------------------------------------------

By Section 3.2 in \cite{SS20},
when $R^{\mathrm{sp}}_\calQ$ is not a Q-matrix, the quitting game $\calQ$ admits a stationary approximate equilibrium, without the need to use correlation.%
\footnote{The $\ep$-equilibrium constructed in \cite{SS20}[Section 3.2] when $R^{\mathrm{sp}}_\calQ$ is not a Q-matrix was not stationary, because they defined $I_\calQ^{\mathrm{sp}}$ as the set of players whose minmax value is non-positive,
and so punishment was executed by switching to a minmax strategy.
In our case, any player in $I_\calQ^{\mathrm{sp}}$ can be punished by a perturbation of the vector in which no player stops, and hence, as we explain below, there is a stationary $\ep$-equilibrium.}

The idea of the proof is as follows.
Since $R^{\mathrm{sp}}_\calQ$ is not a Q-matrix,
there is a vector $q^* \in [-1,1]^{|I^{\mathrm{sp}}_{R_\calQ}|}$ such that $\LCP(R^{\mathrm{sp}}_\calQ,q^*)$ has no solution.
Extend $q^*$ to a vector $q$ in $[-1,1]^n$,
by letting $q^i$ be any real number larger than all payoffs of player~$i$, for
every $i \not\in I^{\mathrm{sp}}_{R_\calQ}$:
\[ q^i := \left\{
\begin{array}{ll}
q^{*,i}, & i \in I^{\mathrm{sp}}_{R_\calQ},\\
1, & i\not\in I^{\mathrm{sp}}_{R_\calQ}.
\end{array}
\right.
\]
The large value of $q$ in coordinates that correspond to free-riders guarantees that those players do not
stop when the per-stage stopping probability of the other players is low.

For every $\lambda \in (0,1]$, 
consider an auxiliary quitting game, 
where
(a) the payoffs are $\lambda$-discounted,
and (b) the stage payoff is $q$.
Formally, the payoff under a stationary strategy profile $x$ is
\[
\label{def:payoff:discounted}
\gamma_\lambda(x) := \E_x\left[ (1-\lambda)^\theta X_{J_\theta} + \bigl(1-(1-\lambda)^\theta\bigr) q \right]. \]
By \cite{Fink64},
a stationary equilibrium exists:
For every $\lambda \in (0,1]$ there is a stationary strategy profile $x^*_\lambda \in [0,1]^n$ such that 
for every $i \in I$ and every strategy $x^i$ of player $i$ (not necessarily stationary),
$\gamma^i_\lambda(x^i,x^{*,-i}_\lambda) \leq \gamma^i_\lambda(x^*_\lambda)$.
By \cite{Neyman03},
the function $\lambda \mapsto x^*_\lambda$ can be taken to be semi-algebraic.

Since $\lambda \mapsto x^*_\lambda$ is semi-algebraic,
the limit
$x^{*}_0 := \lim_{\lambda \to 0} x^*_\lambda$ exists,
as well as all other limits we consider below.%
\footnote{The use of semi-algebraic properties here is for convenience only. 
We could have instead consider a subsequence of $(x^*_{\lambda})_{\lambda> 0}$ as $\lambda$ goes to 0 where all limits that we need converge.}
\cite[Section 3.2]{SS20} show that $x^{*}_0$ is a 0-equilibrium, or a perturbation of this strategy profile is an $\ep$-equilibrium.
The argument involves four steps.

First, standard continuity arguments show that for every mapping $\lambda \mapsto x_\lambda$ such that $x_0 := \lim_{\lambda \to 0} x_\lambda$ exists and $\prob_{x_0}(\theta < \infty) = 1$,
we have 
\begin{equation}
\label{equ:limit:abs}
\lim_{\lambda \to 0} \gamma_\lambda(x_\lambda) = \gamma(x_0), 
\end{equation}
see, e.g., \cite[Lemma 4]{VT89} or \cite[Lemma 6.6]{Solan99}.

Second, it is proved that $x^*_0 \neq (0,\dots,0)$.
Indeed, assume by contradiction that $x^*_0 = (0,\dots,0)$.
To derive a contradiction, we show that $\LCP(R^{\mathrm{sp}}_\calQ,q^*)$ has a solution.
Set
\[ z^i_\lambda := \frac{x^{*,i}_\lambda}{\lambda+\sum_{j \in I} x^{*,j}_\lambda}, \ \ \ i \in I, 
\qquad  
z_\lambda^0 := \frac{\lambda}{\lambda+\sum_{j \in I} x^{*,j}_\lambda}.
\]
Set also 
\[ z^i_0 := \lim_{\lambda \to 0} z^i_\lambda, \ \ \ \forall i \in I,
\ \ \ 
z^0_0 := \lim_{\lambda \to 0} z^0_\lambda, \ \ \ 
w := \lim_{\lambda \to 0} \gamma_\lambda(x^*_\lambda).
\]
As in \cite[Lemma 2]{VT89},
\begin{equation}
\label{equ:45}
w = \lim_{\lambda \to 0} \gamma_\lambda(x^*_\lambda) = \lim_{\lambda \to 0} \left( z_\lambda^0 q + \sum_{i \in I} z_\lambda^i R_{\bullet i} \right) = z^0_0 q + \sum_{i \in I} z^i_0 R_{\bullet i},
\end{equation}
and
\begin{equation}
\label{equ:45:a}
 w^i = \lim_{\lambda \to 0} \gamma_\lambda^i(x^*_\lambda) 
\geq \lim_{\lambda \to 0} \gamma_\lambda^i(1^i,x_\lambda^{*,-i}) = R_{ii}, 
\end{equation}
where $1^i$ is the strategy of player~$i$ in which she always stops: $1^i_t \equiv 1$ for every $t \in \dN$.
Moreover, if there is a strict inequality in \eqref{equ:45:a},
then $x_\lambda^{*,i} = 0$ for every $\lambda > 0$ sufficiently small,
so that $z^i_0 = 0$.
It follows that
$z^i_0 = 0$ for every $i \not\in I^{\mathrm{sp}}_{R_\calQ}$.
Therefore,
$(w,(z^0_0,(z^i_0)_{i\in I}))$ satisfies \eqref{eq:LCP},
and hence it is a solution of $\LCP(R^{\mathrm{sp}}_\calQ,q^*)$.

From this, the existence of an $\ep$-equilibrium follows.
Suppose first that there are at least two players $j \in I$ such that $x^*_j > 0$.
Since $x^*_{\lambda}$ is a $\lambda$-discounted equilibrium,
and by \eqref{equ:limit:abs},
for every $i \in I$ and every stationary strategy $x^i$ of that player,
\[ \gamma^i(x^*) = \lim_{\lambda\to 0} \gamma^{i}_\lambda(x^*_\lambda)
\geq \lim_{\lambda\to 0} \gamma^*_\lambda(x^i, x^{*,-i}_\lambda) = \gamma^i(x^i,x^{*,-i}_0). \]
Since $i$ and $x^i$ are arbitrary, 
this implies that $x^*_0$ is a 0-equilibrium of $\calQ$.

Suppose next that there is exactly one player $i_0 \in I$ such that $x^{*,i_0}_{0} > 0$.
Then necessarily $i_0 \in I^{\mathrm{sp}}_{R_\calQ}$.
Hence,
$X_{\{i_0\}}^{i_0} \geq 0$ or
there is a player $i_1$ such that $X_{\{{i_0}\}}^{i_0} \geq X_{\{i_1\}}^{i_0}$.

If $X_{\{i_0\}}^{i_0} \geq 0$,
then $x^*_0$ is a 0-equilibrium.
Indeed,
as above,
no player $i \neq i_0$ can profit by deviating,
and since $X_{\{i_0\}}^{i_0} \geq 0$,
player~$i_0$ cannot profit by never stopping,
in which case her payoff is 0.

Assume then that $X_{\{i_0\}}^{i_0} \geq X_{\{i_1\}}^{i_0}$ for some $i_1 \neq {i_0}$.
Consider the following stationary perturbation of $x^*_0$, denoted $\widehat x$, where $\ep \in (0,1/2)$:
\[ \widehat x^i := \left\{
\begin{array}{ll}
0, & \hbox{if } i \neq i_0,i_1,\\
x_0^{*,i_0}, & \hbox{if } i={i_0},\\
\ep x_0^{*,i_0}, & \hbox{if } i=i_1.
\end{array}
\right.
\]
Under $\widehat x$, the game eventually terminates by $i_0$, $i_1$, or both.
Moreover,
\[ \frac{\prob_{\widehat x}(J_\theta = \{i_1\} \hbox{ or } J_\theta = \{i_0,i_1\})}{\prob_{\widehat x}(J_\theta = \{i_0\})} = \frac{\ep}{1-\ep x_0^{*,i_0}} \leq \frac{\ep}{1-\ep} \leq 2\ep.
\]
It follows that 
$|\gamma(\widehat x) - X_{\{i_0\}}| \leq 4\ep$.
As above, $X^i_{\{i\}} \leq X^i_{\{i_0\}}$ for every $i \neq i_0$, and hence no player can profit more than $6\ep$ by stopping at stage 1.
Moreover, player $i_1$ cannot profit more than $4\ep$ by never stopping.
Since $X_{\{i_0\}}^{i_0} \geq X_{\{i_1\}}^{i_0}$,
player~$i_0$ cannot profit by never stopping as well.
It follows that $\widehat x$ is a $10\ep$-equilibrium.

%%=====================================================================
\subsection{Stopping Games on a Finite Tree}
\label{sec:trees}
%%=====================================================================

In this section, we define stopping games on finite trees, in which play restarts from the root whenever a leaf is reached.
This class of games is a subclass of stochastic games (\cite{Sha53}).

\subsubsection{Definition of games on a finite tree}

\begin{definition}[Stopping games on a finite tree]
\label{def:tree-game}
A \emph{stopping game on a finite tree}
(or simply a \emph{game on a tree}) is a tuple
$\calT = \bigl(I, S, S_1, r,
    (C_s)_{s \in S \setminus S_1},
    (p_s)_{s \in S \setminus S_1},
    (X_s)_{s \in S \setminus S_1}\bigr)$,
where:
\begin{enumerate}[(i)]
\item $I = \{1,\dots,n\}$ is the set of players.
\item $(S, S_1, r, (C_s)_{s \in S \setminus S_1})$ is a finite connected rooted
  tree: $S$ is a finite nonempty set of \emph{nodes},
  $S_1 \subseteq S$ is the nonempty set of \emph{leaves},
  $r \in S$ is the \emph{root}, and for each non-leaf node
  $s \in S_0 := S \setminus S_1$, the set
  $C_s \subseteq S \setminus \{r\}$ is the nonempty set of
  \emph{children} of $s$
  (so that $(C_s)_{s \in S_0}$ is a partition of $S \setminus\{r\}$).
%  We write $\mathrm{depth}(s)$ for the length of the unique path from
%  $r$ to $s$.
\item For each $s \in S_0$, $p_s$ is a probability distribution over
  $C_s$.
\item For each $s \in S_0$, the \emph{payoff} at $s$ is a vector
  $X_s = (X_{J,s}^i)_{\varnothing \neq J \subseteq I, i \in I} \in \dR^{n \cdot (2^n - 1)}$.
  Here, $X^i_{J,s}$ is the payoff to player~$i$ when the set of players
  who stop at node~$s$ is~$J$.
\end{enumerate}
\end{definition}
 
The game on a tree is played in \emph{rounds} until it terminates.
A round starts at the
root~$r$.
In any given stage,
when the game is at some node $s \in S_0$, all players simultaneously and
independently decide whether to \emph{stop} or \emph{not to stop}.
Let $J \subseteq I$ be the set of players who decided to stop.
\begin{itemize}
\item If $J \neq \varnothing$, the game terminates,
  and the terminal payoff to each player~$i \in I$ is $X^i_{J,s}$.
\item If $J = \varnothing$, a child $s' \in C_s$ is drawn according to
  $p_s$, and the play moves to $s'$.
\end{itemize}
If $s' \in S_1$ is a leaf, the round ends without termination and a new
round begins at the root.
Thus, players cannot stop at leaves; the game is \emph{periodic} in the sense that
reaching a leaf starts a new round.

\begin{remark}[Games on a tree that do not start at the root]
Below we will also consider cases in which the game on a tree starts at some node $s \in S_0$ other than the root. However, all rounds after the first round always start at the root.
\end{remark}

Similarly to stopping games,
the maximal payoff that player~$i$ can obtain by stopping alone in $\calT$ is
\begin{equation}
\label{equ:y:finite:tree}
y^i_\calT := \max_{s \in S_0} X^i_{\{i\},s}. 
\end{equation}
We will consider here only games on finite trees that satisfy the following genericity assumption.

\begin{assumption}
\label{ass:tree}
For each $i \in I$, 
if $X^i_{J,s} = X^i_{J',s'}$
for some distinct pairs $(J,s)$ and $(J',s')$,
then $X_{s} = X_{s'}$.
\end{assumption}

When the game on a tree $\calT$ satisfies Assumption~\ref{ass:tree},
for each $i \in I$
we can denote by $Z_{\calT,i}$ the vector $X_{\{i\},s}$,
where $s$ is any node for which there is an equality in \eqref{equ:y:finite:tree}.
Denote by $R_\calT$ the $n \times n$ matrix whose $i$'th column is $Z_{\calT,i}$, 
for each $i \in I$.
The set of supporters for $R_\calT$ is
$I^{\mathrm{sp}}_{R_\calT}$.
Let $R^{\mathrm{sp}}_\calT$ be the principal submatrix of $R_\calT$ indexed by $I^{\mathrm{sp}}_{R_\calT}$.
\label{def:Z:tree}

We will also denote by
\[ 
\label{def:tree:delta}
\delta_\calT := \min\bigl\{y^i_\calT - X^i_{\{i\},s}
\colon i\in I, s \in S_0 \hbox{ such that } y^i_\calT - X^i_{\{i\},s} > 0\bigr\} > 0 \]
the minimal positive gap between $y^i_\calT$ and a strictly smaller payoff that player $i$ can obtain by stopping alone.
If the set over which the minimum is taken is empty,
set $\delta_\calT = \infty$.
By Assumption~\ref{ass:tree},
when $\delta_\calT = \infty$, the payoff vector when player~$i$ quits alone is the same at all nodes in $S_0$.

\subsubsection{Strategies, payoffs, and equilibrium}
\label{subsection:tree:strategy}

A \emph{stationary strategy} of player~$i$ in $\calT$ is a function
$\alpha^i : S_0 \to [0,1]$, where $\alpha^i(s)$ is the probability that
player~$i$ stops at node~$s$.
A non-stationary strategy is a function from past histories to $[0,1]$.
%We identify a stationary strategy profile
%$\alpha = (\alpha^i)_{i \in I}$ with a vector in $[0,1]^{n \cdot |S_0|}$.

\begin{notation}[$\prob_{s,\alpha}$, $\E_{s,\alpha}$, $\prob_{s,\mathbf{0}}(M)$]
\label{def:tree:prob}
Denote by $\prob_{s,\alpha}$ the distribution induced by
the strategy profile~$\alpha$
when the initial node (in the initial round) is $s$.
Denote by $\E_{s,\alpha}$ the corresponding
expectation operator.
For every subset $M \subseteq S$,
denote by $\prob_{s,\mathbf{0}}(M)$ 
the probability that, under the strategy profile $\mathbf{0} = (0^i)_{i \in I}$ in which no player ever stops, 
the play reaches some node in $M$ during the first round when the initial node is $s$.
\end{notation}

\begin{notation}[${\theta_{\mathrm{fr}}}$, $\pi(\alpha)$, $\rho^i(\alpha)$, $\gamma^i(\alpha)$]
\label{def:tree:theta}
Let $\theta_{\mathrm{fr}}$ denote the last stage of the first round when the initial node is the root $r$;
that is, the first stage in which a leaf is reached when the players follow the strategy profile $\mathbf 0$.
%Let ${\widehat \theta}$ denote the stopping stage within the first round.
%Thus, ${\widehat \theta} = \infty$
%if no player stops during the first round.
%Denote by $s_{\widehat \theta}$ the node at stage ${\widehat \theta}$ and by
%$J_{\widehat \theta}$ the set of players stopping at stage ${\widehat \theta}$.
For every stationary strategy profile $\alpha$, 
the probability of termination in one round is
\[  \pi(\alpha) := \prob_{r,\alpha}({\theta} < \theta_{\mathrm{fr}}), \]
the expected payoff to each player~$i$ in one round,
unconditional on termination in that round, is
\[  \rho^i(\alpha) := \E_{r,\alpha}[X_{J_{\theta},s_{\theta}}^{i} \mathbf{1}_{{\theta} < \theta_{\mathrm{fr}}}], \]
and the expected payoff under the stationary strategy profile $\alpha$ is
\[  \gamma^i(\alpha) := \rho^i(\alpha)/\pi(\alpha), \]
where by convention $0/0 = 0$.
\end{notation}

The quantity $\gamma^i(\alpha)$ is well-defined and satisfies
$\pi(\alpha) \cdot \gamma^i(\alpha) = \rho^i(\alpha)$.
When we wish to emphasize the dependence of $\pi(\alpha)$, $\rho^i(\alpha)$, and $\gamma^i(\alpha)$ on $\calT$, we will write
$\pi(\alpha;\calT)$, $\rho^i(\alpha;\calT)$, and $\gamma^i(\alpha;\calT)$, respectively.

\begin{definition}[$\ep$-equilibrium]
A stationary strategy profile $\alpha^* = (\alpha^{*,i})_{i \in I}$ is an
\emph{$\ep$-equilibrium} of $\calT$ if for every player
$i \in I$ and every stationary strategy $\alpha^i$,
\[
  \gamma^i(\alpha^*) \geq \gamma^i(\alpha^i, \alpha^{*,-i}) - \ep.
\]
\end{definition}

Since $\alpha^{*,-i}$ is stationary,
the optimization problem faced by player $i$ is reduced to a Markov decision problem.
Hence,
if player~$i$ has no deviation to a stationary strategy that yields a profit higher than $\ep$,
then she has no deviation to any strategy that yields a profit higher than $\ep$.

\subsubsection{Terminating approximate equilibria in games on finite trees}

We now state a result that 
asserts that when the matrix $R^{\mathrm{sp}}_\calT$ is not close to any Q-matrix,
the game in a tree $\calT$ admits an approximate equilibrium where the per-round probability of termination is bounded away from 0.
We first restrict attention to a bounded family of games.

\begin{definition}[The family $\mathbb{T}_\ep$]
\label{def:family:t}
Let $\varepsilon > 0$.
Denote by $\mathbb{T}_\ep$ the family of all stopping games on a finite
tree
$\calT = \bigl(I,S,S_1,r,(C_s)_{s\in S_0},
(p_s)_{s\in S_0},(X_s)_{s\in S_0}\bigr)$
that satisfy the following properties:
\begin{itemize}
\item[T1)] $\|X_s\|_\infty<1$ for every $s\in S_0$.
\item[T2)] Every payoff $X^i_{J,s}$ is an integer multiple of $\ep$.
\item[T3)] Assumption~\ref{ass:tree} holds.
\item[T4)] $R^{\mathrm{sp}}_\calT$ is not $\ep$-close to a Q-matrix.
\item[T5)] Every node is reached with positive probability under
$\mathbf0$.
\item[T6)] For every $g \in \range(X)$,
\[
 \prob_{r,\mathbf0}
 \bigl(X_{s_t}=g
       \text{ for some }t < \theta_{\mathrm{fr}}\bigr)
 = 1.
\]
\end{itemize}
\end{definition}

The next result states that all games in $\mathbb T_\ep$ admit a stationary $54\sqrt\ep$-equilibrium where the per-round probability of termination is bounded away from 0.
Moreover, if some player obtains in this $54\sqrt\ep$-equilibrium a low payoff,
then the per-round probability of termination by some other player is bounded away from 0.

\begin{proposition}[Terminating approximate equilibrium for games in $\mathbb T_\ep$]
\label{prop:absorbing-eq}
For every $\ep>0$ sufficiently small there is a constant $d>0$,
depending only on $\ep$ and the number of players, such that every
$\calT\in\mathbb T_\ep$ 
that does not satisfy the conditions of Lemma~\ref{lemma:simple:eq}
admits a stationary $54\sqrt\ep$-equilibrium
$\alpha^*$ satisfying the following properties:
\begin{enumerate}[(i)]
\item $\pi(\alpha^*)\ge d$.
\item If $\gamma^i(\alpha^*)<y^i_\calT-7\sqrt\ep$,
then $\pi(0^i,\alpha^{*,-i})\ge d^2$.
\item If $\gamma^i(\alpha^*)\le-7\sqrt\ep$ and
$\pi(0^i,\alpha^{*,-i})<d^2$, then there are
$s\in S_0$ and $(\beta^j)_{j\ne i} \in [0,1]^{n-1}$, not all zero, for which
\[
\max\{A,B\}\le\gamma^i(\alpha^*)+\ep^2,
\]
where, 
\[
 A=\sum_{J\subseteq I\setminus\{i\}}
 X^i_{J\cup\{i\},s}
 \prod_{j\in J}\beta^j
 \prod_{j\in I\setminus(J\cup\{i\})}(1-\beta^j),
\]
and
\[
 B=\frac{\displaystyle\sum_{\varnothing\ne J\subseteq I\setminus\{i\}}
 X^i_{J,s}\prod_{j\in J}\beta^j
 \prod_{j\in I\setminus(J\cup\{i\})}(1-\beta^j)}
 {\displaystyle\sum_{\varnothing\ne J\subseteq I\setminus\{i\}}
 \prod_{j\in J}\beta^j
 \prod_{j\in I\setminus(J\cup\{i\})}(1-\beta^j)}.
\]
Moreover, there is at most one player satisfying the two inequalities in the
antecedent of this assertion.
\end{enumerate}
\end{proposition}

Condition~(ii) in Proposition~\ref{prop:absorbing-eq} ensures that if under $\alpha^*$ a player obtains a low payoff,
and if that player decides never to stop,
then the probability that some \emph{other} player stops is still bounded away from 0.

Condition~(iii) in Proposition~\ref{prop:absorbing-eq} ensures that if under $\alpha^*$ a player obtains a negative payoff,
and if, when that player decides not to stop the termination probability in a single round is low,
then there is a node $s \in S_0$
and $(\beta^j)_{j \neq i} \in [0,1]^{n-1}$, not all of them 0,
such that player~$i$'s payoff under $\alpha^*$ is not much lower than $A$ and $B$.
This node will allow us to force player $i$ to follow $\alpha^{*,i}$, 
by threatening her that if the game does not terminate in some given number of stages,
she will be punished.

The proof of Proposition~\ref{prop:absorbing-eq} follows the proof of \cite[Proposition~5.5]{SS04}.
However, the latter heavily uses the fact that the setup of \cite{SS04} includes only two players.
We will therefore need to adapt the proof,
using ideas from \cite{SS20}, 
to take into account the presence of more than two players.
To obtain the uniform bound in Proposition~\ref{prop:absorbing-eq} 
we need the variant of the linear complementarity problem of Definition~\ref{def:lcp}.
The proof requires several concepts and intermediate results that are independent of the rest of the paper,
and is therefore deferred to Section~\ref{sec:proof:absorbing-eq}.

%%%%%%%%%%%%%%%%%%%%%%%
\subsection{$\tau$-Functions and a Stochastic Variation of Ramsey's Theorem}
\label{sec:nt}

\cite{Ramsey} proved that every complete infinite graph
whose edges are colored by finitely many colors
admits a complete infinite monochromatic subgraph.
In this section, we present a stochastic variation of this result due to \cite{SS04}.

Consider the probability space $\{0,1\}^\dN$
that corresponds to tossing a fair coin infinitely many times.
This probability space can be presented as an infinite complete binary tree,
where a node of depth $t$ corresponds to the outcomes in the first $t+1$ tosses.
A finite portion of this tree is defined by a root
(which is defined by a finite (possibly empty) sequence of outcomes $(\omega_0,\dots,\omega_{t-1}) \in \{0,1\}^t$)
and finitely many leaves.

Every pair of bounded stopping times $\tau < \tau'$ and every $\omega \in \{0,1\}^\dN$ define a finite portion of the binary tree:
the root is $(\omega_0,\dots,\omega_{\tau(\omega)})$,
and the leaves are all nodes of the form
$(\omega'_0,\dots,\omega'_{\tau'(\omega')})$,
for any $\omega' \in \{0,1\}^\dN$
that extends the root.

Given an increasing sequence of bounded stopping times 
$0 \leq \tau_0 < \tau_1 < \dots$,
every $\omega \in \{0,1\}^\dN$ passes through countably many portions of the binary tree,
where the root of each one is a leaf of the former.

Suppose that we are given a finite set $\calA$ of \emph{colors},
and that we assign to each such finite portion of the binary tree a color in $\calA$.
\cite{SS04} proved that for every $\ep > 0$,
there is an increasing sequence of stopping times 
$0 \leq \tau_0 < \tau_1 < \dots$ such that,
with probability at least $1-\ep$,
the \emph{same} color is assigned to the finite portions that $\omega$ passes through, except, possibly, for the first portion.
Moreover, they proved the result for every filtered probability space
and not only for the complete binary tree.

We will need a trivial extension of this result,
where the increasing sequence of stopping times cannot be arbitrary.
We will now present the extension we need.

\begin{definition}[$\tau$-functions]
\label{def:NT}
Let $(\Omega,\calF,\prob)$ be a probability space,
let $\mathbb{F} = (\calF_t)_{t \in \dN}$ be a filtration,
and let $(n_k)_{k\in \dN}$
be an increasing sequence of non-negative integers.

An $\mathbb{F}$-adapted stopping time $\tau$ is \emph{$(n_k)_{k\in \dN}$-consistent}
if $\tau \in \{ n_k, k \in \dN\}$ a.s.

A \emph{$\tau$-function} (for $(n_k)_{k\in \dN}$)
is a function $f$ that assigns, to every two $(n_k)_{k\in \dN}$-consistent bounded stopping times $\tau < \tau'$,
an $\calF_{\tau}$-measurable random variable $f_{\tau,\tau'}$
that satisfies the following consistency condition:
For every three $(n_k)_{k\in \dN}$-consistent bounded stopping times $\tau < \tau',\tau''$
and every set $D \in \calF_\tau$ such that 
$\{\tau' = \tau''\}$ on $D$,
we have $f_{\tau,\tau'} = f_{\tau,\tau''}$ on $D$.
\end{definition}

In words, the value of a $\tau$-function at
$(\tau, \tau')$ depends on $\tau'$ only through the information available
at time~$\tau$ about the value of $\tau'$: if $\tau'$ and $\tau''$ agree
on an $\calF_{\tau}$-measurable set, then $f_{\tau,\tau'}$ and $f_{\tau,\tau''}$ agree on that
set as well.

The following is the stochastic variation of Ramsey's theorem due to \cite{SS04} that we will need.

\begin{theorem}[\cite{SS04}, Theorem~4.3]
\label{thm:ramsey}
Let $(\Omega,\calF,\prob)$ be a probability space,
let $\mathbb{F}$ be a filtration,
and let $(n_k)_{k\in \dN}$
be an increasing sequence of integers.
Let $f$ be a $\tau$-function (for $(n_k)_{k\in \dN}$) with a finite range: 
There is a finite set $\calA$ such that the range of $f_{\tau,\tau'}$ is in $\calA$ a.s.
Then for every $\ep > 0$, there exists an increasing sequence of
$(n_k)_{k\in \dN}$-consistent
bounded stopping times $(\phi_k)_{k \in \dN}$ such that
$\prob(f_{\phi_0,\phi_1} = f_{\phi_1,\phi_2} = \cdots) > 1 - \ep$.
\end{theorem}

\cite{SS04} proved Theorem~\ref{thm:ramsey} in the special case that $n_k = k$ for every $k \in \dN$.
The statement of Theorem~\ref{thm:ramsey} reduces to that special case 
by considering the filtration $\mathbb{H} = (\calH_k)_{k \in \dN}$ defined by $\calH_k := \calF_{n_k}$ for every $k \in \dN$,
and applying the result of \cite{SS04} to the filtration $\mathbb H$.

We conclude this section by mentioning that,
when $\Omega$ is countable,
\cite{Xu11} provided a version of Theorem~\ref{thm:ramsey}
in which,
with probability at least $1-\ep$,
all values $(f_{\phi_i,\phi_j})_{0 \leq i < j}$ are equal.

%%%%%%%%%
\subsection{Partial games and \texorpdfstring{$\delta$}{delta}-approximations}
\label{sec:approximations}

In this section,
for every two bounded stopping times $\tau < \tau'$,
we study the stopping games restricted to the stages between $\tau$ and $\tau'$ (excluded).
We will then approximate,
for every $\omega \in \Omega$,
the part of the game induced between two stopping times $\tau < \tau'$ (excluded) at $\omega$ by a game on a tree.
We start by defining strategies restricted to the stages between $\tau$ and $\tau'$ (excluded).

\begin{definition}[$(\mathbb{F}, \tau, \tau')$-strategies]
\label{def:f:strategy}
Let $\Gamma =
(\Omega, I, \calF,\prob,{\mathbb F},X)$ be a stopping game,
let $\tau \leq \tau'$ be two bounded $\mathbb F$-adapted stopping times,
and let $i \in I$.
An $(\mathbb{F}, \tau, \tau')$-\emph{strategy} for player $i$ is a
sequence $x^i = (x^i_k)_{k \geq 0}$ of functions such that
$x^i_k : \{\tau \leq k < \tau'\} \to [0,1]$ is measurable with respect to the sigma-algebra $\calF_k$ restricted to $\{\tau \leq k < \tau'\}$,
for every $k \geq 0$.
\end{definition}

An $(\mathbb{F}, \tau, \tau')$-strategy prescribes to player~$i$
what to play at each stage $\tau \leq k < \tau'$.
Whenever necessary, we extend such strategies to strategies in $\Gamma$, by setting them to 0 whenever $k < \tau$ or $k \geq \tau'$.
For an $(\mathbb{F}, \tau, \tau')$-strategy profile $x=(x^i)_{i \in I}$,
define the conditional termination probability and expected payoffs:
\begin{align*}
  \pi(x; \mathbb{F}, \tau, \tau')
    &:= \Prob_{x}\bigl(\tau \leq \theta < \tau'
      \mid \calF_{\tau}\bigr), \\
  \rho^i(x; \mathbb{F}, \tau, \tau')
    &:= \E_{x}\bigl[\bfone_{\{\tau \leq \theta < \tau'\}}
      X^i_{J_\theta,\theta} \mid \calF_{\tau}\bigr], \\
  \gamma^i(x; \mathbb{F}, \tau, \tau')
    &:= \frac{\rho^i(x;\mathbb{F},\tau, \tau')}
      {\pi(x;\mathbb{F},\tau, \tau')},
\end{align*}
with the convention $0/0 = 0$.
These are $\calF_{\tau}$-measurable random variables.

We next define approximations of the stopping game induced between the stopping times $\tau < \tau'$ by a finite collection of stopping games on finite trees.

\begin{definition}[$\delta$-approximations]
\label{def:delta-approx}
Let $\tau < \tau'$ be two bounded stopping times and let
$\delta > 0$.
A \emph{$\delta$-approximation} of $\Gamma$ between $\tau$ and
$\tau'$ is a pair $\bigl((\mathcal{G}_k)_{k \in \dN}, (q_{G,k})_{\{k \in \dN, G \in \calG_k \colon k < \tau' \hbox{ on } G\}}\bigr)$,
where, for every $k \in \dN$,
\begin{enumerate}[(1)]
\item $\mathcal{G}_k$ is a finite partition of the event
  $\{\tau \leq k \leq \tau'\}$ into $\calF_k$-measurable sets, 
  each of positive measure,
  such that $\tau$ and $\tau'$ are adapted to
  $(\mathcal{G}_k)$: for every $k \geq 0$, the events $\{\tau = k\}$
  and $\{\tau' = k\}$ are unions of atoms of $\mathcal{G}_k$;
\item For every atom $G \in \mathcal{G}_k$
such that $k < \tau'$ on $G$, 
$q_{G,k}$ is a
  probability distribution over the atoms of $\mathcal{G}_{k+1}$
  contained in $G$;
\end{enumerate}
such that, for every
$k \in \dN$, on $\{\tau \leq k < \tau'\}$,
\begin{enumerate}[(1)]
\setcounter{enumi}{2}
\item The payoff $X_k$ is $\mathcal{G}_k$-measurable;
\item Every atom $G \in \mathcal{G}_k$ is a
  union of atoms of $\mathcal{G}_{k+1}$;
\item For every atom $G \in \mathcal{G}_k$ and almost every
  $\omega \in G$,
\[ \sum_{\{G' \in \mathcal{G}_{k+1} \colon G' \subseteq G\}} |\prob(G' \mid \calF_k)(\omega) -
  q_{G,k}(G')| < \delta/2^{k+2}. \]
\item
For every $k\in \dN$, every
$G\in\mathcal G_k$ with $\prob(G)>0$,
and every atom $G'\in\mathcal G_{k+1}$
contained in $G$
such that $q_{G,k}(G')>0$, 
we have $\prob(G')>0$.
\end{enumerate}
\end{definition}

The $\delta$-approximation approximates the stopping game $\Gamma$ between stages $\tau$ and $\tau'$ (excluded) in the following sense.
On the set $\tau \leq k < \tau'$,
on each atom $G$ of $\calG_k$,
\begin{itemize}
\item 
the payoff is constant on each atom of the partitions $(\calG_k)_{k \in \dN}$ (Condition (3) in Definition~\ref{def:delta-approx}),
and
\item 
the transition at $G$ is close
to the transition in $\Gamma$ at stage $k$ (Condition (5) in Definition~\ref{def:delta-approx}).
\end{itemize}
This leads us to define for every $\delta$-approximation a corresponding game on a tree.

\begin{notation}[$\calT(\tau, \tau'; (\mathcal{G}_k), (q_{G,k}), G)$]
\label{def:tree:approx}
To every $\delta$-approximation $\bigl((\mathcal{G}_k)_{k \in \dN}, (q_{G,k})_{k \in \dN, G \in \calG_k}\bigr)$ of
$\Gamma$ between $\tau$ and $\tau'$ and every%
\footnote{Here, $\calG_\tau$ is the partition that contains,
for each $k \in \dN$,
all sets $G \in \calG_k$ such that $\tau = k$ on $G$.
}
atom
$G \in \mathcal{G}_{\tau}$, define a game on a
finite tree $\calT(\tau, \tau'; (\mathcal{G}_k), (q_{G,k}), G)$
as follows:
\begin{itemize}
\item The root is $G$.
\item The nodes are all pairs $(k,F)$ such that $F \in \mathcal{G}_k$, $F \subseteq G$, and $\tau \leq k \leq \tau'$ on $F$.
That is, the nodes correspond to the atoms of the partitions $(\calG_k)_k$ that lie between $\tau$ and $\tau'$ (including),
where an atom that appears in several of the partitions corresponds to several nodes.
\item The leaves are those nodes $(k,F)$ such that
  $k = \tau'$ on $F$.
\item The payoff at a non-leaf node $(k,F)$ is
the unique vector $g$ such that 
$X_k(\omega) = g$ a.s.~on $F$. 
This vector is unique since $\prob(F) > 0$.
\item The children of $(k,F)$ are all nodes
  $(k+1,F')$ with $F' \subseteq F$.
\item The transition from $(k,F)$ is $q_{F,k}$.
\end{itemize}
\end{notation}

%For every filtration $\mathbb G = (\calG_k)_{k \in \dN}$ 
%and every adapted stopping time $\tau$,
%denote by $\calG_\tau$ the partition that contains,
%for each $k \in \dN$,
%all sets $G \in \calG_k$ such that $\tau = k$ on $G$.

%%%
\subsection{A $\tau$-function of approximating trees}
\label{sec:nt:approx}

Let $\tau < \tau'$ be two bounded stopping times,
and suppose that $\bigl(\mathcal{G}_k, (q_{G,k})_{G \in \calG_k}\bigr)_{k \in \dN}$
is a $\delta$-approximation of $\Gamma$ between stages $\tau$ and $\tau'$.
Let $x^i=(x_t^i)_{t\ge 0}$ be an $(\mathbb{F},\tau,\tau')$-strategy 
for player~$i$
such that, for every $t \in \dN$,
the random variable $x_t^i$ is almost surely constant on each atom of $\calG_k$.

For 
each $G \in \calG_\tau$,
the strategy $x^i$ induces a strategy $\alpha^i$
in the game on a tree
$\calT(\tau, \tau'; (\mathcal{G}_k), (q_{G,k}), G)$,
which is defined as follows:
For every non-leaf node $(k,F)$,
$\alpha^i(F)$
is the mixed action that $x^i_k$ coincides with a.s.~on $F$.

Conversely, any collection of strategies $\bigl(\alpha^i(\,\cdot\,;G)\bigr)_{G \in \calG_\tau}$ in the games on finite trees\\ 
$(\calT(\tau, \tau'; (\mathcal{G}_k), (q_{G,k}), G))_{G \in \calG_\tau}$ induces an $(\mathbb{F},\tau,\tau')$-strategy by setting,
\begin{equation}
\label{equ:defin:e}
x_t^i(\omega):=\alpha^i(F;G),
\qquad\text{whenever } \omega\in F\in\mathcal G_t \text{ and } \tau(\omega)\le t<\tau'(\omega).
\end{equation}
Therefore, an $(\mathbb{F},\tau,\tau')$-strategy that is constant on the approximating atoms is equivalent to a family of strategies on a finite tree.
%that is measurable and constant on the atoms of the approximating partitions.
Abusing notation, we will treat the collection $\bigl(\alpha^i(\,\cdot\,;G)\bigr)_{G \in \calG_\tau}$ as an $(\mathbb{F},\tau,\tau')$-strategy.

Under this identification, the quantities
$\pi(\cdot;\mathbb F,\tau,\tau')$ and $\rho^i(\cdot;\mathbb F,\tau,\tau')$
are approximated by the corresponding quantities in the games on a tree:
For every collection of strategy profiles $(\alpha(\omega))_{\omega \in \Omega}$ in the games on a tree $(\calT_{\tau,\tau'}(\omega))_{\omega \in \Omega}$ that are constant on atoms of the partitions $(\calG_k)_k$,
  and each $i \in I$,
  \begin{equation}
  \label{equ:app:1}
    \bigl|\rho^i(\alpha; \mathbb{F}, \tau, \tau')
      - \rho^i(\alpha; \calT_{\tau,\tau'})\bigr|
    < \delta/2^{\tau+1}, \hbox{  a.s.}
  \end{equation}
and
\begin{equation}
\label{equ:app:2}
    \bigl|\pi(\alpha; \mathbb{F}, \tau, \tau')
      - \pi(\alpha; \calT_{\tau,\tau'})\bigr|
    < \delta/2^{\tau+1}, \hbox{  a.s.}.
  \end{equation}

We now state an approximation theorem, which is the
multiplayer version of \cite[Corollary~6.6]{SS04}.
This result states that the stopping game between stages $\tau$ and $\tau'$ (excluded)
can be $\delta$-approximated by a game on a tree.
Moreover, for any $\ep$-equilibrium of the approximating game on a tree,
no player can profit much in the original game $\Gamma$ by deviating.
 
\begin{theorem}
\label{cor:approx}
Let $\Gamma$ be a multiplayer stopping game
whose payoff process has finite range and
is bounded by 1.
let $(n_k)_{k \in \dN}$ be an increasing sequence of non-negative integers,
let $(\eta_k)_{k\in \dN}$ be any sequence of positive real numbers,
and let $\delta,\ep > 0$ be sufficiently small.
There is a consistent $\tau$-function $(\calT_{\tau,\tau'})$ that
assigns to every two bounded $(n_k)_{k \in \dN}$-consistent stopping times $\tau < \tau'$ and
every $\omega \in \Omega$
a game on a tree
$\calT_{\tau,\tau'}(\omega)$,
such that the following properties hold for every two bounded $(n_k)_{k \in \dN}$-consistent stopping times $\tau < \tau'$:
\begin{enumerate}[(1)]
\item 
There is a $\delta$-approximation $\bigl(\mathcal{G}_k, (q_{G,k})_{G \in \calG_k}\bigr)_{k \in \dN}$ 
of $\Gamma$ between stages $\tau$ and $\tau'$, such that $\calT_{\tau,\tau'}(\omega) = \calT(\tau, \tau'; (\mathcal{G}_k), (q_{G,k}), G)$ 
a.s.
%is constant on each atom of $\calG_\tau$.

%where $\eta_\tau(\omega) = \eta_{\tau(\omega)}$ for every $\omega \in \Omega$.
\item Let $(\alpha^*(\omega))_{\omega \in \Omega}$ be $\ep$-equilibria of $(\calT_{\tau,\tau'}(\omega))_{\omega \in \Omega}$ that are constant on the atoms of the partition $\calG_\tau$,
and suppose that for every $\omega \in \Omega$
there is $a(\omega) = (a_i(\omega))_{i\in I}\in [-1,1]^I$,
such that
the expected payoff under $\alpha^*(\omega)$ lies in the box
  $\prod_{i \in I} [a^i, a^i + \ep]$. 
  Then for every $i\in I$ and every $(\mathbb{F}, \tau,\tau')$-strategy $x^i$ of player $i$,
    \[
      \rho^i(x^i,\alpha^{*,-i}; \mathbb{F}, \tau,\tau')
      \leq (a^i + 2\ep)\pi(x^i,\alpha^{*,-i}; \mathbb{F}, \tau,\tau')
        + \eta_\tau, \hbox{  a.s.}
    \]
%\item More generally, let $(\alpha^*(\omega))_{\omega\in\Omega}$
%be stationary profiles of $(\calT_{\tau,\tau'}(\omega))_{\omega \in \Omega}$, constant on atoms
%of $\calG_\tau$, and let $b(\omega)\in[-1,1]^n$ be
%$\calG_\tau$-measurable.  Suppose that, for every player $i$ and every
%stationary strategy $z^i$ in the corresponding tree,
%\[
% \rho^i(z^i,\alpha^{*,-i};\calT_{\tau,\tau'})
% \le
% b^i\,
% \pi(z^i,\alpha^{*,-i};\calT_{\tau,\tau'}).
%\]
%Then, for every $(\mathbb F,\tau,\tau')$-strategy $x^i$ of player $i$,
%\begin{equation}
%\label{eq:approx-general}
% \rho^i(x^i,\alpha^{*,-i};\mathbb F,\tau,\tau')
% \le
% b^i\,
% \pi(x^i,\alpha^{*,-i};\mathbb F,\tau,\tau')
% +\eta_\tau,
% \qquad\text{a.s.}
%\end{equation}
\end{enumerate}
\end{theorem}

%The additional assertion in part~(3) is proved by the same backward
%induction as \cite[Lemma~6.4]{SS04}.  At each node one compares the
%best continuation value under the full filtration with the best value
%under the finite approximation.  The argument is player-by-player and
%does not depend on the number of opponents.

The construction in \cite[Section~6]{SS04} is purely
measure-theoretic and involves only the filtration and the payoff
process.
It therefore extends without substantive changes to the multiplayer case, delivering Theorem~\ref{cor:approx}.

We note that in their definition of $\delta$-approximation, \cite{SS04} did not require property (6).
However, since $\range(X)$ is finite, 
their construction can be adapted to satisfy this property.

%%%%%%%%%%%%%%%%%%%
\section{Proof of Theorem~\ref{thm:main}}
\label{sec:proof}

In this section, we prove Theorem~\ref{thm:main}.
In Section~\ref{sec:reduction}, we show that without loss of generality we can assume that the payoff process is bounded (Section~\ref{sec:reduction:bounded}) and has a finite range (Section~\ref{section:finite:range}).
Moreover, we show that without loss of generality we can assume that the payoffs are generic and a multiple of a fixed real (Section~\ref{section:reduction:unique}).
We also show that without loss of generality we can assume that any payoff vector that appears infinitely often along the play with positive probability,
appears infinitely often with probability 1 (Section~\ref{section:ii:payoff}).
The assumptions on the game after these reductions are summarized in Section~\ref{sec:reduction:summary}.

We then turn to construct an approximate equilibrium.
The construction is divided according to properties of the set of vector payoffs that appear infinitely often.

In Section~\ref{sec:simple} we identify a condition that ensures the existence of an approximate equilibrium where 
stopping is done mainly by a single player.
In Section~\ref{sec:Acase}, we study the case when the 
matrix $R^{\mathrm{sp}}$ is $\ep$-close to some Q-matrix, and
in Section~\ref{sec:Bcase} we do the same when this matrix is not $\ep$-close to any Q-matrix.
Throughout we will assume that $\ep \in (0,1/36^2)$ is sufficiently small.

%% =====================================================================
\subsection{Reduction Steps}
\label{sec:reduction}
%% =====================================================================

\subsubsection{Reduction to bounded payoffs}
\label{sec:reduction:bounded}

In this subsection, we show that it is enough to prove Theorem~\ref{thm:main} when the payoff process is bounded.

Let $\Gamma = (\Omega,I,\calF,\prob,\mathbb{F}, X)$ be a stopping game that satisfies the conditions of Theorem~\ref{thm:main},
and let $\ep> 0$.
Define the random variable $W$ by 
\[ W(\omega) := \max_{i \in I, \varnothing \neq J\subseteq I}\left(\sup_{t \in \dN}  |X^{i}_{J,t}(\omega)|\right),
\ \ \ \forall \omega \in \Omega. \]
Since the payoff process is integrable,
$W$ is integrable as well.
Hence,
there exists $M < \infty$ such that
$\E[\mathbf{1}_{W > M} W] \leq \ep$.

Consider an auxiliary stopping game $\widetilde\Gamma = (\Omega,I,\calF,\prob,\mathbb{F}, \widetilde X)$
which is similar to $\Gamma$,
except for the payoff process $\widetilde X$ which is given by
\[ \widetilde X^{i}_{J,t} := \left\{
\begin{array}{ll}
M, & \hbox{if } X^{i}_{J,t} > M,\\
X^{i}_{J,t}, & \hbox{if } -M \leq X^i_{J,t} \leq M,\\
-M, & \hbox{if } X^{i}_{J,t} < -M.
\end{array}
\right.
\]

The game $\widetilde \Gamma$ satisfies the conditions of Theorem~\ref{thm:main}
and its payoff process is bounded by $M$.
Denote by $\widetilde \gamma^i(x)$ the expected payoff of player~$i$ in $\widetilde \Gamma$ under strategy profile $x$.
Then, for every strategy profile $x$,
\[ |\gamma^i(x) - \widetilde \gamma^i(x)| \leq 
\E_x[\mathbf{1}_{W > M} W] \leq \ep, 
\ \ \ \forall i \in I.\]
Hence, every $\ep$-equilibrium of $\widetilde \Gamma$ is a $3\ep$-equilibrium of $\Gamma$.
Replacing $X$ by $\widetilde X$, we may therefore assume without loss of generality that the payoff process is uniformly bounded.

Since the payoffs are bounded, 
by multiplying the payoff function by a constant,
we can further assume that the payoffs are bounded by 1.

\subsubsection{Reduction to finite range}
\label{section:finite:range}

In this subsection, we argue that it is sufficient to prove Theorem~\ref{thm:main} when the payoff process takes values that are all multiples of a fixed real positive number.

Let $\Gamma = (\Omega,I,\calF,\prob,\mathbb{F}, X)$ be a stopping game with a conditionally atomless filtration and a payoff process
that is bounded by 1,
and let $\ep > 0$.
Let $\widetilde \Gamma = (\Omega,I,\calF,\prob,\mathbb{F}, \widetilde X)$ be the auxiliary stopping game 
which is similar to $\Gamma$,
except for the payoff process $\widetilde X$,
where payoffs are rounded down to the closest integer multiple of $\ep$.
That is,
\[ \widetilde X^i_{J,t}(\omega) = \ep \times \left\lfloor \frac{X^i_{J,t}(\omega)}{\ep} \right\rfloor, \ \ \ \forall i \in I, \forall \varnothing \neq J \subseteq I, \forall t \in \dN, \forall \omega \in \Omega. \]
Then 
$\prob(\|X -\widetilde X\|_\infty \leq \ep) = 1$.
By Lemma~\ref{lemma:simple}(O1),
every $\ep$-equilibrium of $\widetilde\Gamma$ is a $7\ep$-equilibrium of $\Gamma$.%
\footnote{In fact, every $\ep$-equilibrium of $\widetilde\Gamma$ is a $2\ep$-equilibrium of $\Gamma$.}
Replacing $X$ by $\widetilde X$, we may therefore assume without loss of generality that the payoff process has finite range.

\subsubsection{Reduction to generic payoffs that are an integer multiple of a fixed real}
\label{section:reduction:unique}

Let $\ep > 0$,
and let $\Gamma = (\Omega,I,\calF,\prob,\mathbb{F}, X)$ be a stopping game with a conditionally atomless filtration and a payoff process
that is bounded by 1 and where all payoffs are an integer multiple of a fixed real number.
We will show that without loss of generality we can assume that if $g^i_{J} = g'^i_{J'}$ for some $i \in I$, $\varnothing \subset J,J' \subseteq I$, and $g,g' \in \range(X)$, then $g = g'$.

Denote by 
$\range(X) \subseteq \dR^{n(2^n-1)}$ the finite range of $X$:
\[ \range(X) := \bigl\{ g \in \dR^{n(2^n-1)} \colon 
\prob(X_{t} = g) > 0 \hbox{ for some } t \in \dN\bigr\}. \]

%We argue that 
%after an approximation argument, it suffices to consider stopping games that satisfy the following two properties:
%\begin{itemize}
%\item 
%All payoffs are integer multiples of some fixed real number $\eta > 0$.
%\item
%For every $g,g' \in \range(X)$,
%either $g = g'$,
%or $g^i_J \neq g'^i_{J'}$ for every $i \in I$ and $\varnothing \neq J,J' \subseteq I$.
%\end{itemize}
Since $\range(X)$ is finite,
there is $\eta > 0$, and for every $g \in \range(X)$ there is $\beta(g) \in (-1,1)^{n(2^n-1)}$ such that 
the collection $(\beta(g))_{g \in \range(X)}$ satisfies the following:
\begin{itemize}
\item $\beta(g)^i_J$ is an integer multiple of $\eta$, for every $i \in I$ and every $\varnothing \subset J \subseteq I$.
\item $\|g - \beta(g)\|_\infty \leq \ep$ for every $g \in \range(X)$,
and
\item $\beta(g)^i_J \neq \beta(g')^i_{J'}$,
for every $i \in I$, 
every $\varnothing \neq J,J' \subseteq I$,
and every $g' \in \range(X)$
such that $J \neq J'$ or $g \neq g'$.
\end{itemize} 
Let $\widetilde \Gamma = (\Omega,I,\calF,\prob,\mathbb{F}, \widetilde X)$ be a stopping game that is similar to $\Gamma$, except for the payoff process $\widetilde X$,
that is defined by
\[ \widetilde X_{t}(\omega) := \beta(X_{t}(\omega)),
\ \ \ 
\forall t \in \dN.
\]
Then $\prob(\|X -\widetilde X\|_\infty \leq \ep) = 1$.
By Lemma~\ref{lemma:simple}(O1),
every $\ep$-equilibrium of $\widetilde \Gamma$ is a $7\ep$-equilibrium of $\Gamma$.
Finally, the game $\widetilde\Gamma$ satisfies the property that 
for every $g,g' \in \range(X)$,
either $g = g'$,
or $g^i_J \neq g'^i_{J'}$ for every $i \in I$ and $\varnothing \neq J,J' \subseteq I$.
Replacing $X$ by $\widetilde X$, we may therefore assume without loss of generality that the game $\Gamma$ satisfies this property.

\subsubsection{Reduction to the case where all payoffs appear i.o.~with probability 1}
\label{section:ii:payoff}

Let $\Gamma = (\Omega,I,\calF,\prob,\mathbb{F}, X)$ be a stopping game with a conditionally atomless filtration and a payoff process
that is bounded by 1,
has a finite range,
and in addition satisfies the genericity assumption of Section~\ref{section:reduction:unique}.
Let $\ep > 0$.

The set of all vectors in the range of $X$ that occur i.o.~is
\[ \calH := \{ g \in \range(X) \colon X_t = g \hbox{ i.o.}\}. \]

We will need the following technical result,
whose proof appears at the end of this subsection.

\begin{lemma}
\label{lemma:approx:H}
Let $\calD$ be a finite $(\vee_{t \in \dN} \calF_t)$-measurable partition of $\Omega$ such that $\calH$ is constant on each element of $\calD$.
There are an $\mathbb F$-adapted process $(\widetilde X_t)_{t \in \dN}$
such that 
\[ \prob(X_t \neq \widetilde X_t \hbox{ for some } t \in \dN) \leq \ep, \]
a non-negative integer $t_0 \in \dN$,
and a finite $\calF_{t_0}$-measurable partition $\widetilde\calD$,
such that $\widetilde \calH$ is constant on each element of $\widetilde\calD$,
where 
\[ \widetilde\calH := \{ g \in \range(\widetilde X) \colon \widetilde X_t = g \hbox{ i.o.}\}. \]
\end{lemma}

Lemma~\ref{lemma:approx:H} states that by changing $X$ on a small set,
we can assume that $\calD$ is $\calF_{t_0}$-measurable, for some $t_0 \in \dN$.

Let then $\widetilde X$, $t_0$, and $\widetilde D$ be as in the conclusion of Lemma~\ref{lemma:approx:H}.
By Lemma~\ref{lemma:simple}(O2,O3),
if for every $D \in \widetilde\calD$ the game 
$\Gamma_D = (\Omega,I,\calF,\prob(\cdot \mid D),\mathbb{F}, \widetilde X)$ that starts at stage $t_0$
admits an $\ep$-equilibrium,
then the game $(\Omega,I,\calF,\prob,\mathbb{F}, \widetilde X)$ admits an $\ep$-equilibrium.
By Lemma~\ref{lemma:simple}(O1),
in this case the game $\Gamma$ admits a $7\ep$-equilibrium.

It is therefore sufficient to prove Theorem~\ref{thm:main} when,
in addition to the conditions at the beginning of this section,
there is a set $H \subseteq \range(X)$ such that $\calH = H$ a.s.

\bigskip

Some vectors in $\range(X)$ may appear along $X$ finitely many times.
We next show that we can assume that this does not happen.
Increase $t_0$ if necessary so that
\[ \prob(X_t \not\in H \hbox{ for some } t \geq t_0) < \ep. \]
Let $\widetilde \Gamma = (\Omega,I,\widetilde\calF,\prob,\widetilde{\mathbb{F}}, \widetilde X)$ be the auxiliary stopping game 
which is similar to $\Gamma$ when that game starts at stage $t_0$,
except for the payoff process $\widetilde X$,
where payoffs not in $H$ are replaced by some payoff in $H$: select an arbitrary $g^* \in H$, and define for every $t \in \dN$,
\begin{align*}
\widetilde\calF_t &:= \calF_{t_0 + t},\\
\widetilde X_t &:= \left\{
\begin{array}{ll}
g^*, & \hbox{if } X_{t_0+t} \not\in H,\\
X_{t_0+t}, & \hbox{if } X_{t_0+t} \in H.
\end{array}
\right.
\end{align*}
The process $\widetilde X$ is $\widetilde{\mathbb{F}}$-adapted.
The game $\widetilde\Gamma$ satisfies that the range of $\widetilde X$ is $H$,
and all vectors in $H$ appear along $X$ i.o.~with probability 1.
By Lemma~\ref{lemma:simple}(O1),
every $\ep$-equilibrium in $\widetilde\Gamma$
is a $7\ep$-equilibrium in the subgame of $\Gamma$ that starts at stage $t_0$.
By Lemma~\ref{lemma:simple}(O2),
such a $7\ep$-equilibrium can be extended to a $7\ep$-equilibrium of $\Gamma$.

\bigskip

\begin{proof}[Proof of Lemma~\ref{lemma:approx:H}]
Write $\calD=\{D^1,\ldots,D^L\}$ and
$K:=|\range(X)|$.
We may assume that $\prob(D^\ell)>0$ for every $\ell \in \{1,\dots,L\}$.
Denote by $H^\ell$ the value of $\calH$ on $D^\ell$.
Since $X$ has finite range, each $H^\ell$ is nonempty and contains at most $K$ values.

For every $n\in \dN$, choose $t_n \in \dN$
and an $\calF_{t_n}$-measurable partition
$\widetilde\calD_n
=\{\widetilde D_n^1,\ldots,\widetilde D_n^L\}$
such that
\[
\prob(D^\ell\triangle\widetilde D_n^\ell)
<\frac{\ep}{L2^{n+2}},
\qquad \forall \ell=1,\ldots,L.
\]
Choose a strictly increasing sequence of non-negative integers $(T_n)_{n\in \dN}$ such that
\[
T_n\geq t_n,
\qquad T_{n+1}-T_n\geq K,
\qquad \forall n \in \dN.
\]
Set $t_0:=T_0$ and $\widetilde\calD:=\widetilde\calD_0$.
For each $\ell \in \{1,\dots,L\}$, fix an ordering
\[
H^\ell=\{h_1^\ell,\ldots,h_{K_\ell}^\ell\},
\qquad \hbox{where } K_\ell:=|H^\ell|\leq K.
\]
Define $\widetilde X_t:=X_t$ for $t<T_0$.
For each $n\in \dN$ and
$T_n\leq t<T_{n+1}$, define $\widetilde X_t$
as follows.
On $\widetilde D_0^\ell\cap\widetilde D_n^\ell$,
set
$\widetilde X_t:=X_t$.
On $\widetilde D_0^\ell\setminus\widetilde D_n^\ell$,
set
$\widetilde X_t
:=h_{1+((t-T_n)\bmod K_\ell)}^\ell$.

Since $\widetilde\calD_0$ and
$\widetilde\calD_n$ are $\calF_{T_n}$-measurable,
the process $(\widetilde X_t)_{t \in \dN}$ is $\mathbb F$-adapted.
Moreover, 
\[ \bigl\{ w\in \Omega \colon \widetilde X_t \neq X_t \hbox{ for some } t \in \dN\bigr\} \subseteq 
\bigcup_{1 \leq \ell \leq L, n \in \dN} (\widetilde D^\ell_0 \setminus \widetilde D^\ell_n) 
\subseteq \bigcup_{1 \leq \ell \leq L, n \in \dN} (D^\ell \triangle \widetilde D^\ell_n) \]
Hence,
\begin{align*}
&\prob(\widetilde X_t\neq X_t
       \text{ for some }t\in\dN)
\leq
\sum_{n\in \dN}\sum_{\ell=1}^L
\prob(D^\ell\triangle\widetilde D_n^\ell)
<\ep.
\end{align*}
It remains to prove that $\widetilde\calH$ is constant
on every element of $\widetilde\calD$.

By the Borel--Cantelli lemma, almost surely there is
$n_0$ such that $\mathbf1_{\widetilde D_n^\ell}
=\mathbf1_{D^\ell}$
for every $n\geq n_0$ and every
$\ell=1,\ldots,L$.
Fix $\ell$ and consider
$\omega\in\widetilde D_0^\ell$ for which this
property holds.

If $\omega\in D^\ell$, then
$\omega\in\widetilde D_n^\ell$ for all sufficiently
large $n$. Hence $\widetilde X_t=X_t$ for all
sufficiently large $t$, and therefore
$\widetilde\calH(\omega)
=\calH(\omega)=H^\ell$.
If $\omega\notin D^\ell$, then
$\omega\notin\widetilde D_n^\ell$ for all
sufficiently large $n$.
Hence, for every $n \in \dN$ sufficiently large, the process $\widetilde X$
cycles through all the elements of $H^\ell$
during $T_n \leq t < T_{n+1}$.
Therefore,
$\widetilde\calH(\omega)=H^\ell$.
We conclude that
$\widetilde\calH=H^\ell$ a.s.~on $\widetilde D_0^\ell$,
for every $\ell=1,\ldots,L$, and the claim follows.
\end{proof}

\subsubsection{Summary of the reductions}
\label{sec:reduction:summary}

From now on,
we fix a stopping game $\Gamma = (\Omega,I,\calF,\prob,\mathbb{F}, X)$ that satisfies the following conditions:
\begin{itemize}
\item[AS1)] The filtration $\mathbb F$ is conditionally atomless.
\item[AS2)] Payoffs are bounded by 1: 
$\|g\|_\infty < 1$ for every $g \in \range(X)$ (Section~\ref{sec:reduction:bounded}).
\item[AS3)] Payoff range is finite:
$\range(X)$ is finite, and all payoffs are integer multiples of a given real number (Sections~\ref{section:finite:range} and~\ref{section:reduction:unique}).
\item[AS4)] 
For every $g,g' \in \range(X)$,
either $g = g'$,
or $g^i_J \neq g'^i_{J'}$ for every $i \in I$ and $\varnothing \neq J,J' \subseteq I$ 
(Section~\ref{section:reduction:unique}).
\item[AS5)] 
All payoff vectors in the range of $X$ 
appear i.o.~with probability 1
(Section~\ref{section:ii:payoff}).
\end{itemize}
We also fix $\ep>0$ sufficiently small such that all payoffs in the
game are integer multiples of $\ep$ and all the inequalities required
in Section~\ref{sec:Bcase} and in the proof of
Proposition~\ref{prop:absorbing-eq} hold.
We will prove that $\Gamma$ admits a $O(\sqrt\ep)$-equilibrium.

\subsubsection{Further notations: $y^i$, $Z_i$, $R$, $R^{\mathrm{sp}}$, $I^{\mathrm{sp}}$}

For each player $i \in I$, 
the maximal payoff player~$i$ can attain by stopping alone is
\[ 
\label{def:y:stopping}
y^i = \max\bigl\{ g^i_{\{i\}} \colon g \in \range(X)\bigr\}. \]

The reduction in Section~\ref{section:reduction:unique} implies that there is a single vector $g \in \range(X)$ such that $g_{\{i\}}^i = y^i$;
set 
$Z_i = g_{\{i\}} \in \dR^I$.
\label{def:z:stopping}
Denote by $R$ the $n \times n$ matrix 
whose $i$'th column is $Z_i$ (viewed as a column vector).
Note that $R$ need not be derived from a single element in $\range(X)$,
since for $i \neq j$, $Z_i$ and $Z_j$ might be realized by different elements in $\range(X)$. 

Denote by $I^{\mathrm{sp}}$ the set of supporters for $R$.
As in Lemma~\ref{lemma:abnormal},
if $I^{\mathrm{sp}} = \emptyset$, then 
the strategy profile in which no players ever stops is 
a 0-equilibrium of $\Gamma$.
We will therefore assume from now on that $I^{\mathrm{sp}} \neq \emptyset$.
Denote by $R^{\mathrm{sp}}$ the principal submatrix of $R$ indexed by $I^{\mathrm{sp}}$.

%=====================================================================
\subsection{$\ep$-equilibria in simple cases}
\label{sec:simple}

In this section, we identify conditions that ensure the existence of an approximate equilibrium where stopping is done mainly by a single player.

\begin{lemma}
\label{lemma:simple:eq}
Suppose that there is a player $i_0 \in I$ and a real number
$\beta \in (0,1]$ such that the following properties hold,
where $g$ is the unique vector in $\range(X)$ that satisfies $g_{\{i_0\}} = Z_{i_0}$:
\begin{itemize}
\item[i)] $y^{i_0} \geq 0$ or there is $\widetilde g \in \range(X)$ and $i_1 \neq i_0$ such that $\widetilde g_{\{i_1\}}^{i_0} \leq y^{i_0}$.
\item[ii)]
$Z^i_{i_0} \geq (1-\beta) g^i_{\{i\}} + \beta g^i_{\{i,i_0\}}$,
for every $i \neq i_0$.
\item[iii)]
For every $i \neq i_0$ and $\widetilde g \in \range(X) \setminus \{g\}$
we have
$Z^i_{i_0} \geq \widetilde g^i_{\{i\}}$.
\end{itemize}
% Suppose there is a supporter $i \in I^{\mathrm{sp}}$ such that
% (a) $Z_i^j \geq Z_j^j$ for every $j \neq i$,
% and (b) there exists $\beta \in [0,1]$
% such that for the unique $g \in H$ with $g_{\{i\}} = Z_i$ we have
% $Z^j_i \geq (1-\beta) g^j_{\{j\}} + \beta g^j_{\{i,j\}}$ for every $j \neq i$.
Then an $\ep$-equilibrium $x^*$ exists for every $\ep > 0$.
\end{lemma}

Roughly, 
in the approximate equilibrium that we will construct,
player~$i_0$ will stop with probability $\beta$ whenever the payoff vector is $g$.
Condition (i) ensures that players $I \setminus \{i_0\}$ have a way to punish player $i_0$ if she does not stop.
Condition (ii) (resp., Condition (iii)) implies that each player $i \neq i_0$ cannot profit by deviating when the vector payoff is $g$ (resp. different from $g$).

Note that Condition~(iii) is strictly weaker than $Z_{i_0}^i \geq y^i$ when the maximal payoff that players $i$ and $i_0$ can obtain by stopping alone is attained at the same vector $g \in \range(X)$.

\bigskip
\begin{proof}[Proof of Lemma~\ref{lemma:simple:eq}]
The proof is standard, and hence we only provide a sketch.
Fix a player $i_0 \in I$ and $\beta \in (0,1]$ for which properties (i)--(iii) hold. 

\bigskip
\noindent\textbf{Case 1:
$y^{i_0} \geq 0$.}

Let $x^*$ be the following strategy profile:
\begin{itemize}
\item $x^{*,i_0}_t = \beta$ for every $t$ such that $X_{\{i_0\},t} = Z_{i_0}$,
and $x^{*,i_0}_t = 0$ otherwise.
\item $x^{*,i}_t = 0$ for every $i \neq i_0$ and every $t \in \dN$.
\end{itemize}
That is, under $x^*$,
player $i_0$ stops with probability $\beta$ whenever by stopping alone she gets her maximal payoff $y^{i_0}$,
and all other players never stop.

Since $\beta > 0$,
under $x^*$ the play terminates with probability 1
and the terminal payoff is $Z_{i_0}$.
Since $y^{i_0} \geq 0$, player~$i_0$ cannot profit by never stopping,
and since $Z_{i_0}^{i_0} = y^{i_0}$ is player~$i_0$'s maximal payoff by stopping alone,
she cannot profit by any deviation.
Condition~(ii) guarantees that no player $i \neq i_0$ can profit by stopping when $X_t = g$,
and Condition~(iii) guarantees that no player $i \neq i_0$ can profit by stopping when $X_t \neq g$.
Hence $x^*$ is a 0-equilibrium.

\bigskip
\noindent\textbf{Case 2:
$y^{i_0} < 0$ and
there is $\widetilde g \in \range(X)$ and $i \neq i_0$ such that $\widetilde g_{\{i\}}^{i_0} \leq y^{i_0}$.}

Let $\widehat x^*$ be the strategy profile where 
\begin{itemize}
\item 
Players $I \setminus \{i_1\}$ follow their strategies defined in Case~1.
\item
For every $t \in \dN$ such that $X_{t} = g$,
player~$i_1$ 
stops with probability $\ep\beta$ at the first stage $t' > t$ such that $X_{t'} = \widetilde g$.
In all other stages, Player~$i_1$ stops with probability 0.
\end{itemize}

Player~$i_1$ stops with probability $\ep\beta$ (at most) once for each time player~$i_0$ stops with probability $\beta$.
Hence, under $\widehat x^*$ the conditional probability that the play terminates by player~$i_0$ alone is at least $1-\ep$.
Therefore, $\|\gamma(\widehat x^*) - Z_{i_0}\|_\infty \leq 2\ep$.

Moreover, player~$i_1$ stops with probability $\ep\beta$ infinitely often. 
Hence, for any strategy $x^{i_0}$ of player~$i_0$
we have $\gamma^{i_0}(x^{i_0},\widehat x^{*,-i_0}) \leq y^{i_0} + 2\ep\beta$.
Thus, player~$i_0$ cannot profit more than $4\ep$ by deviating from $\widehat x^*$.

Any player $i \neq \{i_0,i_1\}$ can profit by stopping at a stage in which player $i_1$ stops.
Since in such a stage player~$i_1$ stops with probability $\ep\beta$, 
the profit of such a deviation is bounded by $2\ep$.
Thus, property (iii) ensures that
no player $i \neq \{i_0,i_1\}$ can profit more than $4\ep$ by deviating from $\widehat x^*$ as well.

By never stopping, player~$i_1$ can profit at most $2\ep$,
and by stopping in a stage where $X_t \neq g$ she can profit at most $4\ep$.
Thus, $\widehat x^*$ is a $4\ep$-equilibrium.
\end{proof}

\bigskip

\begin{lemma}
\label{lemma:simple:2}
Suppose there is a player $i_0 \in I$,
a vector $g \in \range(X)$,
a vector $(\beta^j)_{j \neq i_0} \in [0,1]^{n-1}$ such that not all $(\beta^j)_{j \neq i_0}$ are 0,
and $\ep > 0$,
such that the following conditions hold:
\begin{itemize}
\item[i)] $y^{i_0} + \ep \geq \max\{A, B\}$,
where
\[ A = \sum_{J \subseteq I \setminus \{i_0\}} g^{i_0}_{J \cup \{i_0\}} \prod_{j \in J} \beta_j \prod_{j \not\in J \cup \{i_0\}} (1-\beta_j), \]
and
\[ B = \frac{\sum_{\varnothing \neq J \subseteq I \setminus \{i_0\}} g^{i_0}_{J} \prod_{j \in J} \beta_j \prod_{j \not\in J \cup \{i_0\}} (1-\beta_j)}{\sum_{\varnothing \neq J \subseteq I \setminus \{i_0\}} \prod_{j \in J} \beta_j \prod_{j \not\in J \cup \{i_0\}} (1-\beta_j)}. \]
\item[ii)] 
$Z^i_{i_0} \geq Z^i_i$ for every $i \neq i_0$.
\end{itemize}
Then a $4\ep$-equilibrium exists.
\end{lemma}

The quantity $B$ represents the expected payoff conditional on termination at the current stage, 
when the current payoff vector is $g$,
player~$i_0$ stops with probability 0,
and each player $j \neq i_0$ stops with probability $\beta_j$.
The quantity $A$ is the analog expected payoff 
when 
player~$i_0$ stops with probability 1,
and each player $j \neq i_0$ stops with probability $\beta_j$.
Thus, 
Condition~\ref{lemma:simple:2}(i) requires that there is $g \in \range(X)$ such that players $I \setminus \{i_0\}$ have a way to ensure that if the play terminates when the payoff vector is $g$,
then player~$i_0$'s payoff would be low, whether she stops or not.

\bigskip

\begin{proof}
We will construct a $4\ep$-equilibrium under which player~$i_0$ terminates the game,
and if she does not, the other players punish player~$i_0$ by stopping with probabilities $(\beta_j)_{j \neq i_0}$ whenever the payoff vector is $g$.

Let $N \in \dN$ be sufficiently large such that $(1-\ep)^N < \ep/8$.
Let $T \in \dN$ be sufficiently large such that 
\[ \prob\bigl(\#\{ t \leq T \colon X_{\{i_0\},t} = Z_{i_0}\} \geq N\bigr) \geq 1-\ep/8. \]

Let $x^* = (x^{*,j})_{j \in I}$ be the following strategy profile:
\[ x^{*,i_0}_t = \left\{
\begin{array}{ll}
\ep, & X_{\{i_0\},t} = Z_{i_0},\\
0, & X_{\{i_0\},t} \neq Z_{i_0},
\end{array}
\right.
\qquad
 x^{*,j}_t = \left\{
\begin{array}{ll}
0, & t \leq T,\\
\beta_j, & X_t = g, t > T,\\
0, & X_t \neq g, t > T,
\end{array}
\right.
\qquad
\forall j \neq i_0.
\]

Under $x^*$,
with probability at least $1-\ep/4$ the play terminates by player~$i_0$,
and hence 
\[ \gamma^j(x^*) \geq Z_{i_0}^j - \ep/2,
\qquad \forall j \in I. \]
Condition (i) ensures that player~$i_0$ cannot profit more than $3\ep/2$ by deviating.
By Condition (ii) and since player~$i_0$ stops with probability at most $\ep$ at each stage,
no player $j \neq i_0$ can profit more than $3\ep$ by deviating.
Therefore, $x^*$ is a $3\ep$-equilibrium.
\end{proof}

%%%%%%%%%%%%%%%%%%
\subsection{Construction when $R^{\mathrm{sp}}$ is $\ep$-close to a Q-Matrix}
\label{sec:Acase}
%% =====================================================================

In this section, we assume that $R^{\mathrm{sp}}$ is $\ep$-close to a Q-matrix,
and we construct a $14\varepsilon$-equilibrium $x^*$.
The strategy profile $x^*$ will follow the sunspot equilibrium $\xi^*$ of
Section~\ref{sec:ss:Q}, with the public correlation device replaced by
randomness extracted from the conditionally atomless filtration.
We will need the following result,
which asserts that such a replacement is possible.

\begin{lemma}
\label{lem:simulate-history}
Let $(\eta_k)_{k\geq 0}$ be a sequence of positive numbers with
$\sum_{k\geq 0}\eta_k<\ep/2$.
For every $k \in \dN$ and every $h_k \in I^k$,
let $\nu(h_k) = (\nu_i(h_k))_{i \in I} \in \Delta(I)$.
There exists an increasing sequence of
integers
$0=n_0<m_1<n_1<m_2<\cdots$
and for every $k \geq 1$ and every $h_k \in I^k$
there exists an $\calF_{m_{k}}$-measurable event $C(h_k)$,
such that 
for every $k \in \dN$ and every $h_k \in I^k$,
\begin{enumerate}[(a)]
\item $C(\varnothing) = \Omega$,
where $\varnothing$ is the empty history.
\item The sets $(C(h_k,i))_{i \in I}$ are disjoint and 
$C(h_k) = \bigcup_{i \in I} C(h_k,i)$.
\item For every $A\in \calF_{n_k}$,
  \[
  \sum_{i \in I} \bigl|\Prob(A\cap C(h_k,i))
      -\Prob(A \cap C(h_k))\nu_i(h_k)\bigr|
  \leq \eta_k\prob(C(h_k)).
  \]
\item 
$
%  \Prob\left(
  \Prob\Bigl(\forall i \in I \hbox{ there is } m_{k+1} \leq t < n_{k+1} \hbox{ such that }
  X_{\{i\},t}=Z_i
%  \text{ for some } t \in [m_{k+1},n_{k+1}), \forall i \in I
%  \Bigm| \calF_{m_{k+1}}
\Bigr)
  \geq 1-\eta_k$.
  %\right) \geq 1-\eta_k .
\end{enumerate}
\end{lemma}

The probability distribution $\nu(h_k)$ will be the one that governs the public correlation device.
Condition (c) will imply that the events $(C(h_k,i))_{i \in I}$
can be used in lieu of the public signal at stage $k$:
These events are almost independent of the information at stage $n_k$,
and their probability is close to the probability that, according to the public signal, player~$i$ is the one who should stop with positive probability at stage $k$.
Condition (d) requires that with high probability,
there is a stage between $m_{k+1}$ and $n_{k+1}$ (excluded)
in which the terminal payoff when player~$\widehat i_k$ stops alone is $Z_i$,
so that $\xi^*$ can be implemented (with high probability).

\bigskip

\begin{proof}
Assume by induction that 
we already defined 
$0=n_0<m_1<n_1<m_2<\cdots < n_{k}$
and, for each $\ell \in \{1,\dots,k\}$
we defined an $\calF_{m_{\ell}}$-measurable family $(C(h_\ell))_{i \in I}$ that satisfy (a)--(d).

Fix for a moment $h_k\in I^k$ such that $\prob(C(h_k)) > 0$,
and restrict attention to the event $C(h_k)$
(so that the probability distribution is $\prob(\cdot \mid C(h_k))$).
Apply Lemma~\ref{lem:simulate} 
with $t_0=n_k$, $n=|I|$,  
%$\mu=(\Leb(\widehat i_k = i\mid h_k))_{i \in I}$, 
$\mu = \nu(h_k)$,
and $\delta=\eta_k/|I|^{k+1}$.
Since the number of sequences $h_k \in I^k$ is finite, there exists $m_{k+1}>n_k$ and, for each $h_k$, an
$\calF_{m_{k+1}}$-measurable partition
$(C(h_k,i))_{i\in I}$ satisfying (a)--(c).

By (AS5), $\prob\bigl( \left\{ X_{\{i\},t}=Z_i \text{ i.o.}\right\} \bigr) = 1$.
Hence, and since $I$ is finite, one can choose 
$n_{k+1}>m_{k+1}$ such that (d) holds for every $i\in I$.
\end{proof}

\bigskip
We turn to construct the $14\ep$-equilibrium.
Let $(\eta_k)_{k\geq 0}$ be a sequence of positive reals such that
$\sum_{k\geq 0}\eta_k \leq \ep/3$.

\bigskip
\noindent\textbf{Step~1: An appropriate sunspot $4\ep$-equilibrium in $\calQ_{R}$.}

Since $R^{\mathrm{sp}}$ is $\ep$-close to a Q-matrix,
Theorem~\ref{thm:SS20-main} and the discussion in
Sections~\ref{sec:ss:Q} and~\ref{sec:ss:Q:close} imply that 
on the probability space $[0,1]^\dN$ there exist
\begin{itemize}
\item an adapted process $(\widehat i_k)_{k \geq 0}$ with values in the set
  $I\cup \{\varnothing\}$,
\item and an adapted process $(\widehat\delta_k)_{k \geq 0}$ with values in $(0,{\ep}/{3}]$,
\end{itemize}
such that the strategy
profile $\xi^*$ in $\calQ^{\mathrm{sun}}_{R}$ given by
\[
\xi_{k}^{*,i}(\psi_0,\dots,\psi_k):=
\left\{
\begin{array}{ll}
\widehat\delta_k(\psi_0,\dots,\psi_k), & \text{if } i=\widehat i_k(\psi_0,\dots,\psi_k),\\
0, & \text{if } i\neq \widehat i_k(\psi_0,\dots,\psi_k),
\end{array}
\right.
\]
is a sunspot $4\ep$-equilibrium in the quitting game $\calQ$.
Moreover, for every $k \geq 0$ the random
variable $\widehat i_k$ depends on
$(\psi_0,\dots,\psi_k)$ only through
$(\widehat i_0(\psi_0),\dots,\widehat i_{k-1}(\psi_0,\dots,\psi_{k-1}),\psi_k)$,
and the random
variable $\widehat\delta_k$ depends on
$(\psi_0,\dots,\psi_k)$ only through
$(\widehat i_0(\psi_0),\dots,\widehat i_{k}(\psi_0,\dots,\psi_k))$.
Furthermore,
$\prob_\xi(\theta < \infty \mid \psi_0,\dots,\psi_k,\theta \geq k) = 1$ a.s., for every $k \in \dN$.

\bigskip
\noindent\textbf{Step~2: Defining events that mimic the public signals.}

Apply Lemma~\ref{lem:simulate-history}
with 
%$(\ep\eta_k)_{k \in \dN}$ (for $(\eta_k)_{k \in \dN}$) and 
$(\Leb(\widehat i_k = i \mid h_k))_{i \in I}$ (for $\nu(h_k)$)
for every $k \in \dN$ and every $h_k \in I^k$,
to obtain an increasing sequence of bounded
stopping times
$0=n_0<m_1<n_1<m_2<n_2<\cdots$
and events $(C(h_k,i))_{i \in I}$, 
for every $k \in \dN$ and every $h_k \in I^k$,
that satisfy the conclusion of the lemma.

\bigskip
\noindent\textbf{Step 3: We can assume without loss of generality that $\prob(X_{\{i\},t}=Z_{i} \text{ for some } m_{k+1} \leq t < n_{k+1}) = 1$, for every $i \in I$.}

For every $k \in \dN$ define
\[
E_k(i):=\left\{
X_{\{i\},t}=Z_{i} \text{ for some } m_{k+1} \leq t < n_{k+1}\}
\right\}.
\]
Then $E_k(i) \in \calF_{n_{k+1}-1}$.
By Lemma~\ref{lem:simulate-history}(d), $\prob\left(\bigcap_{i \in I} E_k(i)\right) \geq 1-\eta_k$.
Let $\widetilde \Gamma = (\Omega,I,\calF,\prob,\mathbb{F}, \widetilde X)$ be the auxiliary stopping game 
which is similar to $\Gamma$,
except for the payoff process $\widetilde X$,
which is set to $Z_i$ on $\Omega \setminus E_k(i)$ for $t=n_{k+1}-1$:
\[ \widetilde X_{J,t}(\omega) :=
\left\{
\begin{array}{ll}
Z_i, & \hbox{if } J = \{i\}, \omega \not\in E_k(i), t = n_{k+1}-1,\\
X_{J,t}(\omega), & \hbox{otherwise}.
\end{array}
\right.
\]
Since $E_k(i)$ is $\calF_{n_{k+1}-1}$-measurable,
the process $\widetilde X$ is $\mathbb F$-adapted. 
Moreover, $\prob\bigl(\bigcap_{i \in I, k \in \dN} E_k(i)\bigr) \geq 1-\sum_{k \in \dN} \eta_k \geq 1-\ep/3$.

The set $E_k(i)$ for the game $\widetilde\Gamma$ is $\Omega$.
By Lemma~\ref{lemma:simple}(O1),
every $\ep$-equilibrium of $\widetilde\Gamma$ is a $7\ep$-equilibrium of $\Gamma$.
Hence we can assume without loss of generality that $E_k(i) = \Omega$ for every $k \in \dN$ and every $i \in I$.

\bigskip
\noindent\textbf{Step~4: Definition of $x^*$.}

We next define a strategy $x^{*,i}$ in $\Gamma$
for each player $i \in I$.
\begin{itemize}
\item Between stages $n_k$ and $m_{k+1}$ (excluded) no player stops:
  \[ x^{*,i}_t := 0, \ \ \ \forall i \in I, \forall n_k \leq t < m_{k+1}. \]

\item 
Between stages $m_{k+1}$ and $n_{k+1}$ (excluded), player $i$ 
  stops with positive probability only on the event $C(h_k,i)$; she stops with probability $\widehat\delta_k(h_k,i)$ at the first stage $t$ between $m_{k+1}$ and $n_{k+1}$ (excluded) in which $X_{\{i\},t}=Z_{i}$.
  By Step~3, such a stage exists.
  Formally,
  \[ x^{*,i}_t := \left\{
  \begin{array}{ll}
  \widehat\delta_k(h_k,i), & 
  \hbox{on } C(h_k,i), \hbox{ whenever } t = \min\bigl\{m_{k+1} \leq t' < n_{k+1} \colon X_  {\{i\},t'}=Z_{i}\bigr\},\\
  0, & \hbox{otherwise.}
  \end{array}
  \right.
  \]
\end{itemize}

This completes the definition of $x^{*,i}$.

Intuitively, stages $n_k,\dots,m_{k+1}-1$ are used to let the filtration select the player instead of the correlation device,
and stages $m_{k+1},\dots,n_{k+1}-1$ allow the selected players to stop as indicated by the sunspot approximate equilibrium.

\bigskip
\noindent\textbf{Step~5: The payoff under $x^*$.}

We next compare the continuation payoffs generated by $x^*$ in $\Gamma$ to those of
$\xi^*$ in $\calQ_{R}$.

For every sequence $h_k=(j_0,\dots,j_{k-1}) \in I^k$, let
$\gamma_k^i(h_k,x^*)$ denote player $i$'s payoff in $\Gamma$ under $x^*$ on the set $C(h_{k})$, conditional on no player having stopped before stage $n_k$.
Similarly, let $\gamma_{\calQ,k}^i(h_k,\xi^*)$ denote player $i$'s continuation
payoff in $\calQ_{R}$ under $\xi^*$,
conditional on
$(\widehat i_0,\dots,\widehat i_{k-1})=h_k$.

We argue that
for every player $i\in I$, every $k\geq 0$, and every history
$h_k\in I^k$,
\begin{equation}
\label{equ:payoff:diff:1}
\bigl|\gamma_k^i(h_k,x^*)-\gamma_{\calQ,k}^i(h_k,\xi^*)\bigr|
\leq 6\sum_{\ell\geq k}\eta_\ell .
\end{equation}
Indeed,
fix $k \geq 0$ and $h_k \in I^k$.
We compare block $k$ of the construction of $x^*$ to stage $k$ of $\xi^*$.

Under the sunspot $\ep$-equilibrium $\xi^*$, conditional on $h_k$, the
recommended player at stage $k$ is distributed according to
$(\Leb(\widehat i_k = i\mid h_k))_{i \in I}$.
By Lemma~\ref{lem:simulate-history}(c),
under $x^*$, conditional on the same $h_k$ and on no stopping before stage
$n_k$,
the probability of the event $C(h_k,i)$ is 
$\eta_k$-close to $\Leb(\widehat i_k = i\mid h_k)$.

Conditional on $C(h_k,i)$, the stopping probability during
block $k$ under $x^*$ equals $\widehat\delta_k(h_k,i)$.
Moreover, if stopping occurs in block $k$, then
the realized payoff vector is $Z_i$.

Therefore, the law of the absorption outcome in block $k$ under $x^*$
differs from that under $\xi^*$ by at most $3\eta_k$ in total
variation.
Since payoffs are between $-1$ and $1$,
the one-stage effect on
expected payoff is at most $6\eta_k$.
If termination does not occur in block $k$,
then we proceed to block $k+1$.
Hence
\[
|\gamma_k^i(h_k,x^*)-\gamma_{\calQ,k}^i(h_k,\xi^*)|
\leq 6\eta_k + 
\sum_{j \in I} \nu_j(h_k)
\bigl(1-\widehat\delta(h_k,j)\bigr)
 \bigl|\gamma_{k+1}^i((h_{k},j),x^*)-\gamma_{\calQ,k+1}^i((h_{k},j),\xi^*)\bigr|.
\]
Since under $x^*$ the play terminates with probability 1,
iterating this inequality yields \eqref{equ:payoff:diff:1}.

\bigskip
\noindent\textbf{Step~6: $x^*$ is a $14\ep$-equilibrium.}

Fix $i\in I$ and a strategy $x^i$ of player~$i$.
Define a strategy $\xi^i$ in $\calQ_{\mathrm{sun}}$ by
\[
\xi^i_k(h_k,j):=
\prob_{x^i,x^{*,-i}}
\bigl(n_k\leq t^i<n_{k+1}
\mid C(h_k,j),\theta\geq n_k\bigr),
\qquad \forall j\in I,
\]
with the convention that the conditional probability is 0 when the conditioning event
has probability 0.
This strategy mimics $x^i$ in the quitting game $\calQ_{\mathrm{sun}}$.
By Lemma~\ref{lem:simulate-history}(c),
\begin{align}
\label{equ:967}
\sum_{h_k\in I^k}\sum_{j\in I}
\Bigl|
\prob_{x^i,x^{*,-i}}
\bigl(C(h_k,j)\cap\{\theta\geq n_k\}\bigr)
-\nu_j(h_k)\prob_{x^i,x^{*,-i}}
\bigl(C(h_k)\cap\{\theta\geq n_k\}\bigr)
\Bigr|
\leq\eta_k.
\end{align}
%Indeed, integrate the estimate in part~(c)
%against the $\calF_{n_k}$-measurable conditional
%probability of survival to $n_k$, and sum over
%$h_k$ and $j$.

Conditional on $C(h_k,j)\cap \{\theta\geq n_k\}$,
under $(x^i,x^{*,-i})$,
a player in $I \setminus \{i\}$
stops in block $k$ with probability
at most $\ep/3$, independently of player~$i$'s
stopping time.
Hence, under $(x^i,x^{*,-i})$,
the probability that player~$i$ and some other player
stop in the same block,
with no stopping in an earlier block,
is at most $\ep/3$.

On $C(h_k,j)$, player~$j$ is the only player
who stops with positive probability under $x^*$
during block $k$.
Apart from the events considered in the preceding
paragraph, in which player~$i$ and an opponent
have their stopping stages in the same block,
if the game terminates during block $k$,
it terminates by a single player.
If player~$j\neq i$ stops alone in $\Gamma$, player~$i$ receives
$Z_j^i=R_{ij}$, as in $\calQ_{\mathrm{sun}}$.
If player~$i$ stops alone, her payoff in $\Gamma$
is at most $R_{ii}$, which is her payoff upon
stopping alone in $\calQ_{\mathrm{sun}}$.
%
%Fix $N\geq1$.
%For this comparison, count the payoff in $\Gamma$
%only when $\theta<n_N$, and the payoff in
%$\calQ_{\mathrm{sun}}$ only when $\theta<N$;
%otherwise, count the payoff as zero.
%For each $k<N$, the definition of $\xi^i_k(h_k,j)$
%gives the same probability that player~$i$ stops
%in block $k$ as in stage $k$ of $\calQ_{\mathrm{sun}}$,
%conditional on the corresponding history and
%on no earlier stopping.
%
%Multiply each difference of probabilities in \eqref{equ:967} by player~$i$'s expected payoff
%in $\calQ_{\mathrm{sun}}$ from stage $k$ onward,
%given $(h_k,j)$ and no earlier stopping,
%counting zero if stopping does not occur before $N$.
%Applying this comparison successively for
%$k=0,\ldots,N-1$, and including the error from
%stopping by two players in the same block, gives
It follows that for every $N \in \dN$,
\[
\begin{aligned}
\E_{x^i,x^{*,-i}}
 \left[X^i_{J_\theta,\theta}
       \mathbf1_{\{\theta<n_N\}}\right]
\leq
 \E_{\xi^i,\xi^{*,-i}}
 \left[X^i_{J_\theta,\theta}
       \mathbf1_{\{\theta<N\}}\right]
 +2\sum_{k=0}^{N-1}\eta_k + \frac{2\ep}{3}.
\end{aligned}
\]
Letting $N$ tend to infinity and using bounded
convergence, we obtain
\[
\gamma^i(x^i,x^{*,-i})
\leq
\gamma_{\mathrm{sun}}^i(\xi^i,\xi^{*,-i})
+2\ep
\leq \gamma^i(x^*) + 6\ep\leq \gamma^i(x^*) + 8\ep.
\]

%% =====================================================================
\subsection{Construction when $R^{\mathrm{sp}}$ is not $\ep$-close to a Q-Matrix}
\label{sec:Bcase}

We now assume that $R^{\mathrm{sp}}$ is not $\ep$-close to any
Q-matrix
and that the conditions in Lemmas~\ref{lemma:simple:eq} 
and~\ref{lemma:simple:2} do not hold. 
By (AS3), there is $\kappa > 0$ such that 
all payoffs are an integer multiple of $\kappa$.
Without loss of generality, 
assume that $\ep$ is such that $\kappa/\ep$ is an integer,
so that all payoffs are an integer multiple of $\ep$.

The construction is based on Sections~5.3 and
Sections 7.2.4--7.2.5 of \cite{SS04},
but it is more complicated because of the presence of more than two players.

Let $\delta_0(\ep)\in(0,\ep]$ be a constant
satisfying 
Proposition~\ref{prop:absorbing-eq}.
Fix $\delta \in \bigl(0,\min\{\delta_0(\ep)/2, \ep/4, 1/4\}\bigr)$ such that $1/\delta$ is an integer,
and fix $\eta \in \bigl(0,\min\{\delta,\ep\}/4\bigr)$.
Choose a sequence $(\eta_k)_{k\geq0}$ of positive numbers such that
\begin{equation}
\label{eq:zeta-sum}
   \sum_{k\in \dN}\eta_k<\eta^3.
\end{equation}

\bigskip
\noindent\textbf{Step 1: Defining blocks containing all vectors in $\range(X)$.}

By (AS5), 
$\prob(X_{t}=g\text{ i.o.})=1$
for every $g \in \range(X)$.
Hence there is an increasing sequence of integers $0=n_0<n_1<n_2<\cdots$
such that
\[
 \prob\Bigl(
  \forall g\in \range(X) \hbox{ there exists } n_k \leq t < n_{k+1} - |\range(X)|
        \text{ such that }X_{t}=g
 \Bigr)>1-\eta_k .
\]
As in Step 3 of the proof of Lemma~\ref{lem:simulate-history}, by changing the payoff process on
an event of probability at most $\sum_{k\in \dN}\eta_k < \eta^3$ without harming the fact that the payoff process is adapted, we may assume that 
\begin{equation}
\label{equ:875}
\prob\Bigl(
\forall g\in \range(X)\  \hbox{ there exists } n_k \leq t < n_{k+1}
        \text{ such that }X_{t}=g
 \Bigr)=1.
\end{equation}
Indeed, we order the elements in $\range(X)$ according to some fixed order,
and for every $k \in \dN$,
on the complement of the event
\[ \left\{\forall g\in \range(X) \hbox{ there exists } n_k \leq t < n_{k+1} - |\range(X)|
        \text{ such that }X_{t}=g \right\}, \]
we set $X_{n_{k+1}-|\range(X)| + \ell}$ to be equal to the $\ell$'th element in $\range(X)$,
for every $0 \leq \ell < |\range(X)|$. 

\bigskip
\noindent\textbf{Step 2: Defining approximating trees.}

By Theorem~\ref{cor:approx} applied with $\delta^2$ (instead of $\delta$) and $54\sqrt\ep$ (instead of $\ep$), there is a consistent
(in the sense
of Definition~\ref{def:NT})
$\tau$-function $(\tau,\tau') \longmapsto \calT_{\tau,\tau'}$ for $(n_k)_{k \in \dN}$
that assigns to every two $(n_k)_{k\in\dN}$-consistent bounded stopping
times $\tau<\tau'$ 
and every $\omega \in \Omega$ a game on a tree
$\calT_{\tau,\tau'}(\omega)$ satisfying the conclusions of
Theorem~\ref{cor:approx} with $\delta^2$ instead of $\delta$ and
$54\sqrt\ep$ instead of $\ep$.
In particular, the sets $\bigl(\mathcal{G}_k, (q_{G,k})_{G \in \calG_k}\bigr)_{k \in \dN}$ that define $\calT_{\tau,\tau'}$ are a $\delta^2$-approximation.
%Apply Theorem~\ref{cor:approx} with 
%$\delta_{\rm app} \in (0,\min\{\eta^4,\ep/(4n)\})$
%and obtain a consistent $\tau$-function
%$(\tau,\tau')\longmapsto\calT_{\tau,\tau'}$.

By \eqref{equ:875},
$Z_{\calT_{n_k,n_{k+1}}(\omega),i} = Z_i$ a.s.,
for every $i \in I$ and every $k \in \dN$.
Therefore, in particular,
$R^{\mathrm{sp}}_{\calT_{n_k,n_{k+1}}(\omega)} = R^{\mathrm{sp}}$ a.s., for every $k \in \dN$.

We argue that 
for every two bounded and $(n_k)_{k \in \dN}$-consistent stopping times $\tau < \tau'$ and almost every $\omega \in \Omega$, 
after deleting
nodes that are reached with probability zero, 
the tree
$\calT_{\tau,\tau'}(\omega)$ belongs to $\mathbb T_\ep$.
Indeed, (T1) follows from (AS2),
(T2) follows from (AS3),
(T3) follows from (AS4),
(T4) follows from Step 1 and since in this section we assumed that $R^{\mathrm{sp}}$ is not $\ep$-close to any Q-matrix,
(T5) holds after
deleting unreachable nodes,
and (T6) follows from Step 1 and from property (6) in Definition~\ref{def:delta-approx}.

%To verify (T6), fix a tree
%$\calT_{\tau,\tau'}(\omega)$ and a path
%from its root to a leaf that has positive
%probability under the approximating
%transition probabilities.
%By property~(6) in
%Definition~\ref{def:delta-approx},
%the atom $G$ corresponding to the leaf
%has positive probability under $\prob$.
%Since $\tau$ and $\tau'$ are
%$(n_k)_{k\in\dN}$-consistent and
%$\tau<\tau'$, the set of stages between $\tau$ and $\tau'$ (excluded)
%contains at least one
%complete block $\{n_k,\dots,n_{k+1}-1\}$.
%By \eqref{equ:875}, every payoff vector
%$g\in\range(X)$ occurs in this interval
%almost surely.
%Since $X_t$, $\tau$, and $\tau'$
%are constant on the appropriate atoms
%of the approximation, every such
%payoff vector occurs along the
%chosen path.
%Thus every path having positive
%probability contains every
%$g\in\range(X)$, which proves (T6).

%For every such tree apply the canonical construction of
%Proposition~\ref{prop:absorbing-eq}.  We denote its color by
%$c(\calT)$, its trimming numbers by $\lambda_j(\calT)$, and the
%corresponding high-box profiles by $\sigma^j(\calT)$.

\bigskip
\noindent\textbf{Step 3: Defining a consistent $\tau$-function.}

We next define a $\tau$-function (for $(n_k)_{k \in \dN}$)
with values in a finite product set $\calA \times (I \cup \{\varnothing\})$.
%Moreover,
%we associate to each pair $(\tau,\tau')$ of $(n_k)_{k \in \dN}$-consistent bounded stopping times a strategy profile $\alpha^*_{\tau,\tau'}$.
To do so, we first discretize the set of possible continuation payoffs.
Let $\calA \subset [-1,1]^n$ be a finite set such that 
(i) all the coordinates of all points in $\calA$ are multiple of $\delta$
(which, by definition, is the inverse of an integer), and
(ii)
$\bigcup_{a \in \calA}\Pi(a) = [-1,1)^n$,
where
\[
\Pi(a) := \prod_{i\in I}[a^i,a^i+\delta), 
\qquad \forall a \in \calA.
\]

Fix two $(n_k)_{k \in \dN}$-consistent bounded stopping times $\tau<\tau'$.
Recall that the function $\omega \mapsto \calT_{\tau,\tau'}(\omega)$ is constant on each atom of $\calG_\tau$.
Fix such an atom $G \in \calG_\tau$, so that $\tau$ is constant on $G$.
In particular, $\calT_{\tau,\tau'}(\omega)$ is a.s.~equal on $G$ to some fixed game on a tree $\calT \in \mathbb T_\ep$.

By Proposition~\ref{prop:absorbing-eq} with $d=\delta_0(\ep)$, 
the game $\calT_{\tau,\tau'}$ admits a stationary
$54\sqrt\ep$-equilibrium $\alpha^*$ that satisfies properties (i)--(iii) of that proposition.
The payoff under $\alpha^*$,
$\gamma(\alpha^* ; \calT)$,
lies in some of the sets $(\Pi(a))_{a \in \calA}$.
If there is a player $i\in I$ such that $a^i < -8\sqrt\ep$ and $\pi(0^i,\alpha^{*,-i};\calT) < (\delta_0(\ep))^2$,
then by Proposition~\ref{prop:absorbing-eq}(iii) there is only one such player, denoted $i(\alpha^*)$.
If there is no such player $i$, set $i(\alpha^*) = \varnothing$.

Among all pairs $(a,i) \in \calA \times (I \cup \{\varnothing\})$
such that, for some $54\sqrt\ep$-equilibrium $\alpha^*$ of $\calT_{\tau,\tau'}$ that satisfies properties (i)--(iii) of Proposition~\ref{prop:absorbing-eq} with $d=\delta_0(\ep)$,
the payoff lies in $\Pi(a)$ and $i(\alpha^*) = i$,
choose the minimal one according to a fixed ordering
of $\calA \times (I \cup \{\varnothing\})$.
On $G$, define $(f_{\tau,\tau'},i_{\tau,\tau'})$ to be equal to that pair $(a,i)$.

Thus, we obtain, for every $\tau<\tau'$, an $\calF_\tau$-measurable
random variable
\[
(f_{\tau,\tau'},i_{\tau,\tau'}) : \Omega \to \calA \times (I \cup \{\varnothing\}).
\]
Since 
$(\tau,\tau') \longmapsto \calT_{\tau,\tau'}$ is consistent in the sense
of Definition~\ref{def:NT}, 
and since $(f_{\tau,\tau'}(\omega),i_{\tau,\tau'}(\omega))$ is the minimal pair $(a,i)$ 
that satisfies the properties described above,
the $\tau$-function $(f,i)$ is consistent in the sense
of Definition~\ref{def:NT}.

\bigskip
\noindent\textbf{Step 4: Selecting an increasing sequence of stopping times.}

Applying Theorem~\ref{thm:ramsey} to the $\tau$-function
$(f,i)$, we obtain an increasing sequence of
$(n_k)_{k\in \dN}$-consistent bounded stopping times
$(\phi_k)_{k \in \dN}$
such that
\begin{equation}
\label{eq:Bcase-const-color}
\prob\bigl(
f_{\phi_0,\phi_1}=f_{\phi_1,\phi_2}=f_{\phi_2,\phi_3}=\cdots,
i_{\phi_0,\phi_1}=i_{\phi_1,\phi_2}=i_{\phi_2,\phi_3}=\cdots
\bigr) > 1-\ep .
\end{equation}
Let
\[
F := \bigl\{
f_{\phi_0,\phi_1}=f_{\phi_1,\phi_2}=f_{\phi_2,\phi_3}=\cdots,
i_{\phi_0,\phi_1}=i_{\phi_1,\phi_2}=i_{\phi_2,\phi_3}=\cdots
\bigr\}
\]
and  
\[ a := f_{\phi_0,\phi_1}. \]

Recall that the approximating tree 
$\calT_{\phi_k,\phi_{k+1}}(\omega)$ is constant on each atom of $\calG_{\phi_k}$.
For every $k\in\dN$ and every $\omega \in \Omega$, choose a stationary $54\sqrt\ep$-equilibrium
$\alpha^*_{\phi_k,\phi_{k+1}}(\omega)$
of the approximating tree 
$\calT_{\phi_k,\phi_{k+1}}(\omega)$
that satisfies conclusions
(i)--(iii) of
Proposition~\ref{prop:absorbing-eq} with $d=\delta_0(\ep)$
and realizes 
$(f_{\phi_k,\phi_{k+1}},
i_{\phi_k,\phi_{k+1}})$,
and such that the function $\omega \mapsto \alpha^*_{\phi_k,\phi_{k+1}}(\omega)$ is constant on each atom of $\calG_{\phi_k}$.
The strategy profiles $(\alpha^*_{{\phi_k},{\phi_{k+1}}}(\omega))_{\omega \in \Omega}$ 
define an $(\mathbb F, \phi_k,\phi_{k+1})$ strategy profile 
(see Eq.~\eqref{equ:defin:e}).

\bigskip

It will be convenient to change the payoff process $X$ on $\Omega \setminus F$.
Since this set has probability smaller than $\ep$,
by Lemma~\ref{lemma:simple}(O1) any approximate equilibrium of the game with the modified payoffs is an approximate equilibrium of $\Gamma$, albeit with a larger $\ep$.
We turn to define the modified payoff process.

\bigskip
\noindent\textbf{Step 5: The auxiliary game $\widetilde \Gamma$.}

Let $\widetilde \Gamma  = (\Omega,I,\calF,\prob,\mathbb{F}, \widetilde X)$ be the stopping game that is similar to $\Gamma$,
except for the payoff process $\widetilde X$,
which is equal to $f_{\phi_0,\phi_1}$ on a small set: 
for every $t \in \dN$ and every $\varnothing \neq J \subseteq I$,
\[ \widetilde X_{J,t}(\omega) := 
\left\{
\begin{array}{ll}
f_{\phi_0,\phi_1}(\omega), & \hbox{if } 
(f_{\phi_0,\phi_1}(\omega),i_{\phi_0,\phi_1}(\omega)) \neq (f_{\phi_k,\phi_{k+1}},i_{\phi_k,\phi_{k+1}}) \hbox{ and } \phi_k(\omega) \leq t < \phi_{k+1}(\omega),\\
X_{J,t}(\omega), & \hbox{otherwise}.
\end{array}
\right.
\]
Thus, on $(f_{\phi_0,\phi_1}(\omega),i_{\phi_0,\phi_1}(\omega)) \neq (f_{\phi_k,\phi_{k+1}},i_{\phi_k,\phi_{k+1}})$,
the payoff vector is the same, irrespective of the set of stopping players and of the termimnal stage $t$,
provided $\phi_k(\omega) \leq t < \phi_{k+1}(\omega)$.

Since $f_{\phi_k,\phi_{k+1}}$ is $\calF_{\phi_k}$-measurable for every $k \in \dN$, the process $\widetilde X$ is $\mathbb F$-adapted.
%Denote by $\widetilde \gamma(\cdot)$ the payoff function of $\widetilde \Gamma$.
%The sets of strategies of the players in $\Gamma$ and $\widetilde \Gamma$ coincide.
Since 
\begin{equation}
\label{equ:small:1}
\prob(\widetilde X_t \neq X_t \hbox{ for some } t \in \mathbb N) \leq \prob(\Omega \setminus F) < \ep, 
\end{equation}
by Lemma~\ref{lemma:simple}(O1),
any $(116\sqrt\ep+3\ep)$-equilibrium of $\widetilde\Gamma$ is a $(116\sqrt\ep+11\ep)$-equilibrium of $\Gamma$.

For every $k\in\dN$ and every $\omega\in\Omega$,
let $\widetilde\calT_{\phi_k,\phi_{k+1}}(\omega)$ be the game on a tree which is similar to $\calT_{\phi_k,\phi_{k+1}}(\omega)$,
except that the payoff is derived from $\widetilde X$ rather than from $X$:
If
$(f_{\phi_k,\phi_{k+1}}(\omega),
 i_{\phi_k,\phi_{k+1}}(\omega))
=
((f_{\phi_0,\phi_{1}}(\omega),i_{\phi_0,\phi_1}(\omega))$,
then $\widetilde\calT_{\phi_k,\phi_{k+1}}(\omega)=\calT_{\phi_k,\phi_{k+1}}(\omega)$.
Otherwise, the terminal payoff 
at every non-leaf node of $\widetilde\calT_{\phi_k,\phi_{k+1}}(\omega)$ is $f_{\phi_0,\phi_{1}}(\omega)$.

In the latter case, we redefine 
$\alpha^{*}_{\phi_k,\phi_{k+1}}(\omega)$,
and set it to be $(1^i)_{i \in I}$;
that is, all players stop with probability 1 at all nodes.

Almost surely the strategy profile $\alpha^*_{\phi_k,\phi_{k+1}}(\omega)$ is a $54\sqrt\ep$-equilibrium in $\widetilde \calT_{\phi_k,\phi_{k+1}}(\omega)$.
Indeed, on $\{(f_{\phi_0,\phi_1}(\omega),i_{\phi_0,\phi_1}(\omega)) = (f_{\phi_k,\phi_{k+1}}, i_{\phi_k,\phi_{k+1}})\}$ the games $\widetilde \calT_{\phi_k,\phi_{k+1}}(\omega)$ and $\calT_{\phi_k,\phi_{k+1}}(\omega)$ coincide,
and on its complement the terminal payoff is constant,
and since there are at least two players, the strategy profile in which all players always stop with probability 1 is a 0-equilibrium.

\bigskip
\noindent\textbf{Step 6: Definition of a strategy profile $x^*$.}

Define a strategy profile $x^*$ in $\Gamma$ from stage ${\phi_0}$ onward by gluing together the strategies $(\alpha^*_{{\phi_k},{\phi_{k+1}}}(\omega))_{\omega \in \Omega, k \in \dN}$:
\[ x^{*,i}_t(\omega) := \alpha^{*,i}_{{\phi_k(\omega)},{\phi_{k+1}(\omega)}}(\omega)
\hbox{ on the set } {\phi_k(\omega)} \leq t < {\phi_{k+1}(\omega)}.
\]
We note that the redefined strategy profile $x^*$ is $\mathbb F$-adapted.

We will prove that (a proper modification of) $x^*$ is a $(116\sqrt\ep+11\ep)$-equilibrium of $\Gamma$
when the game starts in stage $\phi_0$.
By Lemma~\ref{lemma:simple}(O3),
this will imply that $\Gamma$ admits a $(116\sqrt\ep+11\ep)$-equilibrium.

\bigskip
\noindent\textbf{Step 7: The payoff under $x^*$ from stage $\phi_0$ in $\widetilde \Gamma$ is close to $a$.}

Consider the stopping game $\widetilde\Gamma$ that starts at stage $\phi_0$.
Denote by $\widetilde\rho^i$ 
the function $\rho^i$ for $\widetilde\Gamma$.%
\footnote{The function $\pi$ is the same for $\Gamma$ and $\widetilde\Gamma$.}
For every strategy profile $x$, denote
\[
\widetilde\gamma^i(x\mid\calF_{\phi_0},
\theta\geq\phi_0)
:=
\E_x\left[
\widetilde X^i_{J_\theta,\theta}
\mathbf{1}_{\{\theta<\infty\}}
\mid\calF_{\phi_0},\theta\geq\phi_0
\right].
\]

Since $\pi(\alpha^*_{{\phi_k},{\phi_{k+1}}}(\omega);\calT_{\phi_k,\phi_{k+1}})\ge\delta$ for almost every $\omega \in \Omega$ and every $k \in \dN$,
and since the game on a tree $\calT_{\phi_k,\phi_{k+1}}$ is derived from a $\delta^2$-approximation,
$\pi(x^*;\mathbb F, \phi_k,\phi_{k+1}) \geq \delta(1-\delta)$.
Hence,
the play under $x^*$ terminates with probability $1$.

For every $k \in \dN$ and almost every $\omega \in \Omega$,
$\gamma(\alpha^*_{{\phi_k},{\phi_{k+1}}}(\omega);\widetilde\calT_{\phi_k,\phi_{k+1}}) \in \Pi(a(\omega))$,
and hence 
\[ \rho^i(\alpha^*_{{\phi_k},{\phi_{k+1}}}(\omega);\widetilde\calT_{\phi_k,\phi_{k+1}}) \geq a^i\pi(\alpha^*_{{\phi_k},{\phi_{k+1}}}(\omega);\widetilde\calT_{\phi_k,\phi_{k+1}}), \qquad a.s. \]

On $\{(f_{\phi_0,\phi_1}(\omega),i_{\phi_0,\phi_1}(\omega)) = (f_{\phi_k,\phi_{k+1}}, i_{\phi_k,\phi_{k+1}})\}$, 
the game on a tree
$\widetilde\calT_{\phi_k,\phi_{k+1}} = \calT_{\phi_k,\phi_{k+1}}$ is derived from a $\delta^2$-approximation.
On the complementary event, all players stop immediately and
$\widetilde\rho^i \equiv a^i$.
In both cases, for every $k \in \dN$,
\begin{equation}
\label{equ:estimate:1}
\widetilde\rho^i(x^*;\mathbb F, \phi_k,\phi_{k+1}) \geq a^i\pi(x^*;\mathbb F, \phi_k,\phi_{k+1}) - {2\delta^2}/{2^k}, \qquad a.s. 
\end{equation}
Multiplying both sides of \eqref{equ:estimate:1}
by $\mathbf{1}_{\{\theta\geq\phi_k\}}$,
summing over all $k\in\dN$,
and taking conditional expectations given
$\calF_{\phi_0}$ and $\{\theta\geq\phi_0\}$,
we obtain
\begin{align*}
\widetilde\gamma^i
(x^*\mid&\calF_{\phi_0},\theta\geq\phi_0)
=
\E_{x^*}\left[
\sum_{k\in\dN}
\mathbf{1}_{\{\theta\geq\phi_k\}}
\widetilde\rho^i
(x^*;\mathbb F,\phi_k,\phi_{k+1})
\;\middle|\;
\calF_{\phi_0},\theta\geq\phi_0
\right]
\\
&\geq
a^i
\E_{x^*}\left[
\sum_{k\in\dN}
\mathbf{1}_{\{\theta\geq\phi_k\}}
\pi(x^*;\mathbb F,\phi_k,\phi_{k+1})
\;\middle|\;
\calF_{\phi_0},\theta\geq\phi_0
\right]
-4\delta^2
\\
&=
a^i
\prob_{x^*}
(\theta<\infty\mid\calF_{\phi_0},
\theta\geq\phi_0)
-4\delta^2
\\
&\geq a^i-\ep,
\end{align*}
where the two equalities follow from
the tower property of conditional
expectations, 
and the last inequality
follows since under $x^*$ the game terminates with probability 1
and by the choice of $\delta$.

\bigskip
\noindent\textbf{Step 8: Bounding the payoff under deviations from $x^*$.}

Consider any $i \in I$ and any deviation $x^i$ of player~$i$ in $\widetilde \Gamma$.
Then for every $k \in \dN$,
\begin{equation}
\label{eq:Bcase-approx}
\widetilde\rho^i(x^i,x^{*,-i}_{\phi_k,\phi_{k+1}};\mathbb F, \phi_k,\phi_{k+1})
\le
(a^i+108\sqrt\ep)
\pi(x^i,x^{*,-i}_{\phi_k,\phi_{k+1}};\mathbb F, \phi_k,\phi_{k+1}) + \eta_{\phi_k}.
\end{equation}
Indeed, if $\widetilde\calT_{\phi_k,\phi_{k+1}}=\calT_{\phi_k,\phi_{k+1}}$,
\eqref{eq:Bcase-approx} follows from Theorem~\ref{cor:approx},
while if $\widetilde\calT_{\phi_k,\phi_{k+1}}\neq\calT_{\phi_k,\phi_{k+1}}$,
then in the game $\widetilde\calT_{\phi_k,\phi_{k+1}}$ the payoff under any strategy profile that terminates with positive probability is $f_{\phi_0,\phi_1}$.

Multiplying both sides of
\eqref{eq:Bcase-approx} by
$\mathbf{1}_{\{\theta\geq\phi_k\}}$,
summing over all $k\in\dN$, and taking
conditional expectations given
$\calF_{\phi_0}$ and $\{\theta\geq\phi_0\}$,
we obtain
\begin{align}
\label{equ:nonq:eq2}
&\widetilde\gamma^i
(x^i,x^{*,-i}\mid
\calF_{\phi_0},\theta\geq\phi_0)
\nonumber\\
&=
\E_{x^i,x^{*,-i}}\left[
\sum_{k\in\dN}
\mathbf{1}_{\{\theta\geq\phi_k\}}
\widetilde\rho^i
(x^i,x^{*,-i};
\mathbb F,\phi_k,\phi_{k+1})
\;\middle|\;
\calF_{\phi_0},\theta\geq\phi_0
\right]
\nonumber\\
&\leq
(a^i+108\sqrt\ep)
\E_{x^i,x^{*,-i}}\left[
\sum_{k\in\dN}
\mathbf{1}_{\{\theta\geq\phi_k\}}
\pi(x^i,x^{*,-i};
\mathbb F,\phi_k,\phi_{k+1})
\;\middle|\;
\calF_{\phi_0},\theta\geq\phi_0
\right]
+\sum_{t\in\dN}\eta_t
\nonumber\\
%&=
%(a^i+108\sqrt\ep)
%\prob_{x^i,x^{*,-i}}
%(\theta<\infty\mid
%\calF_{\phi_0},\theta\geq\phi_0)
%+\sum_{t\in\dN}\eta_t
%\nonumber\\
&\leq
(a^i+108\sqrt\ep)
\prob_{x^i,x^{*,-i}}
(\theta<\infty\mid
\calF_{\phi_0},\theta\geq\phi_0)
+\ep,
\end{align}
where the equality and the second inequality follow
from the tower property of conditional
expectations,
and the first inequality follows from
\eqref{eq:Bcase-approx} and the fact that
$(\phi_k)_{k\in\dN}$ is strictly increasing,
so that
$\sum_{k\in\dN}
\mathbf{1}_{\{\theta\geq\phi_k\}}
\eta_{\phi_k}
\leq\sum_{t\in\dN}\eta_t$.
%The last inequality follows from
%\eqref{eq:zeta-sum}.

%Multiplying both sides of \eqref{eq:Bcase-approx} by $\mathbf{1}_{\theta \geq \phi_k}$,
%summing over all $k \in \dN$,
%and taking expectation,
%\begin{align} 
%\label{equ:nonq:eq2}
%\widetilde\gamma^i(x^i,x^{*,-i} \mid \calF_{\phi_0}, \theta \geq \phi_0) &= 
%\E\left[\sum_{k \in \dN}
%\mathbf{1}_{\theta \geq \phi_k}
%\widetilde\rho^i(x^i,x^{*,-i}_{\phi_k,\phi_{k+1}};\mathbb F, \phi_k,\phi_{k+1})
%\mid \calF_{\phi_0}, \theta \geq \phi_0\right]\\
%&\le
%(a^i+108\sqrt\ep)\prob_{x^i,x^{*,-i}}(\theta < \infty \mid \calF_{\phi_0}, \theta 
%\geq \phi_0)+\ep.
%\nonumber
%\end{align}
It follows that if $a^i \geq -8\sqrt\ep$ or if $\prob_{x^i,x^{*,-i}}(\theta < \infty \mid \calF_{\phi_0}, \theta \geq \phi_0)  =1$,
then the deviation does not yield a profit higher than $116\sqrt\ep+\ep$.
%
%&\gamma^i(x^i,x^{*,-i} \mid \calF_{\phi_0}, \theta \geq \phi_0) \prob_{x^i,x^{*,-i}}(\theta < \infty \mid \calF_{\phi_0}, \theta 
%\nonumber
%\geq \phi_0) \\
%&\leq a^i+2\ep + \sum_{k \in \dN}\eta_k\\
%\nonumber
%&\leq a^i \prob_{x^i,x^{*,-i}}(\theta < \infty \mid \calF_{\phi_0}, \theta \geq \phi_0) + 3\ep.
%\end{align}
%

\bigskip
\noindent\textbf{Step 9: On
$\{i_{\phi_0,\phi_1}=\varnothing\}$,
$x^*$ is a $(116\sqrt\ep+\ep)$-equilibrium
of $\widetilde\Gamma$ from stage $\phi_0$.}

Fix a player $i\in I$.
If $a^i\geq-8\sqrt\ep$, the conclusion
follows from Step~8.

Suppose that $a^i<-8\sqrt\ep$.
Since $i_{\phi_0,\phi_1}=\varnothing$,
for every $k \in \dN$ such that $\widetilde\calT_{\phi_k,\phi_{k+1}}=\calT_{\phi_k,\phi_{k+1}}$
we have
$\pi(0^i,\alpha_{\phi_k,\phi_{k+1}}^{*,-i};
\calT_{\phi_k,\phi_{k+1}})
\geq\delta_0^2$.
And if $\widetilde\calT_{\phi_k,\phi_{k+1}}\neq\calT_{\phi_k,\phi_{k+1}}$, every player other
than $i$ stops with probability 1 throughout the block.
Since the approximating trees are obtained
from a $\delta^2$-approximation,
Eq.~\eqref{equ:app:2} implies that
\[
\pi(0^i,x^{*,-i};
\mathbb F,\phi_k,\phi_{k+1})
\geq\delta_0^2-\delta^2>0,
\]
almost surely, for every $k\in\dN$.
Consequently, under every deviation $x^i$
of player~$i$, 
\begin{equation}
\label{equ:931}
\prob_{x^i,x^{*,-i}}
\bigl(\theta\geq\phi_k
\mid\calF_{\phi_0},\theta\geq\phi_0\bigr)
\leq(1-\delta_0^2 + \delta^2)^k,
\qquad k\in\dN.
\end{equation}
Since $\delta_0^2 - \delta^2>0$,
taking the limit in \eqref{equ:931} as $k$ goes to $\infty$ we obtain
$\prob_{x^i,x^{*,-i}}
\bigl(\theta<\infty
\mid\calF_{\phi_0},\theta\geq\phi_0\bigr)
=1$.
By Step~8, player~$i$ cannot gain more
than $116\sqrt\ep+\ep$ by deviating.
Since $i$ was arbitrary, the claim
follows.

% \bigskip
% \noindent\textbf{Step 8: On $\{i_{\phi_0,\phi_1} = \varnothing\}$, 
% $x^*$ is a $(116\sqrt\ep+\ep)$-equilibrium of $\widetilde\Gamma$ from stage $\phi_0$.}

% On $\{i_{\phi_0,\phi_1} = \varnothing\}$, 
% for every player $i \in I$ we have
% $a^i \geq -8\sqrt\ep$ or 
% $\pi(0^i,\alpha^{*,-i}_{\phi_k,\phi_{k+1}};\calT_{\phi_k,\phi_{k+1}})\geq \delta_0^2$
% for every $k \in \dN$.

% By Step~7,
% $x^*$ is a $(116\sqrt\ep+\ep)$-equilibrium of $\widetilde\Gamma$ from stage $\phi_0$.

\bigskip
\noindent\textbf{Step 10: Approximate equilibrium in $\widetilde\Gamma$ on $\{i_{\phi_0,\phi_1} \neq \varnothing\}$.}

Fix $i_0 \in I$ and consider the set $\{i_{\phi_0,\phi_1} = i_0\}$.
Let $K \in \dN$ be sufficiently large such that 
\begin{align}
\label{equ:K}
\prob_{x^*}(\theta < \phi_K \mid \{i_{\phi_0,\phi_1} = i_0\}, \theta \geq \phi_0) &\geq 1-\ep,\\
(1-\delta_0(1-\delta_0)+\delta^2)^K &\leq \ep.
\label{equ:K2}
\end{align}

Let $\widehat x^*$ be the following strategy profile, 
which is defined only at and after stage $\phi_0$:
\begin{itemize}
\item Until stage $\phi_K$ (excluding), $\widehat x^*$ coincides with $x^*$.
\item 
By Proposition~\ref{prop:absorbing-eq}(iii)
applied to the tree $\calT_{\phi_0,\phi_1}$,
there are $g \in \range(X)$ and $(\beta^j)_{j \neq i_0}$,
not all of them 0,
such that 
$\max\{A,B\}\le a^{i_0}+\delta + \ep^2$, where
$A$ and $B$ are given as in Proposition~\ref{prop:absorbing-eq}(iii).
At and after stage $\phi_K$,
%the strategy profile $\widehat x^*$ is defined as follows:
%\begin{itemize}
%\item During blocks $k \geq K$ such that $(f_{\phi_0,\phi_1}(\omega),i_{\phi_0,\phi_1}(\omega)) \neq (f_{\phi_k,\phi_{k+1}},i_{\phi_k,\phi_{k+1}})$,
%all players stop with probability 1 at every stage.
%\item 
%During blocks $k \geq K$ such that $(f_{\phi_0,\phi_1}(\omega),i_{\phi_0,\phi_1}(\omega)) = (f_{\phi_k,\phi_{k+1}},i_{\phi_k,\phi_{k+1}})$,
no player stops when $\widetilde X_t \neq g$;
Whenever $\widetilde X_t=g$,
each player $j \neq i_0$ stops with probability $\beta^j$ and player $i_0$ stops with probability 0.
%\end{itemize}
\end{itemize}

By \eqref{equ:K},
$\|\widetilde\gamma(\widehat x^* \mid \{i_{\phi_0,\phi_1} = i_0\},\theta \geq \phi_0) - \widetilde\gamma(x^* \mid \{i_{\phi_0,\phi_1} = i_0\},\theta \geq \phi_0)\|_\infty \leq 2\ep$,
and hence in particular 
$\widetilde\gamma^j(\widehat x^* \mid \{i_{\phi_0,\phi_1} = i_0\},\theta \geq \phi_0) \geq a^j - 5\ep$ for every $j \in I$.

We next estimate the gains from deviations.
We start with player~$i_0$.
The conditional probability that $\widetilde\calT_{\phi_k,\phi_{k+1}}\neq\calT_{\phi_k,\phi_{k+1}}$ for some $k \in \dN$ is at most
$\prob(F^c\mid\calF_{\phi_0})$.

On $F$ we have
$\widetilde\calT_{\phi_k,\phi_{k+1}}=\calT_{\phi_k,\phi_{k+1}}$
for every $k \in \dN$.
We argue that $a^{i_0} \geq y^{i_0} - 57\sqrt\ep$.
Indeed, let $k \in \dN$.
Since $\pi(0^i,\alpha^{*,-i}_{\phi_k,\phi_{k+1}}) < \delta_0^2$ and by \eqref{equ:875},
by deviating from $\alpha^{*}_{\phi_k,\phi_{k+1}}$ in $\calT_{\phi_k,\phi_{k+1}}$ to the strategy that stops with probability 1 whenever $X_{\{i_0\},s} = Z_{i_0}$,
player~$i_0$ obtains at least $y^{i_0} - 2\delta_0^2$.
Since $\alpha^*_{\phi_k,\phi_{k+1}}$ is a $54\sqrt\ep$-equilibrium, 
\begin{equation}
\label{equ:877}
a^{i_0} \geq \gamma^{i_0}(\alpha^*_{\phi_k,\phi_{k+1}};\calT_{\phi_k,\phi_{k+1}}) - \delta \geq y^{i_0} - \delta - 2\delta_0^2 - 54\sqrt\ep \geq y^{i_0} - 57\sqrt\ep.
\end{equation}

Since $\widetilde\calT_{\phi_k,\phi_{k+1}}=\calT_{\phi_k,\phi_{k+1}}$
for every $k \in \dN$, we have
$X_t = g$ i.o., and hence 
for every strategy $x^{i_0}$ of player $i_0$ in $\widetilde\Gamma$,
under 
$(x^{i_0},\widehat x^{*,-i_0})$
the play terminates with probability 1.

If player~$i_0$ stops at a stage $t\geq \phi_K$ where $X_t \neq g$,
by Proposition~\ref{prop:absorbing-eq}(iii) 
and \eqref{equ:877} her payoff is bounded by $y^{i_0} + 2\delta_0 \leq a^{i_0} + 4\delta + 54\sqrt\ep$.
Hence, by the choice of $g$ and $(\beta^j)_{j \neq i}$, we have
$\widetilde\gamma(x^{i_0},\widehat x^{*,-i_0} \mid F \cap \{i_{\phi_0,\phi_1} = i_0\},\theta \geq \phi_0) \leq a^{i_0}+57\sqrt\ep$.

We turn to player $i \neq i_0$.
Since $i_0$ satisfies Proposition~\ref{prop:absorbing-eq}(iii)
and by the definition of $\delta$-approximation,
$\pi(0^i,\alpha^{*,-i}_{\phi_k,\phi_{k+1}};\widetilde\calT_{\phi_k,\phi_{k+1}}) \geq \delta_0 - (\delta_0)^2 - \delta$ for every $k \in \dN$.
By~\eqref{equ:K2},
for every deviation of player~$i$,
the play terminates before stage $\phi_K$ with probability at least $1-\ep$,
and hence, by Steps~7 and~8, player~$i$ cannot profit more than $116\sqrt\ep + 3\ep$ by deviating.

\bigskip
\noindent\textbf{Step 11: Conclusion of the proof.}

Define $\widetilde x^*$ to coincide with $x^*$
on $\{i_{\phi_0,\phi_1}=\varnothing\}$ and with
the profile $\widehat x^*$ defined in Step~10
on $\{i_{\phi_0,\phi_1} \neq \varnothing\}$.
These events are $\calF_{\phi_0}$-measurable,
and hence the resulting strategy profile is $\mathbb F$-adapted.

By Steps~7--10, every player's conditional gain
from deviating from stage $\phi_0$ onward
is at most
$116\sqrt\ep+\ep
+2\prob(F^c\mid\calF_{\phi_0})$.
Using backward induction, as in
Lemma~\ref{lemma:simple}(O2), 
one can extend
$\widetilde x^*$ to an approximate equilibrium in $\widetilde\Gamma$
with error at most
$116\sqrt\ep+\ep+2\prob(F^c)
\leq116\sqrt\ep+3\ep$.

%%%%%%%%%%%%%%%%%%%
\section{Proof of Proposition~\ref{prop:absorbing-eq}}
\label{sec:proof:absorbing-eq}

In this section, we prove Proposition~\ref{prop:absorbing-eq}.
We start by providing some tools that are adapted from \cite{SS04} to the multiplayer setup.

Throughout, $\calT \in \mathbb{T}_\ep$ is a fixed stopping game on a finite tree with player set
$I=\{1,\dots,n\}$ and non-leaf nodes $S_0$.

%% ---------------------------------------------------------------------
\subsection{The Sum of Strategy Profiles}
\label{subsec:sum}
%% ---------------------------------------------------------------------

\begin{definition}[Sum of strategies and of strategy profiles]
\label{def:sum}
Let $i \in I$.
Given finitely many stationary strategies $\alpha^i_1,\dots,\alpha^i_m$ of 
player~$i$, 
their \emph{sum} is the stationary strategy $\alpha^i := \alpha^i_1 \dot+ \cdots \dot+ \alpha^i_m$
defined by
\[
   \alpha^i(s)
   := 1-\prod_{\ell=1}^m\bigl(1-\alpha^i_\ell(s)\bigr),
   \qquad s\in S_0 .
\]
Given finitely many stationary strategy profiles $\alpha_1,\dots,\alpha_m$ with
$\alpha_\ell=(\alpha_\ell^i)_{i\in I}$ for each $\ell \in \{1,\dots,m\}$, their \emph{sum}
$\alpha=\alpha_1\dot+\cdots\dot+\alpha_m$ is defined
coordinate-wise by $\alpha^i:=\alpha_1^i \dot+ \dots \dot+ \alpha_m^i$, for every
$i\in I$.
\end{definition}

To interpret $\alpha^i:=\alpha_1^i \dot+ \dots \dot+ \alpha_m^i$,
suppose that 
at every node $s \in S_0$,
independently of the past,
for each $\ell=1,\dots,m$ player~$i$ flips an unfair coin with parameter $\alpha_\ell^i(s)$,
and stops if at least one of the $m$ outcomes is 1.
The interpretation of $\alpha=\alpha_1\dot+\cdots\dot+\alpha_m$ is that at each node, 
these lotteries are made independently by the $n$ players.

%For $1\le \ell\le m$,
%let $t_\ell$ be the stage in the first round in which the outcome of the coin that corresponds to $\alpha^i_\ell$ is 1;
%$t_\ell=\infty$ if termination does not occur in the first round.
%
%The first stage in which the game terminates under $\alpha$ is $\theta=\min_{1 \leq \ell \leq m} t_\ell$, with
%stopping set $J=\bigcup_{\{\ell \colon t_\ell=\theta\}}J_{t_\ell}$.
For $1 \leq \ell \leq m$ and $i \in I$,
let $t^i_\ell$ be the first stage in which player~$i$ stops under the strategy $\alpha^i_\ell$
(or $\infty$, if player $i$ never stops under $\alpha^i_\ell$).
The first stage in which some player stops under $\alpha_\ell$ is $t_\ell := \min_{i \in I} t^i_\ell$,
and the first stage in which some player stops under $\alpha_1\dot+\cdots\dot+\alpha_m$ is $\theta := \min_{1 \leq \ell \leq m} t_\ell$.
Let $J_\theta := \{ i \in I \colon \min_{1 \leq \ell \leq m} t^i_\ell = \theta\}$
be the set of players who stop at the termination stage under $\alpha_1\dot+\cdots\dot+\alpha_m$.

Recall that $\theta_{\mathrm{fr}}$ is the last stage of the first round.
Let $N:=\sum_{\ell=1}^m \bfone_{\{t_\ell<\theta_{\mathrm{fr}}\}}$ be the number of strategy profiles among $\alpha_1,\dots,\alpha_m$
that stop in the first round.

%\begin{lemma}[Independence given the branch; {\cite[Lemma~8.1]{SS04}}]
%\label{lem:indep-branch}
%Let $s\in S_0$ have depth $j$. For every $k\neq l$, the event
%$\{t_k\le j\}$ and the random variable %$t_k\bfone_{\{t_k\le j\}}$ are
%independent of $t_\ell$ given $F_s$.
%\end{lemma}

The next result relates $\pi(\alpha)$ to $\sum_{\ell=1}^m \pi(\alpha_\ell)$,
and $\rho^i(\alpha)$ to $\sum_{\ell=1}^m \rho^i(\alpha_\ell)$. 

\begin{lemma}[\cite{SS04}, Lemma 8.2]
\label{lem:sum-estimates}
With the notation above,
\[
   \sum_{\ell=1}^m \pi(\alpha_\ell)-\E_{r,\alpha}\bigl[N\bfone_{\{N\ge2\}}\bigr]
   \le\pi(\alpha)\le\sum_{\ell=1}^m \pi(\alpha_\ell),
\]
and, for every player $i\in I$,
\[
   \Bigl|\rho^i(\alpha)-\sum_{\ell=1}^m \rho^i(\alpha_\ell)\Bigr|
   \le\E_{r,\alpha}\bigl[(N+1)\bfone_{\{N\ge2\}}\bigr].
\]
\end{lemma}

\cite{SS04} proved Lemma~\ref{lem:sum-estimates} in the presence of two players.
Their proof does not involve the number of players,
and hence applies to the multiplayer setup.

%% ---------------------------------------------------------------------
\subsection{Heavy and Light Nodes}
\label{subsec:heavy-nodes}
%% ---------------------------------------------------------------------

A node $s$ is \emph{$\eta$-heavy} with respect to a stationary strategy profile $\alpha$ if 
under $\alpha$, 
conditional on visiting node $s$ in the first round,
the game terminates in the first round with probability at least $\eta$.

\begin{definition}[Heavy and light nodes]
\label{def:heavy}
Let $\alpha$ be a stationary strategy profile in $\calT$ and let $\eta>0$.
A node $s\in S_0$ is \emph{$\eta$-heavy} with respect to $\alpha$ if
$\prob_{s,\alpha}(\theta<\theta_{\mathrm{fr}})\ge\eta$, and \emph{$\eta$-light}
otherwise.
We write $H_\eta(\alpha)$ for the set of $\eta$-heavy nodes of
$\alpha$ (or $H_\eta(\alpha;\calT)$,
when we want to emphasize the tree $\calT$).
\end{definition}

Note that for all strategy profiles $\alpha,\alpha'$ and all 
$\eta\ge\eta'>0$,
\begin{equation}
\label{fact:heavy}
   H_{\eta}(\alpha)\subseteq H_{\eta}(\alpha\dot+\alpha')
   \qquad\text{and}\qquad
   H_{\eta}(\alpha)\subseteq H_{\eta'}(\alpha).
\end{equation}

%%%%%%%%%%%%
\subsection{Absorbing 0-Equilibria in $\calT$}

In this section we prove the following result,
which will be used in the proof of Proposition~\ref{prop:absorbing-eq}.

\begin{lemma}
\label{lemma:heavy}
Let $\calT$ be a game on a tree that satisfies properties (T1)--(T5) in the definition of $\mathbb T_\ep$,
and,
in addition, satisfies 
\begin{equation}
\label{equ:871}
 \prob_{r,\mathbf0}
 \bigl(X_{s_t}=g
       \text{ for some }t < \theta_{\mathrm{fr}}\bigr)
 \geq \frac{1}{2}, \qquad \forall g \in \range(X).
\end{equation}
Let $\delta \in \bigl(0,\min\{\ep/16,\delta_{\LCP}(\ep)/16,1/2^{n+4}\}\bigr)$.
At least one of the following two assertions holds:
\begin{itemize}
\item[A1)] 
There is a player $i_0 \in I$ and 
$\beta \in (0,1]$ such that the following properties hold:
\begin{itemize}
\item[i)] Player~$i_0$ is a supporter, or there is $s \in S_0$ and $i \neq i_0$ such that $X_{\{i\},s}^{i_0} \leq y^{i_0}_\calT$.
\item[ii)]
For the unique $g \in \range(X)$ that satisfies $g_{\{i_0\}} = Z_{\calT,i_0}$ we have
$Z^i_{\calT,i_0} \geq (1-\beta) g^i_{\{i\}} + \beta g^i_{\{i,i_0\}}$,
for every $i \neq i_0$.
\item[iii)]
For every $i \neq i_0$
we have
$Z^i_{\calT,i_0} \geq X_{\{i\},s}^i$ 
whenever $X_s \neq g$.
\end{itemize}
%such that $R_{ii_0} \geq R_{ii}$ for every $i \neq i_0$.
\item[A2)] 
There is a stationary $0$-equilibrium $\alpha^*$ in $\calT$ such that the following properties hold:
\begin{itemize}
\item[i)] 
$H_{\delta^2}(\alpha^*)\neq\varnothing$.
\item[ii)]
If $\gamma^i(\alpha^*) \leq y^i_\calT - 5\sqrt\ep$,
then $\pi(0^i,\alpha^{*,-i}) \geq \ep$.
\item[iii)]
If  
$\gamma^i(\alpha^*) < 0$,
then there is a node $s \in S_0$
and real numbers $(\beta^j)_{j \neq i} \in [0,1]^{n-1}$, not all of them 0,
such that 
$\max\{A, B\} \leq \gamma^i(\alpha^*)$,
where $A$ and $B$ are defined as in Proposition~\ref{prop:absorbing-eq}(iii).
%\item[iii)]
%If $\pi(0^i,\alpha^{*,-i}) = 0$,
%then $i$ is a supporter.
\end{itemize}
\end{itemize}
\end{lemma}

We note that Condition~A2(iii) in Lemma~\ref{lemma:heavy} might occur when $y^i_\calT < 5\sqrt\ep$ for some $i \in I$.

Lemma~\ref{lemma:heavy} is close in spirit to Proposition~\ref{prop:absorbing-eq},
with one main difference:
while in Proposition~\ref{prop:absorbing-eq} we require the probability of absorption under $\alpha^*$ to be bounded away from 0,
in Lemma~\ref{lemma:heavy} we do not have this restriction.

The proof of Proposition~\ref{prop:absorbing-eq} will involve a repeated use of Lemma~\ref{lemma:heavy}, which will generate a sequence of 0-equilibria in truncations of $\calT$.
Under each of these 0-equilibria, the termination probability is positive,
and the sum of these strategy profiles will have a total probability of termination in a single round which will be bounded away from 0.

\bigskip

\begin{proof}[Proof of Lemma~\ref{lemma:heavy}]
The proof extends the construction for quitting games, which was explained in Section~\ref{sec:ss:nonQ}.

\medskip
\noindent\textbf{Step 1: The stage payoff $q$ and the auxiliary game.}

%Eq.~\eqref{equ:871} implies in particular that 
%$y^i_\calT = y^i$ for every $i \in I$, and hence
%$I_\calT^{\mathrm{sp}} = I^{\mathrm{sp}}$
%and $R^{\mathrm{sp}}_\calT = R^{\mathrm{sp}}$.
%In particular, 
By (T4),
$R^{\mathrm{sp}}_\calT$ is not $\ep$-close to any Q-matrix.
By Lemma~\ref{lemma:Qmatrix}, there is a vector $q^*\in[-1,1]^{I^{\mathrm{sp}}_\calT}$ such that the
approximate problem \eqref{equ:lcp:approx} for $(R^{\mathrm{sp}}_\calT,q^*)$ has no solution.
Extend $q^*$ to a vector $q\in[-1,1]^n$ by
\[
q^i :=
\begin{cases}
q^{*,i}, & i\in I^{\mathrm{sp}}_\calT,\\
1, & i\notin I^{\mathrm{sp}}_\calT.
\end{cases}
\]

For every $\lambda\in(0,1]$,
every node $s \in S_0$,
and every strategy profile $\alpha$, 
denote by $\gamma_{q,\lambda,s}(\alpha)$ the $\lambda$-discounted payoff in $\calT$ under $\alpha$,
when the initial node is $s$ (rather than the root $r$)
and the stage payoff 
as long as no player stops is $q$: 
\[
\gamma_{q,\lambda,s}(\alpha)
:= \E_{s,\alpha}\left[(1-\lambda)^{\theta}X_{J_\theta,\theta}
   \bfone_{\{\theta<\infty\}}
   +\bigl(1-(1-\lambda)^{\theta}\bfone_{\{\theta<\infty\}}\bigr)q\right].
\]
Recall that when the initial node is $s \neq r$, except for the first round of the game,
all subsequent rounds start at $r$.

\medskip
\noindent\textbf{Step 2: A limit of discounted equilibria.}

By \cite{Fink64}, for every $\lambda\in(0,1]$ the game $\calT$
admits a stationary $\lambda$-discounted equilibrium $\alpha^*_\lambda$,
that is, 
\[
\gamma^i_{q,\lambda,s}(\alpha^*_\lambda)
\ge
\gamma^i_{q,\lambda,s}(\alpha^i,\alpha^{*,-i}_\lambda),
\qquad\forall i\in I,\ \forall s\in S_0,\ \forall \alpha^i.
\]
By \cite{Neyman03}, the selection
$\lambda\mapsto\alpha^*_\lambda$ may be taken to be semi-algebraic,
so that the limits
\[
\alpha^* := \lim_{\lambda\to0}\alpha^*_\lambda,
\qquad
\hbox{and}
\qquad
\gamma^*_s := \lim_{\lambda\to0}\gamma_{q,\lambda,s}(\alpha^*_\lambda),
\qquad \forall s\in S_0,
\]
exist.

\medskip
\noindent\textbf{Step 3: $H_{\delta^2}(\alpha^*)\neq\varnothing$.}

Assume by contradiction that $H_{\delta^2}(\alpha^*)=\varnothing$; that is,
\begin{equation}
\label{equ:proof:light}
\prob_{s,\alpha^*}(\theta<\theta_{\mathrm{fr}}) < \delta^2,
\qquad\forall s\in S_0.
\end{equation}
We will derive an approximate solution to the linear complementarity problem \eqref{equ:lcp:approx} for
$(R^{\mathrm{sp}}_\calT,q^*)$, contradicting the choice of $q^*$.
Two consequences of Eq.~\eqref{equ:proof:light} that we will use are:
\begin{itemize}
\item[P1)] 
At $s$, each player stops with low probability:
$\alpha^{*,i}(s) < \delta^2$ for every $s \in S_0$ and every $i \in I$.
\item[P2)]  
The payoff under $\alpha^*$ is hardly affected by the initial node:
$|\gamma^*_s - \gamma^*_{s'}| < 2\delta^2$,
for every $s,s' \in S_0$.
Indeed,
if the play is not terminated in the first round when the initial node is $s$ (or $s'$) (an event with probability at least $1-\delta^2$),
then the process repeats itself from node $r$ and the limit payoff is $\gamma^*_r$.
\end{itemize}
Since $|S|<\infty$,
Property (P1) implies that for every $\lambda > 0$ sufficiently small,
\begin{equation}
\label{lemma:alpha:1}
\alpha^{*,i}_\lambda(s) < \delta^2, \ \ \ \forall s \in S_0, \forall i\in I.
\end{equation}

\bigskip
\noindent\textbf{Step 3.1: Simultaneous stops are negligible.}

Fix for a moment a node $s \in S_0$.
By \eqref{lemma:alpha:1},
for every $J \subseteq I$ with $|J| \geq 2$, every $i \in J$, and every $\lambda > 0$ sufficiently small,
\begin{equation}
\label{equ:a:1}
\prod_{j \in J} \alpha^{*,j}_\lambda(s) \prod_{j \not\in J} (1-\alpha^{*,j}_\lambda(s)) \leq \frac{\delta^2}{1-\delta^2} \alpha^{*,i}_\lambda(s)\prod_{j \neq i} (1-\alpha^{*,j}_\lambda(s)). 
\end{equation}
That is, 
the probability that $J$ is the set of players who stop at $s$
is significantly smaller than the probability that player~$i$ stops alone.
Summing \eqref{equ:a:1} over all subsets $J$ of $I$ 
with $|J| \geq 2$,
we obtain that for every $\lambda > 0$ sufficiently small,
\begin{equation}
\label{equ:a:2}
\prob_{r,\alpha^*_\lambda}\bigl(s_\theta = s, |J_\theta| \geq 2\bigr) \leq \frac{\delta^2 2^n}{1-\delta^2}\prob_{r,\alpha^*_\lambda}\bigl(s_\theta=s, |J_\theta| = 1\bigr). 
\end{equation}
Summing \eqref{equ:a:2} over all nodes $s \in S_0$,
the discounted probability of a simultaneous stop is small:
\begin{equation}
\label{equ:proof:multi}
z_{\ge2} := \lim_{\lambda\to0}\E_{r,\alpha^*_\lambda}
   \bigl[(1-\lambda)^\theta\bfone_{\{|J_\theta|\ge2\}}\bigr]
\le\frac{\frac{\delta^2 2^n}{1-\delta^2}}{1 + \frac{\delta^2 2^n}{1-\delta^2}} 
\leq \delta^2 2^n.
\end{equation}

\bigskip
\noindent\textbf{Step 3.2: Each player~$i$ may stop only when $X^i_{\{i\},s} = y^i_\calT$.}

Let $i \in I$ be arbitrary, and let $s \in S_0$ be a node such that $X^i_{\{i\},s} = y^i_\calT$. By (P1),
\[ \gamma^{*,i}_s = 
\lim_{\lambda \to 0} \gamma^i_{q,\lambda,s}(\alpha^{*}_\lambda)
\geq \lim_{\lambda \to 0} \gamma^i_{q,\lambda,s}(1^i,\alpha^{*,-i}_\lambda) \geq y^i_\calT - 2\delta^2. \]
By (P2), for every $s' \in S_0$,
\begin{equation}
\label{equ:874}
\lim_{\lambda \to 0} \gamma^i_{q,\lambda,s'}(\alpha^{*}_\lambda) = 
\gamma^{*,i}_{s'} \geq \gamma^{*,i}_s -2\delta^2 \geq y^i_\calT - 4\delta^2. 
\end{equation}
Let $s' \in S_0$ be a node such that $X^i_{\{i\},s'} < y^i_\calT$. Since payoffs are multiple of $\ep$, we have  $X^i_{\{i\},s'} \leq y^i_\calT- \ep$.
Hence, by (P1), since $6\delta^2 < \ep$, and by \eqref{equ:874}, 
\[ \lim_{\lambda \to 0} \gamma^i_{q,\lambda,s'}(1^i,\alpha^{*,-i}_\lambda) \leq X^i_{\{i\},s'} + 2\delta^2
\leq y^i_\calT- \ep + 2\delta^2 < \lim_{\lambda \to 0} \gamma^i_{q,\lambda,s'}(\alpha^{*}_\lambda).
\]
Hence, $\alpha^{*,i}_{\lambda}(s') = 0$ for every $\lambda > 0$ sufficiently small: player~$i$ stops with probability 0 at such nodes $s'$. 
Therefore, if $\alpha^{*,i}(s) > 0$, then $X^i_{\{i\},s} = y^i_\calT$.

\bigskip
\noindent\textbf{Step 3.3: Free-riders do not stop under $\alpha^*_\lambda$.}

Fix a free-rider $i \not\in I^{\mathrm{sp}}_\calT$.
By \eqref{equ:proof:multi},
since $q^i=1 > \min_{j \neq i} Z_{\calT,j}^i$,
since by Step~3.1 simultaneous stops are negligible,
and since by Step 3.2 each player $j$ stops only at nodes where her payoff is $y^j_\calT = Z_{\calT,j}^j$,
when the other players follow $\alpha^{*,-j}_\lambda$,
never stopping yields player~$i$ at least 
$\min_{j \neq i} Z_{\calT,j}^i - \delta^22^n/(1-\delta^2)$:
\[ \lim_{\lambda \to 0} \gamma_{q,\lambda,s}^i(0^i,\alpha^{*,-i}_\lambda) \geq \min_{j \neq i} Z_{\calT,j}^i - \frac{\delta^22^n}{1-\delta^2} \geq Z_{\calT,i}^i + \ep  - \frac{\delta^22^n}{1-\delta^2}, \]
where the last inequality holds since payoffs are an integer multiple of $\ep$ and since $i$ is a free-rider.
By \eqref{equ:proof:light},
stopping at some node $s$ yields 
player~$i$ at most $y^i_\calT + 2\delta^2$:
\[ \lim_{\lambda \to 0} \gamma_{q,\lambda,s}^i(1^i,\alpha^{*,-i}_\lambda) \leq y^i_\calT + 2\delta^2 = Z_{\calT,i}^i + 2\delta^2. \]
By the choice of $\delta$, for every $\lambda > 0$ sufficiently small,
\[ \lim_{\lambda \to 0} \gamma_{q,\lambda,s}^i(0^i,\alpha^{*,-i}_\lambda) > \lim_{\lambda \to 0} \gamma_{q,\lambda,s}^i(1^i,\alpha^{*,-i}_\lambda). \]
This implies that at any node, stopping cannot be part of player~$i$'s best response to $\alpha^{*,-i}_\lambda$,
provided $\lambda > 0$ is sufficiently small.

\bigskip
\noindent\textbf{Step 3.4: Construction of $(w,z)$.}

For every $i\in I$ let
\[
z_i := 
\lim_{\lambda\to0}\E_{r,\alpha^*_\lambda}
   \bigl[(1-\lambda)^\theta\bfone_{\{J_\theta=\{i\}\}}\bigr]
\]
be the limit discounted probability that player~$i$ stops alone, and set
$z_0:=1-\sum_{i\in I^{\mathrm{sp}}_\calT}z_i$.
By Step~3.3, $z:=(z_0,(z_i)_{i\in I^{\mathrm{sp}}_\calT})\in\Delta(\{0\}\cup I^{\mathrm{sp}}_\calT)$.

Define $w\in\dR^{I^{\mathrm{sp}}_\calT}$ by
\[
w := z_0q^* + \sum_{i\in I^{\mathrm{sp}}_\calT} z_i(R^{\mathrm{sp}}_\calT)_{\bullet i},
\]
where $(R^{\mathrm{sp}}_\calT)_{\bullet i}$ is the $i$-th column of $R^{\mathrm{sp}}_\calT$.
The discounted weight of simultaneous stops was added to $z_0$.
By Step~3.1, this discounted weight is at most $2^n\delta^2$,
and hence
\begin{equation}
\label{equ:proof:d:err}
\bigl\|w-\gamma^{*}_r\bigr\|_\infty \le 2^{n+1}\delta^2.
\end{equation}

\bigskip
\noindent\textbf{Step 3.5: $(w,z)$ solves \eqref{equ:lcp:approx}.}

Since $w$ is a convex combination of vectors in $[-1,1]^{|I^{\mathrm{sp}}_\calT|}$, 
we have $w \in [-1,1]^{|I^{\mathrm{sp}}_\calT|}$.
Hence, the first condition of \eqref{equ:lcp:approx} holds.

Fix $i\in I^{\mathrm{sp}}_\calT$ and a node $s$ with $X^i_{\{i\},s}=y^i_\calT$.
By \eqref{equ:proof:light} the probability that some
other player(s) stop at $s$ is below $\delta^2$.
Since $\alpha^*_\lambda$ is a $\lambda$-discounted equilibrium, for every $\lambda$ sufficiently small,
\[
\gamma^i_{q,\lambda,s}(\alpha^*_\lambda)
\ge\gamma^i_{q,\lambda,s}(1^i,\alpha^{*,-i}_\lambda)
\ge y^i_\calT-2\delta^2.
\]
Taking the limit as $\lambda$ goes to 0, and using \eqref{equ:proof:d:err} and (P2), we obtain that
$w_i\ge y^i_\calT-2^{n+2}\delta^2\ge y^i_\calT-\delta$, which is the second condition of
\eqref{equ:lcp:approx}.

For complementarity, suppose that $z_i>0$.
Then under $\alpha^*_\lambda$, player~$i$ stops with positive probability, provided $\lambda$ is sufficiently small.%
\footnote{Note that this does not imply that under the limit strategy $\alpha^{*,i}$ player~$i$ stops with positive probability.}
Hence, under $\alpha^*_\lambda$ stopping is
asymptotically optimal for player~$i$ at some node $s \in S_0$.
By Step 3.2, 
at any node $s$ in which she stops with positive probability,
$|\gamma^{*,i}_s - y^i_\calT| < 2^{n+1}\delta^2/(1-\delta^2)$.
By \eqref{equ:proof:d:err} and the triangle inequality, $|w_i-y^i_\calT|\le 2^{n+2}\delta^2/(1-\delta^2) \leq \delta$.
Thus, if $w_i>y^i_\calT+\delta$, then $z_i=0$.
This is the third condition of \eqref{equ:lcp:approx}.

Thus, $(w,z)$ solves \eqref{equ:lcp:approx} for $(R^{\mathrm{sp}}_\calT,q^*)$, contradicting the
choice of $q^*$.
Therefore $H_{\delta^2}(\alpha^*)\neq\varnothing$.

%Since $H_{\delta^2}(\alpha^*)\neq\varnothing$,
%under $\alpha^*$ at least one player stops with positive probability at some node.
%In Step~4 we should ????????

\bigskip
\noindent\textbf{Step 4: If under $\alpha^*$
at least two players stop with positive probability,
then A2 holds.}

When two players stop with positive probability under $\alpha^*$,
the stationary strategy profile 
$\alpha^*$ is a
$0$-equilibrium of $\calT$ from every initial node; see \cite[Case~A in Section~3.3]{VT89}.

By Step~3, A2(i) holds.

To prove that A2(ii) holds,
suppose that $\gamma^i(\alpha^*) \leq y^i_\calT - 5\sqrt\ep$ but $\pi(0^i,\alpha^{*,-i}) < \ep$.
Let $\alpha^i$ be the strategy of player~$i$ in $\calT$ that stops with probability 1 at every node $s$ such that $X^i_{\{i\},s} = y^i_\calT$.
Then 
\[ y^i_\calT - 5\sqrt\ep \geq \gamma^i(\alpha^*) \geq 
\gamma^i(\alpha^i,\alpha^{*,-i}) 
\geq y^i_\calT - 4\ep, \]
where the last inequality holds because $\pi(0^i,\alpha^{*,-i}) < \ep$ and by \eqref{equ:871}.
This is a contradiction.

We turn to prove that A2(iii) holds.
Suppose then that for some player $i \in I$ we have
$\gamma^i(\alpha^{*}) < 0$.
%This last inequality implies that the conditional probability that player $i$ stops alone is at least $1-\delta$.
%In particular, there can be at most one $i \in I$ satisfying the hypotheses of A2(iii).
%
% Since $\alpha^*$ is a 0-equilibrium at all initial nodes,
% and since by the genericity assumption $y^i_\calT > X_{\{i\},s}^i$ whenever $X_{\{i\},s} \neq Z_i$,
% player~$i$ stops with positive probability only when 
% at nodes where $X_{\{i\},s} = Z_i$.
% Hence,
% \[ \|\gamma(\alpha^*) - Z_i\|_\infty \leq 2\delta. \]
%
The set of all nodes at which at least one player other than $i$ stops with positive probability is
\[ D := \left\{ s \in S_0 \colon \sum_{j \neq i} \alpha^{*,j}(s) > 0\right\}, \]
which is nonempty since at least two players stop with positive probability under $\alpha^*$.
We argue that 
\begin{equation}
\label{equ:min:1}
\min_{s \in D} \gamma^{*,i}_s \leq \gamma^i(\alpha^*). 
\end{equation}
Indeed, let $\alpha^i$ be the (not necessarily stationary) strategy
 of player~$i$
that stops with probability 0 before the first visit to $D$,
and coincides with $\alpha^{*,i}$ once a node in $D$ is reached.
Since $\alpha^*$ is a 0-equilibrium,
$\gamma^i(\alpha^*) \geq \gamma^i(\alpha^i,\alpha^{*,-i})$.
Since under $(\alpha^i,\alpha^{*,-i})$ no player stops before reaching $D$,
$\gamma^i(\alpha^i,\alpha^{*,-i})$ is a weighted average of $(\gamma^{*,i}_s)_{s \in D}$.
Eq.~\eqref{equ:min:1} followed.

We finally argue that for any $s \in D$ that attains the minimum in \eqref{equ:min:1},
the conclusion of A2(iii) holds with 
$\beta^j = \alpha^{*,j}(s)$ for $j \neq i$.
To see that $A \leq \gamma^i(\alpha^*)$,
note that $A$ is the payoff when the initial node is $s$,
players $I \setminus \{i\}$ follow $\alpha^*$,
and player $i$ stops with probability 1.
Since $\alpha^*$ is a 0-equilibrium when the initial node is $s$, 
$A \leq \gamma^{*,i}_s \leq \gamma^i(\alpha^*)$.

To see that $B \leq \gamma^i(\alpha^*)$,
let $p:=1-\prod_{j\neq i}(1-\beta^j)>0$
be the probability that at least one player in $I \setminus \{i\}$ stops at $s$ under $\alpha^*$.
Starting from $s$, consider the deviation
in which player~$i$ continues at $s$ and,
if no opponent stops, continues until
the next visit to $D$, after which she
follows $\alpha^{*,i}$.
By the minimality of $s$, this deviation
yields at least
$pB+(1-p)\gamma_s^{*,i}$.
Since $\alpha^*$ is a 0-equilibrium,
$\gamma_s^{*,i}
\geq pB+(1-p)\gamma_s^{*,i}$,
and therefore
$B\leq\gamma_s^{*,i}
\leq\gamma^i(\alpha^*)$.

\bigskip
\noindent\textbf{Step 5: If under $\alpha^*$
one player stops with positive probability, then Conclusion A1 holds.}

Denote by $i_0$ the only player who stops with positive probability under $\alpha^*$. 
Since $\lim_{\lambda \to 0} \alpha^{*,i}_\lambda = 0^i$ for every $i \neq i_0$,
the limit probability of simultaneous stopping is 0.
In particular, $z_{i_0} = 1$.

We argue that 
$\gamma(\alpha^{*,i_0},\mathbf 0^{-i_0}) = \gamma(\alpha^*) = Z_{\calT,i_0}$.
Indeed, let $\alpha^{i_0}$ be the strategy of player $i_0$ that stops with probability 1 at the first node $s$ such that $X_{\{i_0\},s} = Z_{\calT,i_0}$.
Then
\[ \gamma^{i_0}(\alpha^*) = \lim_{\lambda \to 0} \gamma^{i_0}_{q,\lambda,r}(\alpha^{*}_\lambda)
\geq \lim_{\lambda \to 0} \gamma^{i_0}_{q,\lambda,r}(\alpha^{i_0},\alpha^{*,-i_0}_\lambda) = Z_{\calT,i_0}^{i_0}. \]
On the other hand, since player~$i_0$ is the only player whose limit stopping probability is positive,
$\gamma^{i_0}(\alpha^*) \leq y^{i_0}_\calT = Z_{\calT,i_0}^{i_0}$.
Property (T3) ensures that 
$\gamma(\alpha^*) = Z_{\calT,i_0}$.

Similarly, $\gamma(s,\alpha^*)$,
the payoff under $\alpha^*$ when the initial state is $s$,
is $Z_{\calT,i_0}$, for every $s \in S_0$.

Since under $\alpha^*$ only player $i_0$ stops with positive probability,
since by Property (T3) $X^{i_0}_{\{i_0\},s} < Z^{i_0}_{\calT,i_0}$ whenever $X_{\{i_0\},s} \neq Z_{\calT,i_0}$,
and since $\gamma(\alpha^*) = Z_{\calT,i_0}$,
it follows that player~$i_0$ stops with positive probability only at nodes where $X_{\{i_0\},s} = Z_{\calT,i_0}$.
In particular, at all other nodes, no player stops.

Hence, for every player $i \neq i_0$ and every node $s \in S_0$ such that $X_{\{i_0\},s} \neq Z_{\calT,i_0}$,
\[ Z_{\calT,\{i_0\}}^i = \gamma^i(s,\alpha^*)
= \lim_{\lambda \to 0} \gamma^i_{q,\lambda,s}(\alpha^*_\lambda)
\geq \lim_{\lambda \to 0} \gamma^i_{q,\lambda,s}(1^i,\alpha^{*,-i}_\lambda)
= X_{\{i\},s}^i,
\]
and A1(iii) holds.
Similarly,
at every node $s \in S_0$ such that $X_{\{i_0\},s} = Z_{\calT,i_0}$,
\[ Z_{\calT,\{i_0\}}^i = \gamma^i(s,\alpha^*)
= \lim_{\lambda \to 0} \gamma^i_{q,\lambda,s}(\alpha^*_\lambda)
\geq \lim_{\lambda \to 0} \gamma^i_{q,\lambda,s}(1^i,\alpha^{*,-i}_\lambda)
= (1-\alpha^{*,i_0}(s)) X^i_{\{i\},s}
+ \alpha^{*,i_0}(s)X^i_{\{i,i_0\},s}.
\]
Since under $\alpha^*$ the play terminates with probability 1,
there is such a node $s \in S_0$ for which $\alpha^{*,i_0}(s) > 0$.
For this node $s$, A1(ii) holds with $\beta = \alpha^{*,i_0}(s)$.

We turn to prove that A1(i) holds as well.
Indeed, otherwise, $i_0$ is a free-rider and $X_{\{i\},s}^{i_0} > y^{i_0}_\calT$ for every $s \in S_0$ and every $i \neq i_0$.
Since $i_0$ is a free-rider,
$Z_{\calT,i_0}^{i_0} < 0$
%$Z_{\calT,i_0}^{i_0} < Z_{\calT,i}^{i_0}$ for every $i \neq i_0$,
and $q^{i_0} = 1$.
Since $\lim_{\lambda \to 0} \alpha^{*,-i_0}_\lambda = \mathbf 0^{-i_0}$,
the discounted probability of simultaneous stopping by at least two opponents under \((0^{i_0},\alpha_\lambda^{*,-i_0})\) tends to zero.
Hence, since $(\alpha^*_\lambda)_{\lambda > 0}$ are $\lambda$-discounted equilibria,
\[ Z^{i_0}_{\calT,i_0} =\lim_{\lambda \to 0} \gamma^{i_0}_{q,\lambda,r}(\alpha_\lambda^{*}) \geq \lim_{\lambda \to 0} \gamma^{i_0}_{q,\lambda,r}(0^{i_0},\alpha_\lambda^{*,-i_0})
\in \conv\{ q^{i_0},X_{\{i\},s}^{i_0}, s \in S_0, i \neq i_0\} > 
y^{i_0}_\calT = Z^{i_0}_{\calT,i_0},\]
a contradiction.
\end{proof}

%% ---------------------------------------------------------------------
\subsection{Strongly Orthogonal Strategy Profiles and Truncations}
\label{subsec:orthogonal}
%% ---------------------------------------------------------------------

\begin{definition}[Strong orthogonality]
\label{def:strong-orthogonal}
Let $\eta > 0$.
A finite sequence $(\alpha_1,\dots,\alpha_m)$ of stationary strategy profiles is
\emph{$\eta$-strongly orthogonal} if for every $1\le \ell\le m-1$, every
node $s\in H_\eta(\alpha_1\dot+\cdots\dot+\alpha_\ell)$, and every
node $s'$ equal to $s$ or a descendant of $s$, one has
$\alpha_{\ell+1}(s')=\mathbf 0$; that is, $\alpha_{\ell+1}$ has every player
continue from each $\eta$-heavy node of the partial sum onward.
\end{definition}

By~\eqref{fact:heavy}, $(\alpha_1,\dots,\alpha_m)$ is an $\eta$-strongly orthogonal sequence of stationary strategy profiles,
then any subsequence $(\alpha_{k_1},\dots,\alpha_{k_d})$ is $\eta$-strongly orthogonal as well.

A \emph{truncation} of a game on a tree $\calT$, is the game on a tree that is obtained by removing all strict
descendants of a set of nodes and keeping the remaining data fixed.

\begin{definition}[Truncation]
\label{def:subgame}
Let $\calT = (I, S, S_1, r, (C_s), (p_s), (X_s))$ be a game on a tree.
A game on a tree $\calT' = (I, S', S_1', r, (C'_s), (p'_s), (X'_s))$ is a \emph{truncation} of
$\calT$ if $S'\subseteq S$ and for every $s\in S_0':=S'\setminus S_1'$ we have
(a) $C'_s=C_s$, (b) $p'_s=p_s$, and (c) $X'_s=X_s$.
For $D\subseteq S_0$ we denote by $\calT_D$ the truncation of $\calT$ obtained by
deleting every strict descendant of every node of $D$; the nodes of $D$ become
leaves of $\calT_D$.
\end{definition}

%For a node $s\in S_0$ let $\nu_s$ be the probability that
%when no player stops,
%the node $s$ is visited during the first round.

The next result connects equilibria in truncations and sequences of strongly orthogonal stationary strategy profiles.
It shows that if the expected payoff under
every stationary strategy profile in a strongly orthogonal sequence is close to some vector $a \in [-1,1]^n$,
then so is the expected payoff under their sum.

\begin{lemma}[Payoff of a strongly orthogonal sum of equilibria]
\label{lem:strong-orth-payoff}
Let $\delta\in(0,1/6)$,
let $a=(a^i)_{i\in I}\in[-1,1]^n$,
and let $(\alpha_1,\ldots,\alpha_m)$ be a
$\delta^2$-strongly orthogonal sequence
of stationary $\delta$-equilibria of
$\calT$ such that
\[
\gamma(\alpha_\ell;\calT)
\in\prod_{i\in I}[a^i,a^i+\delta],
\qquad \forall \ell=1,\ldots,m.
\]
Denote $\alpha=\alpha_1\dot+\cdots\dot+\alpha_m$.
Then:
\begin{enumerate}[(a)]
\item $a^i-\delta\le\gamma^i(\alpha)$ for every $i\in I$;
\item for every $i\in I$ and every stationary strategy $x'^i$ of
  player~$i$, $\gamma^i(x'^i,\alpha^{-i})\le a^i+8\delta$.
\end{enumerate}
\end{lemma}

Lemma~\ref{lem:strong-orth-payoff} was proved by \cite{SS04} in the context of two-player games
(see \cite[Lemmas~8.13 and~8.14]{SS04}).%
\footnote{\cite[Lemma~8.14]{SS04} does not indicate that $\alpha_\ell$ is an $\ep$-equilibrium of the truncation $\calT_{H_{\ell-1}}$. 
This seems to be an omission in the statement of \cite[Lemma~8.14]{SS04}.}
Their arguments are independent of the number of players,
hence their proof applies to our setup.

The next result compares the probability to visit a node in a certain subset of nodes during the first round in $\calT$ and in a truncation of $\calT$.

Given a truncation $\calT'$ of $\calT$, 
for every $M \subseteq S_0$
denote by $M \cap \calT'$ the nodes in $M$ that lie in $\calT'$ but are not leaves in $\calT'$.

\begin{lemma}
\label{lem:max-nodes-survive}
Let $D\subseteq S_0$ and let $\calT_D$ be the truncation obtained by
cutting the tree at $D$.  
Then, for every subset $M \subseteq S_0$,
\begin{equation}
\label{eq:mu-truncation}
 \prob_{r,\mathbf 0}(M \cap \calT_D)
 \geq
 \prob_{r,\mathbf 0}(M)
 -
 \prob_{r,\mathbf0}(D).
\end{equation}
\end{lemma}

\begin{proof}
If the play visits $M$ in the first round of $\calT$ but fails to visit $M$ in the first round of
$\calT_D$, then the play must first pass through a node of $D$ in the first round.
Eq.~\eqref{eq:mu-truncation} follows by 
calculating the probability of each of these events.
\end{proof}

%%%%%%%%%%%%%
\subsection{Proof of Proposition~\ref{prop:absorbing-eq}}
\label{proof:prop}

We here prove Proposition~\ref{prop:absorbing-eq}.
Fix $\ep>0$ sufficiently small. Choose
$\delta\in \bigl(0,\min\{\ep^3/4,\delta_{\LCP}(\ep)/16,2^{-n-4}\}\bigr)$,
such that $1/\delta$ is an integer,
and set
\[
  \eta:=6\sqrt\ep,\qquad
  \calA:=([-1,1]\cap\delta\dZ)^I,\qquad
  L:=|\calA|=(2/\delta+1)^n,\qquad
  d:=\frac{\delta^4}{2L}.
\]

\medskip
\noindent\textbf{Step 1: Equilibria in truncations.}

Let $D\subseteq S_0$ such that 
$y^i_{\calT_D} = y^i_\calT$ for every $i \in I$,
and suppose that $\beta^*$ is a stationary
$0$-equilibrium of $\calT_D$ with positive termination probability.
Denote by $\alpha^*$ its extension to $\calT$, such that no player stops
at nodes in $D$ and at their descendants. 
Then
\begin{equation}\label{eq:short:invariance}
  \pi(\alpha^*;\calT)=\pi(\beta^*;\calT_D),\qquad
  \gamma(\alpha^*;\calT)=\gamma(\beta^*;\calT_D).
\end{equation}
We argue that for every player $i \in I$ and every deviation $\alpha^i$ of player $i$ in $\calT$, 
\begin{equation}\label{eq:short:barrier}
  \gamma^i(\alpha^i,\alpha^{*,-i};\calT)
  \leq\max\{\gamma^i(\alpha^*;\calT),y^i_\calT\}.
\end{equation}

We will divide the proof of \eqref{eq:short:barrier} to two cases,
according to whether $\beta^{*,-i} = \mathbf 0^{-i}$.

If $\beta^{*,-i} = \mathbf 0^{-i}$,
then, $\gamma^i(\alpha^i,\alpha^{*,-i};\calT) \leq \max\{0,y^i_\calT\}$ for every strategy $\alpha^i$ in $\calT$,
with equality when $\alpha^i = \beta^{*,i}$.
Eq.~\eqref{eq:short:barrier} follows.

Assume next that $\beta^{*,-i} \neq \mathbf 0^{-i}$.
Write 
\[ v:=\gamma^i(\beta^*;\calT_D) = \gamma^i(\alpha^*;\calT)
\qquad \hbox{and} \qquad
w:=\max\{v,y^i_\calT\}. \] 
Given a real number $b$,
consider the version of $\calT_D$ where,
if no player stops during the first round,
the game terminates at the end of the first round with terminal payoff $b$ for player~$i$.
The terminal payoff of the other players is irrelevant.
The maximal payoff for player $i$ in this variant is 
\[ F_{\calT_D}(b) := \max_{\beta^i} \left( \rho^i(\beta^i,\beta^{*,-i};\calT_D) + (1-\pi(\beta^i,\beta^{*,-i};\calT_D))b\right), \]
where $\beta^i$ ranges over all strategies of player~$i$ in $\calT_D$.
Since $\beta^*$ is a 0-equilibrium of $\calT_D$, $F_{\calT_D}(v) =v$.
Since the function $b \mapsto F_{\calT_D}(b)$ is monotone and $1$-Lipschitz,
and since $w \geq v$, we have
$F_{\calT_D}(w) \leq w$.

At nodes that lie in $\calT$ but not in $\calT_D$, 
no player in $I \setminus \{i\}$ stops,
and hence the highest payoff that player~$i$ can obtain by stopping in those nodes is $y^i_\calT$.
Since $w \geq y^i_\calT$, we have $F_\calT(w) \leq w$,
where 
\[ F_{\calT}(b) := \max_{\alpha^i} \left( \rho^i(\alpha^i,\alpha^{*,-i};\calT) + (1-\pi(\alpha^i,\alpha^{*,-i};\calT))b\right)
\]
and $\alpha^i$ ranges over all strategies of player~$i$ in $\calT$.
Since $\beta^{*,-i} \neq \mathbf 0^{-i}$,
we have
$\pi(\alpha^i,\alpha^{*,-i};\calT) > 0$ for every strategy $\alpha^i$ of player $i$ in $\calT$.
Hence, the iterations $((F_\calT)^\ell(w))_{\ell \in \dN}$ converge to $\sup_{\alpha^i}\gamma^i(\alpha^i,\alpha^{*,-i};\calT)$ as $\ell$ goes to infinity.
Eq.~\eqref{eq:short:barrier} follows.

\bigskip

We will also use the following estimate:
If
$\prob_{r,\mathbf0}(D;\calT)\leq h$ and
$\pi(0^i,\beta^{*,-i};\calT_D)\geq p>0$, then under every
deviation $\alpha^i$ of player $i$ in $\calT$, the probability of encountering $D$ before termination is
at most $h/p$. 
Therefore,
\begin{equation}\label{eq:short:transfer}
  \gamma^i(\alpha^i,\alpha^{*,-i};\calT)
  \leq\gamma^i(\alpha^*;\calT)+2h/p.
\end{equation}

\medskip
\noindent\textbf{Step 2: An inductive construction of terminating $0$-equilibria.}

We next construct a sequence of terminating 0-equilibria in a decreasing sequence of truncations of $\calT$
(decreasing w.r.t.~inclusion of the sets of nodes of the truncations).

Set $H_0:=\varnothing$ and $\sigma_0:=\mathbf0$.
Suppose that $H_k$ has been constructed such that 
$\prob_{r,\mathbf0}(H_k)\leq\delta^2$.
By (T6), for every $g\in\range(X)$,
\[
 \prob_{r,\mathbf0}\bigl(
  g\text{ appears at a non-leaf node of }\calT_{H_k}
 \bigr)\geq 1-\delta^2>\tfrac12.
\]
Thus, the range of $X$ on $\calT_{H_k}$ is $\range(X_\calT)$.
In particular, $R_{\calT_{H_k}} = R_\calT$.
Therefore, $R^{\mathrm{sp}}_{\calT_{H_k}} = R^{\mathrm{sp}}_\calT$,
and hence $R^{\mathrm{sp}}_{\calT_{H_k}}$ is not $\ep$-close to any Q-matrix.
Moreover, Conclusion~A1 does not hold for $\calT_{H_k}$,
and hence Conclusion~A2 holds for $\calT_{H_k}$:
there is a stationary $0$-equilibrium $\beta^*_{k+1}$ of $\calT_{H_k}$ that satisfies Conclusion~A2.

Extend $\beta_{k+1}^*$ to a
profile $\alpha_{k+1}^*$ on $\calT$
by letting no player stop at nodes in $H_k$ and their descendants,  and set
\[
 \sigma_{k+1}:=\alpha_1^*\dot+\cdots\dot+\alpha_{k+1}^*,
 \qquad H_{k+1}:=H_{\delta^2}(\sigma_{k+1};\calT).
\]
By A2(i), $H_{k+1} \supsetneq H_k$. 
We continue the inductive construction until the 
first index $N\leq |S_0|$ such that
$\prob_{r,\mathbf0}(H_N)>\delta^2$.

The stationary strategy profiles $(\alpha_1^*,\ldots,\alpha_N^*)$ form a
$\delta^2$-strongly orthogonal sequence, and
\eqref{eq:short:invariance} applies to each pair
$(\beta_k^*,\alpha_k^*)$.

Every $\alpha_k^*$ is an $\eta$-equilibrium of $\calT$.
Indeed, fix $i \in I$ and set
$v_k^i:=\gamma^i(\beta_k^*;\calT_{H_{k-1}})$.
If $v_k^i\geq y^i_\calT-\eta$, 
then by \eqref{eq:short:barrier} player~$i$ cannot profit more than $\eta$ by deviating.
Otherwise, by A2(ii) we have
$\pi(0^i,\beta_k^{*,-i};\calT_{H_{k-1}})\geq\ep$, and hence by
\eqref{eq:short:transfer} 
an upper bound on the amount player~$i$ can gain by deviating is
$2\delta^2/\ep<\eta$.

Choose $a_k\in\calA$ such that
\[
 \gamma(\alpha_k^*;\calT)
 \in\prod_{i\in I}\bigl[a_k^i,a_k^i+\delta\bigr].
\]

\medskip
\noindent\textbf{Step 3: A low-payoff component.}

Suppose that for some $k \in \{1,\dots,N\}$ and $i \in I$ we have
$a_k^i\leq y^i_\calT-\eta$.
Then
$\gamma^i(\beta_k^*;\calT_{H_{k-1}})<y^i_\calT-5\sqrt\ep$.
By A2(ii),
\[
 \pi(\alpha_k^*;\calT)\geq\pi(0^i,\alpha_k^{*,-i};\calT)
 =\pi(0^i,\beta_k^{*,-i};\calT_{H_{k-1}})\geq\ep>d.
\]
Set $\alpha^*:=\alpha_k^*$.
As shown in Step~2, $\alpha^*$ is an $\eta$-equilibrium. 
If any player $j \in I$
satisfies $\gamma^j(\alpha^*;\calT)<y^j_\calT-7\sqrt\ep$,
A2(ii) yields $\pi(0^j,\alpha^{*,-j};\calT)\geq\ep>d^2$.
If $\gamma^j(\alpha^*;\calT)\leq-7\sqrt\ep$, then,
by A2(iii), applied to $\beta_k^*$,
we have $\max\{A,B\}\leq\gamma^j(\beta_k^*;\calT_{H_{k-1}})
=\gamma^j(\alpha^*;\calT)$.
Finally, since $\pi(\alpha^*;\calT)\geq\ep>2d^2$,
there is at most one player with
$\pi(0^j,\alpha^{*,-j};\calT)<d^2$.
Thus all three conclusions hold.

\medskip
\noindent\textbf{Step 4: A common payoff.}
From now on we assume that
$a_k^i>y^i_\calT-\eta$ for every $i \in I$ and every $k=1,\dots,N$.
Since $\prob_{r,\mathbf0}(H_N)\geq \delta^2$,
\[
 \sum_{k=1}^N\pi(\alpha_k^*;\calT)\geq\pi(\sigma_N;\calT)
 \geq\delta^2\prob_{r,\mathbf0}(H_N;\calT)>\delta^4.
\]
Since $|\calA| = L$, 
there is $a^*\in\calA$ such that $\sum_{k:\,a_k=a^*}\pi(\alpha_k^*;\calT)
    >\frac{\delta^4}{L}$.
Write $\{k:a_k=a^*\}=\{k_1<\cdots<k_m\}$ and set
$\alpha^*:=\alpha_{k_1}^*\dot+\cdots\dot+\alpha_{k_m}^*$.
The subsequence $(\alpha_{k_1}^*,\cdots,\alpha_{k_m}^*)$
is $\delta^2$-strongly orthogonal.
By \cite[Corollary~8.9(1)]{SS04},
applied with $\delta^2$ (instead of $\delta$),\begin{equation}\label{eq:short:abs}
 \pi(\alpha^*;\calT)\geq(1-3\delta^2)\sum_{q=1}^m \pi(\alpha_{k_q}^*;\calT)
 >\frac{\delta^4}{2L}=d.
\end{equation}

For every $q \in \{1,\dots,m\}$, the strategy profile $\alpha_{k_q}^*$ is an $\eta$-equilibrium of
the \emph{full} tree $\calT$, and has payoff in
$\prod_{i\in I}\bigl[a^{*,i},a^{*,i}+\delta\bigr]$.
Since the sequence is $\eta^2$-strongly orthogonal and
$\delta<\eta<1/6$,
\cite[Lemma~8.14]{SS04} (applied player by player) yields
\[
 \gamma^i(\alpha^*;\calT)\geq a^{*,i}-\eta,
 \qquad
 \hbox{and }
 \gamma^i(\alpha^i,\alpha^{*,-i};\calT)\leq a^{*,i}+8\eta,
\]
for every player $i \in I$ and every stationary strategy $\alpha^i$ of player~$i$.
Consequently, $\alpha^*$ is a stationary
$9\eta=54\sqrt\ep$-equilibrium.

\cite[Corollary~8.11]{SS04} applied
with $\delta^2$ (instead of $\delta$) gives
\begin{equation}\label{eq:short:lower}
 \gamma^i(\alpha^*;\calT)\geq a^{*,i}-6\delta^2
 >y^i_\calT-7\sqrt\ep.
\end{equation}
Thus, the antecedent of conclusion~(ii) cannot occur in this case.

Finally, suppose that
$\gamma^i(\alpha^*;\calT)\leq-7\sqrt\ep$.
For any $q=1,\dots,m$, \eqref{eq:short:lower} implies that
\[
 \gamma^i(\beta_{k_q}^*;\calT_{H_{k_q-1}})
 \leq a^{*,i}+\delta
 \leq\gamma^i(\alpha^*;\calT)+\delta+6\delta^2<0.
\]
By A2(iii), there is a node $s$ of
$\calT_{H_{k_q-1}}$ (which is also a node of $\calT$)
and $(\beta^j)_{j \neq i} \in [0,1]^{n-1}$, such that
\[
 \max\{A,B\}\leq\gamma^i(\beta_{k_q}^*;\calT_{H_{k_q-1}})
 \leq\gamma^i(\alpha^*;\calT)+\delta+6\delta^2
 <\gamma^i(\alpha^*;\calT)+\ep^2.
\]
For distinct players $i,j \in I$,
\[
 \pi(0^i,\alpha^{*,-i};\calT)+
 \pi(0^j,\alpha^{*,-j};\calT)\geq\pi(\alpha^*;\calT)\geq d>2d^2.
\]
Therefore, there can be 
at most one player who satisfies the two inequalities in the antecedent of (iii).
This completes the proof of the proposition.

%%%%%%%%%%%%%%%%%%%
\section{Proof of Theorem~\ref{theorem:three}}
\label{section:three}

The only place in the proof of Theorem~\ref{thm:main}
where the conditional atomlessness of the filtration takes part is in Section~\ref{sec:Acase}, when $R^{\mathrm{sp}}$ is $\ep$-close to a Q-matrix.
In the presence of three players, 
one can circumvent the use of the filtration.
Indeed, 
in this case,
\cite{Solan99} (see also \cite{Solan03})
extended the cyclic construction of \cite{FTV97}
and showed that the processes $(\widehat i_k)_{k \in \dN}$ and $(\widehat \delta_k)_{k \in \dN}$ that were used in Section~\ref{sec:Acase}
can be chosen to be independent of $(\psi_k)_{k \in \dN}$ ($\widehat i_k$ and $\widehat \delta_k$ only depend on $(\widehat i_0,\dotsm\widehat i_{k-1})$, for each $k \in \dN$).
In particular, in the presence of three players,
one does not need to use the public correlation,
and hence the conditional atomlessness of the filtration is not required.

%%%%%%%%%%%%%%%%%%%%
\section{Discussion}
\label{sec:discussion}

Our proof parallels that of \cite{SS04}
for the existence of an approximate equilibrium in two-player non-zero-sum stopping games.%
\footnote{\cite{SS04} do not present their proof in the terms below, but their proof can be presented that way.}
In both proofs, one identifies a stage $t_0$ and an $\calF_{t_0}$-measurable event $A$,
such that, on $A$, the game that starts at stage $t_0$,
admits a ``simple'' approximate equilibrium.
When there are two players,
this approximate equilibrium involves one player stopping with a small probability whenever her payoff upon stopping alone is maximal.
When there are many players, as in our case,
the event $A$ is characterized by the matrix $R$ being close to a Q-matrix; 
in that case, the simple approximate equilibrium has the players alternately stop with small probability.

On the complement event $\Omega \setminus A$,
both the proof of \cite{SS04} and ours construct an approximate equilibrium in a similar fashion.
For every two bounded stopping times $\tau < \tau'$ we approximate the stopping game induced between the stopping times $\tau$ and $\tau'$ (excluded) by a collection of stopping games on finite trees.
We then show that in this game there is an approximate equilibrium which terminates with a probability bounded away from $0$.
Finally, we use the stochastic variation of Ramsey's theorem
to piece together approximate equilibria in appropriately selected games on finite trees to obtain an approximate equilibrium in the original game.

In the multiplayer setup, 
the main difficulty in this part is to obtain the uniform lower bound on the probability of stopping in the games on finite trees.
To obtain this bound, it was important that on $\Omega \setminus A$,
the matrix $R^{\mathrm{sp}}$ is not close to a Q-matrix.
This is why we needed to consider approximate solutions to the linear complementarity problem.
This uniform bound also requires compactness of the family of linear complementarity problems under consideration. 
For that reason, we needed uniform bounds both on the matrices $R$ and on the vectors $q$.

\section*{Use of AI}
ChatGPT (OpenAI) was used to assist with reviewing mathematical
arguments and estimates, identifying potential gaps and errors,
checking references,
and improving the exposition.
The authors retain full responsibility for the content and
correctness of the paper.

%% =====================================================================
%% References
%% =====================================================================
\bibliographystyle{plainnat}

\end{document}